# A systematic literature review on the application of analytical approaches and mathematical programming in public bus transit network design and operations planning – Part I

**Reza Mahmoudi***

Department of Civil Engineering,

University of Calgary,

AB T2N 1N4, Calgary, Alberta, Canada

E-mail : reza.mahmoudi@ucalgary.ca

**Saeid Saidi**

Department of Civil Engineering,

University of Calgary,

AB T2N 1N4, Calgary, Alberta, Canada

E-mail : ssaidi@ucalgary.ca

**S. Chan Wirasinghe**

Department of Civil Engineering,

University of Calgary,

AB T2N 1N4, Calgary, Alberta, Canada

E-mail : wirasing@ucalgary.ca

*** Corresponding Author**

# Attention

Please note that this literature review encompasses four distinct papers, two of which are published in peer-reviewed journals and two are available on arXiv.org. To optimize your reading, please read the papers in the following order:

1) "A Systematic Literature Review on the Application of Analytical Approaches and Mathematical Programming in Public Bus Transit Network Design and Operations Planning – Part I" (published online on arXiv.org)
2) Mahmoudi, R., Saidi, S. and Wirasinghe, S.C., 2024. A critical review of analytical approaches in public bus transit network design and operations planning with focus on emerging technologies and sustainability. Journal of Public Transportation, 26, p.100100. (https://www.sciencedirect.com/science/article/pii/S1077291X24000201)
3) "A Systematic Literature Review on the Application of Analytical Approaches and Mathematical Programming in Public Bus Transit Network Design and Operations Planning – Part II" (published online on arXiv.org)
4) Mahmoudi, R., Saidi, S. and Emrouznejad, A., 2025. Mathematical Programming in Public Bus Transit Design and Operations: Emerging Technologies and Sustainability–A Review. Socio-Economic Planning Sciences, p.102155. (https://www.sciencedirect.com/science/article/pii/S0038012125000047)

***Note.** For any inquiries or contributions to the unpublished papers, please contact the corresponding author or other authors.*

***Corresponding Author:** Reza Mahmoudi, Postdoctoral Fellow, University of Toronto, reza.mahmoudi@ucalgary.ca, reza.mahmoudi@utoronto.ca, j.mahmoudi.reza@gmail.com*

# A systematic literature review on the application of analytical approaches and mathematical programming in public bus transit network design and operations planning – Part I

**Abstract**

The Public Bus Transit Network Design Problem and Operations Planning (PBTNDP&OP) remains a core research area within transportation, in particular, because of the emergence of new transit technologies and services. Analytical approaches and mathematical programming are the most commonly applied methodologies to study this problem. Many studies utilize either of these two methods, often viewed as competing due to the unique benefits each provides that the other does not. This two-part paper systematically reviews the application of analytical approaches and mathematical programming in PBTNDP&OP, analyzing publications from 1968 to 2021. It begins by comparing analytical methods and mathematical programming through various statistical analyses, including the number of published papers, the most active journals and authors, keyword frequencies, keyword co-occurrence maps, and co-authorship maps. Subsequent analysis of the identified papers includes examinations from multiple perspectives: the problems investigated, modeling methods used, decision variables considered, network structures, and key findings. This is followed by a critical review of selected papers. The paper concludes by discussing the advantages and disadvantages of each approach and suggests potential extensions for future research based on identified gaps in existing studies.

## 1. Introduction

With ever-increasing congestion, public transportation services (PTSs) have become the backbone of the transportation network of most societies. Flexibility, affordability, and sustainability are some of the advantages of PTSs compared to private transportation modes. Among all PTSs, undoubtedly, public bus transit service (PBTS) is the most prevalent PTS all around the world. PBTS, with its different versions (regular buses, BRT, feeder, shuttle, school bus, etc.), is applicable in any city with any transportation network size/structure and any population density. Because of this popularity and applicability, the public bus transit network design problem and operations planning (PBTNDP&OP) have attracted researchers' and transportation managers' attention.

From the mathematical perspective, PBTNDP&OP is assumed to be one of the most complicated transportation problems to analyze (Ibarra-Rojas et al., 2015). The multiplicity of the beneficiaries, decision-makers, and evaluation criteria, interdependencies among the solutions, uncertainty, and missing data are only some of the numerous reasons for the high complexity of PBTNDP&OP. Different methods have been applied in the literature to study these problems and find optimal or near optimal solutions. Analytical approaches (such as continuum approximations) and mathematical programming (optimization models) are two of the most popular mathematical methodologies in the literature of PBTNDP&OP. While some researchers have applied analytical approaches, others have employed the advantages of mathematical programming and heuristic/metaheuristic algorithms.

Despite the significant number of existing studies in this field focused on the different aspects of PBTNDPs&OP, there are many questions that remain to be answered. Because of dynamic technological developments, introducing emerging technologies to the transit networks, the ever-increasing importance of sustainable development and demand changes, even more issues will appear in the future. Hence, reviewing the literature related to the applications of the mentioned methodologies in PBTNDP&OP and comparing the advantages and disadvantages of each of them will provide important findings for future studies. In particular, this study aims to analyze the existing literature and review the definitions, classifications, objectives, constraints, network topology decision variables, modeling approaches, and key findings related to the PBTNDP&OP. In Part I of this literature review paper, the main focus will be on the applications of the analytical approaches in PBTNDP&OP. However, in Part II, a comprehensive literature review on the applications of mathematical programming in the PBTNDP&OP will be presented.

In Part I, the objectives are fourfold: first, classifying the existing studies in PBTNDP&OP based on the mentioned two major methodologies and comparing the advantages/disadvantages of each methodology; second, identifying recent trends in the analytical studies, existing research gaps and possible extensions for future research. In particular, a comparison between the applied approaches will be presented, which can help both researchers and decision-makers to select their modeling method before conducting a study; third, a comprehensive review and criticism of the selected papers; and finally, providing some statistical analysis on the identified published papers. The paper will provide some statistics about the "Number of published papers during the time", "target journals", "keyword frequency", "most cited papers", "most active authors", "map of co-authorship", and "CPM map".

The rest of the paper has been organized as follows: the applied methodology to find the published papers has been described in “Methodology”. “Statistical Analysis” presents some statistical analysis about the identified papers. A comprehensive review, including reviewing understudied problems, applied methods, network structures, and findings, has been presented in section 4. Existing research gaps and possible extensions have been discussed in section 5. “Analytical methods versus mathematical programming”, section 6, discusses the advantages and disadvantages of each methodology, presenting some comparisons between these methodologies. Finally, “Conclusion” includes the concluding remarks.

## 2. Methodology

A systematic search has been done in Scopus, Google Scholar, and Web of Knowledge in order to find a preliminary list of published papers. These are well-known databases for scientific articles. To this end, first, a search algorithm, as shown in Table 1, has been applied. In the next step, the preliminary list has been filtered to make the final list by removing unrelated papers. Only journal papers and conference papers published in English are considered. The PRISMA statement has been applied to conduct a systematic search, reviews, and analyses as shown in Figure 1 (Moher et al., 2009).

**Table 1.** The applied searching algorithm to find published papers.

| **Algorithm:** |
|---|
| #1. "Bus operation planning" [Article title, Abstract, Keywords] |
| #2. "Public transportation" [Article title, Abstract, Keywords] |
| #3. "Bus network design" [Article title, Abstract, Keywords] |
| #4. "Public Transit Network" [Article title, Abstract, Keywords] |
| #5. "Urban Mobility" [Article title, Abstract, Keywords] |
| #6. "Bus" [Article title, Abstract, Keywords] |
| #7. "Transit network design" [Article title, Abstract, Keywords] |
| #8. "Bus route design" [Article title, Abstract, Keywords] |
| #9. "Analytical method" [Article title, Abstract, Keywords] |
| #10. "Analytical approach" [Article title, Abstract, Keywords] |
| #11. "Analytical model" [Article title, Abstract, Keywords] |
| #12. "Continuous approximation" [Article title, Abstract, Keywords] |
| #13. "Continuum approximation" [Article title, Abstract, Keywords] |
| #14. "Discrete optimization" [Article title, Abstract, Keywords] |
| #15. "Bi-level optimization" [Article title, Abstract, Keywords] |
| #16. "Optimal design" [Article title, Abstract, Keywords] |
| #17. "Optimization" [Article title, Abstract, Keywords] |
| #18. "Algorithm" [Article title, Abstract, Keywords] |
| #19. "Heuristics" [Article title, Abstract, Keywords] |
| #20. "Metaheuristics" [Article title, Abstract, Keywords] |
| #21. "Source Type": [Journal] AND [Conference] AND [Thesis] |
| #22. "Document Type": [Article] AND [Article in Press] AND [Review] AND [Conference paper] AND [Thesis] |
| #23. "Article Language": [English] |
| #24. Review the article title, abstract and full text, identify and filter related and suitable papers. |
| #25. (#1 OR #2 OR #3 OR #4 OR #5 OR #6 OR #7 OR #8) AND (#9 OR #10 OR #11 OR#12 OR #13) AND #21 AND #22 AND #23 AND #24 |
| #26. (#1 OR #2 OR #3 OR #4 OR #5 OR #6 OR #7 OR #8) AND (#14 OR #15 OR #16 OR#17 OR #18 OR #19 OR #20) AND #21 AND #22 AND #23 AND #24 |
| #27. Put all identified papers in steps 25 in category 1 |
| #28. Put all identified papers in steps 26 in category 2 |

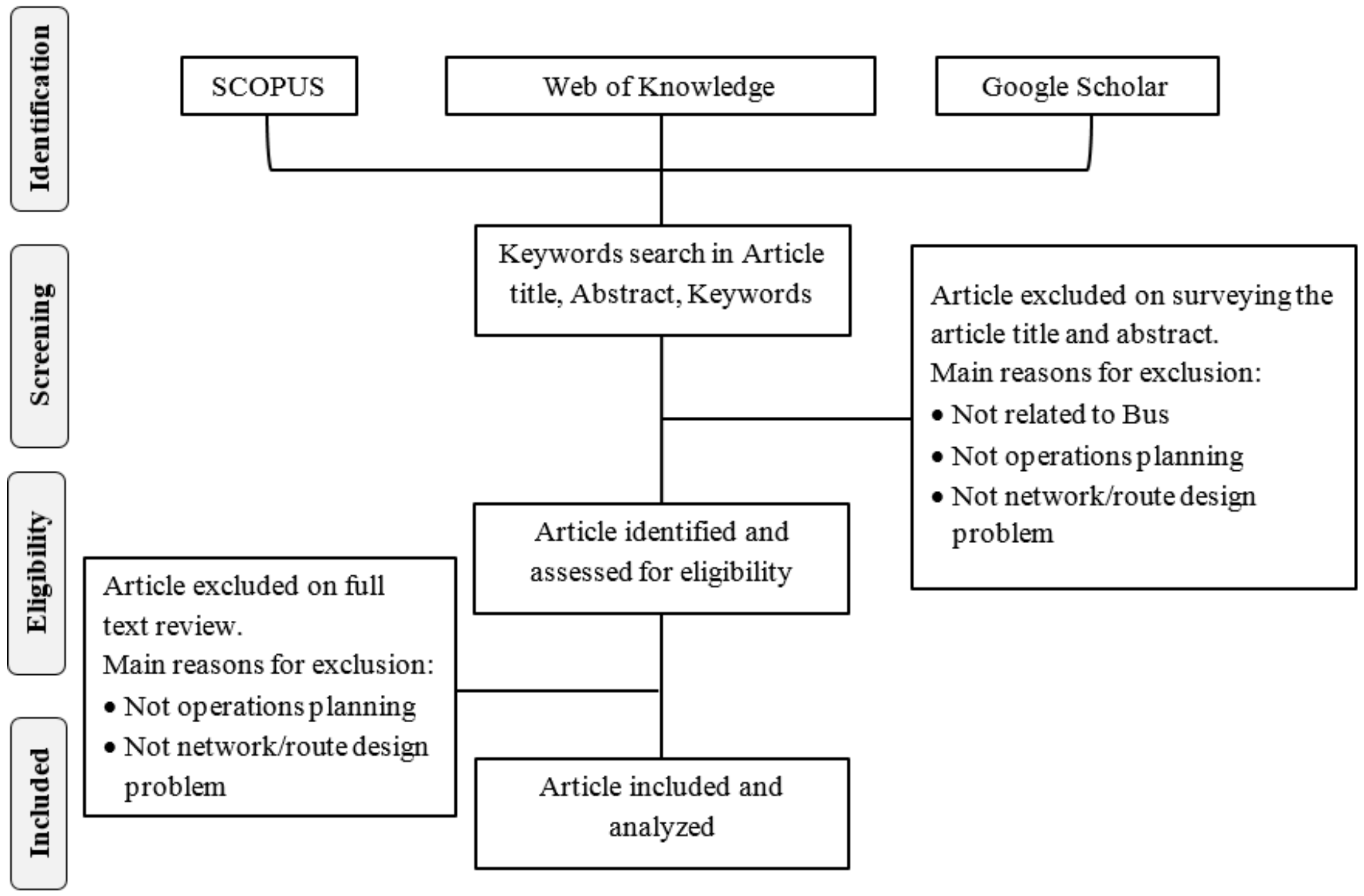

**Figure 1**. The PRISMA statement for systematic search, reviews and meta-analyses.

## 3. Statistical analysis

If a biologist wants to study a specific leaf on a specific tree in a specific forest, he/she won't go directly to see the leaf. He/she will start by studying the forest: the location of the forest, the size of the forest, the forest vegetation, the climate pattern of that specific area during the year, the history of the climate, vegetation changes in that area, etc. After all of these, the biologist will be able to conduct a study about the leaf.

The same is true for any study in other fields. In particular, if the objective is to review literature related to a specific topic, it would be very helpful and beneficial to provide general discussions about the existing literature before going through the details; to show what the forest looks like before discussing the leaf. The objective of this section is to provide some general statistics and trend analysis about the applications of analytical analysis and mathematical programming and algorithms in PBTNDP&OP. Before reviewing the identified papers, it would be helpful to see what statistics reveal about these papers. The statistical analysis can help academia in order to enhance their academic network, find the most active and interesting sources in PBTNDP&OP, find the most recent hot spots in these fields, help junior researchers to find positions for their graduate studies, etc.

Following the proposed algorithm in Table 1 and the PRISMA statement presented in Figure 1, 285 papers are identified for category 1 (Analytical approaches) and 826 papers are found for the second category (mathematical programming and algorithms). Figure 2 shows the number of published papers each year. It should be noted that the papers chosen were either directly related to bus transit services or to general public transportation services, which included bus transit services. Analyzing this figure can lead to seeing interesting facts. First, the total number of published papers in each year has steadily increased over time; additionally, it is worth noting that in only the year 2021, 269 papers related to both categories were published, which is more than the total number of published papers in this area until 2011. This indicates an increasing interest in the applications of these methodologies in this area and shows that it is now (and will likely be in the future) a hot spot in the area of transportation. This observation is not only true for the general topic but also for both the analytical and optimization categories. Second, the results show that until 2001, the number of published papers related to the first category (both in each year and in total) was almost always higher than the number of published papers in the second category. But after that, the number of publications related to mathematical programming is significantly larger than the one related to analytical approaches. This observation can be easily explained by linking it to technological development. Although applying mathematical programming allows the researchers to consider more details, objectives, constraints, and realistic assumptions in their study, most of the developed models in this category were/are highly complicated and impossible to provide close form describable solutions. They are also time-consuming and computationally intensive. Hence, this method was not a popular method in the 20$^{th}$ century. However, by the increasing availability of advanced computers in the late 1990s and early 21st century, most researchers started to apply mathematical programming and optimization algorithms in PBTNDPs&OP. By the end of 2021, the total number of published papers in the second category is almost three times the number of publications in the first category. The main advantages and disadvantages of each method will be discussed later.

It is always important for authors and researchers to find which journals are more interested in their research area and are more likely to publish their studies.

Table 2 presents the most active journals publishing on these topics. Based on this table, “Transportation Research Part B Methodological” is the journal with the greatest interest in this field (in both categories).

Other interesting findings are the authors and countries with the most published articles on these topics. This is especially true for individuals seeking academic positions or research funding, as well as well-known PBTNDP&OP researcher collaborators. In addition, the recently published papers by the pioneers of this area will be a reliable source to see what the latest updates and ongoing hot topics are. According to Table 3, China, the United States, Hong Kong, and Canada are countries with more than 50 published papers, making them attractive destinations for researchers who want to pursue their graduate studies in this field.

Table 4 shows that in the area of analytical approaches, Chan Wirasinghe, with 18 published, papers is one of the leaders in this area. He is being followed by Paul Schonfeld (17 papers), Michael Cassidy (12 papers), and GU Weihua/Wenbo Fan (10 papers). In the second category, Avi Ceder, with 22 papers, is the most active author, followed by Paul Schonfeld (17 papers), Konstantinos Gkiotsalitis/Steven I-Jy Chien/Munoz, J.C. (13 papers), and Zhiyuan Terry Liu (11 papers). Overall, Paul Schonfeld, with 29 published papers, is the leader in the area of PBTNDP and operations planning.

Figure 3 shows the co-authorship network between the authors with more than two published papers in the first category, and Figure 4 presents this network for the second category and between the authors who have published at least four papers. All co-authorship and co-occurrence networks have been drawn using VOSviewer software based on the data base of identified papers. There are some interesting observations in these figures. For example, based on Figure 3, there are nine different co-authorship networks that have focused on the applications of analytical methods to PBTNDP&OP. According to Figure 4, this number for mathematical programming is eleven. The size of the circles is determined by the number of published papers by each author and his or her co-authors. Hence, looking at this figure, the leaders and most active authors in each network can be seen. Note that it is possible to not see some well-known authors' names in the co-authorship networks. The reason for this is that, despite publishing a number of studies, the author does not meet the co-authorship criteria in order to be considered in the network. Looking at these figures, one can find the most active labs in each area and their recent focus. If someone is enthusiastic about the research interests of a pioneer researcher, he/she can consider recent publications of the co-authors or lab members of that researcher.

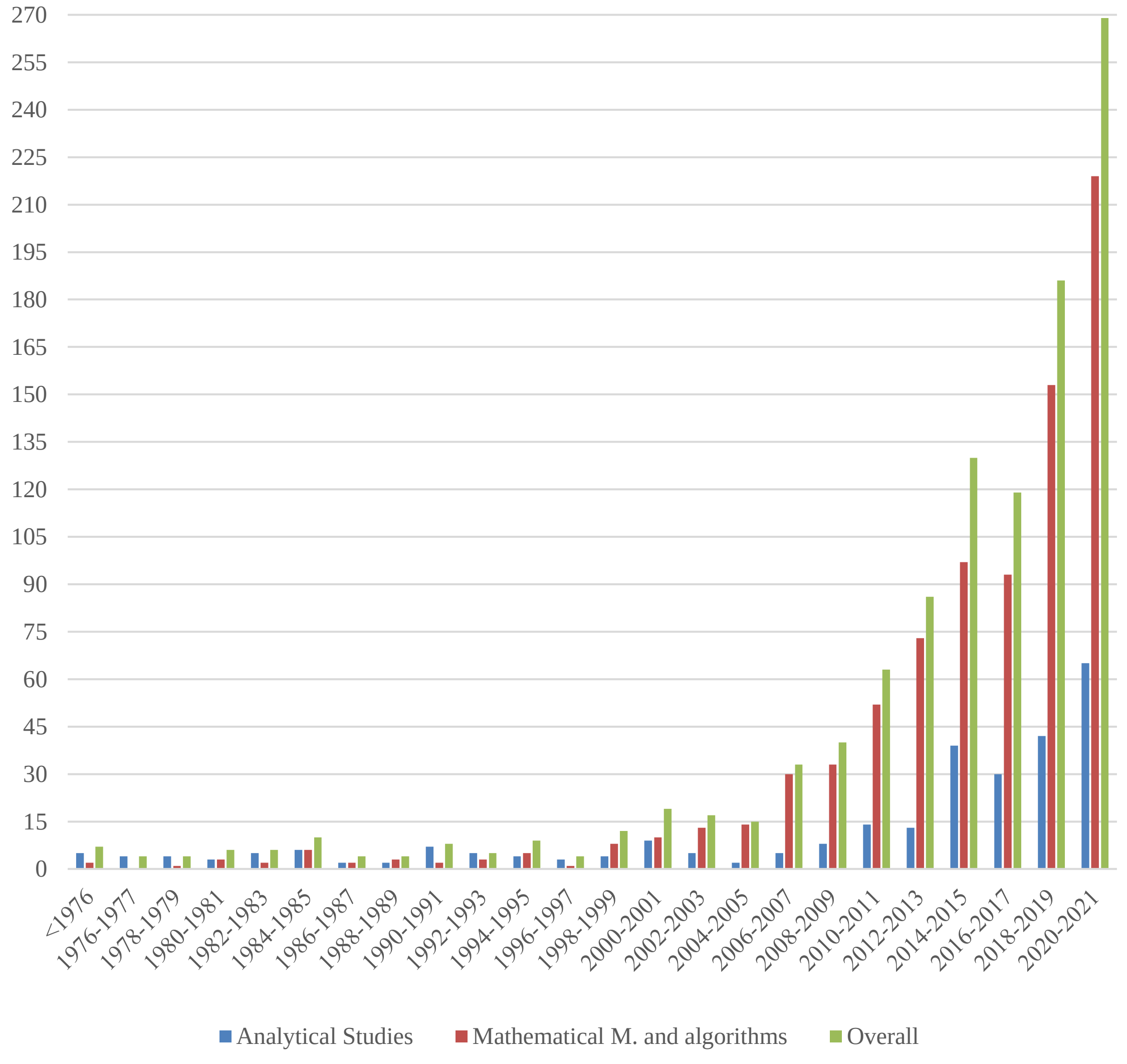

**Figure 2.** Annual number of published articles.

**Table 2.** Most active journals in each category.

| Category | Journal | N. of published papers |
|---|---|---|
| **Analytical models** | Transportation Research Part B Methodological | 39 |
| | Journal of Advanced Transportation | 23 |
| | Transportation Research Part A Policy and Practice | 16 |
| | Transportation Research Record | 17 |
| | Transportation Research Part C Emerging Technologies | 15 |
| | Transportation Science | 12 |
| | Journal of Transportation Engineering | 10 |
| | Transportation Research Part E Logistics and Transportation Review | 8 |
| | Transportation Research Procedia | 7 |
| | IEEE Transactions on Intelligent Transportation Systems | 7 |
| **Mathematical programming and Algorithms** | Transportation Research Part B Methodological | 59 |
| | Transportation Research Part C Emerging Technologies | 52 |
| | Transportation Research Record | 43 |
| | Journal of Advanced Transportation | 37 |
| | Transportation Research Part A Policy and Practice | 26 |
| | Journal of Transportation Engineering | 24 |
| | IEEE Transactions on Intelligent Transportation Systems | 21 |
| | Public Transport | 21 |
| | European Journal of Operational Research | 16 |
| | Transportation Planning and Technology | 15 |
| **Overall** | Transportation Research Part B Methodological | 90 |
| | Transportation Research Part C Emerging Technologies | 64 |
| | Transportation Research Record | 57 |
| | Journal of Advanced Transportation | 52 |
| | Transportation Research Part A Policy and Practice | 37 |
| | Journal of Transportation Engineering | 31 |
| | IEEE Transactions on Intelligent Transportation Systems | 27 |
| | Public Transport | 23 |
| | Transportation Research Part E Logistics and Transportation Review | 20 |
| | Transportation Science | 20 |

**Table 3**. Top countries that have published most papers in each category.

| Analytical models | | Mathematical M. and Algorithms | | Overall | |
|---|---|---|---|---|---|
| **Country** | N. of published paper | Country | N. of published paper | Country | N. of published paper |
| United States | 103 | China | 278 | China | 321 |
| China | 58 | United States | 184 | United States | 264 |
| Canada | 24 | Hong Kong | 49 | Hong Kong | 59 |
| United Kingdom | 18 | Chile | 36 | Canada | 50 |
| Chile | 17 | Italy | 32 | Chile | 49 |
| Australia | 15 | Iran | 31 | Australia | 44 |
| Hong Kong | 15 | Australia | 30 | Taiwan | 39 |
| Taiwan | 14 | Spain | 30 | United Kingdom | 38 |
| Spain | 11 | Singapore | 30 | Spain | 36 |
| Netherlands | 9 | Taiwan | 27 | Iran | 34 |

**Table 4.** Most active authors in first category, second category and in overall.

**Analytical models**

| Author | N. of published papers |
|---|---|
| Wirasinghe, S.C. | 18 |
| Schonfeld, P. | 17 |
| Cassidy, M.J. | 12 |
| Fan, W.<br>Gu, W. | 10 |
| Daganzo, C.F.<br>Gschwender, A.<br>Jara-Díaz, S. | 9 |
| Badia, H.<br>Estrada, M.<br>Quadrifoglio, L.<br>Chien, S.I.J. | 7 |
| Robusté, F. | 6 |
| Dakic, I.<br>Dessouky, M.<br>Fielbaum, A.<br>Madanat, S.<br>Menendez, M.<br>Nie, Y.<br>Luo, S.<br>Newell, G.F. | 5 |

**Mathematical M. and Algorithms**

| Author | N. of published papers |
|---|---|
| Ceder, A. | 22 |
| Schonfeld, P. | 17 |
| Gkiotsalitis, K.<br>Chien, S.I.J.<br>Munoz, J.C. | 13 |
| Liu, Z. | 11 |
| Laporte, G. | 10 |
| Ibeas, A.<br>Lo, H.K.<br>Mahmassani, H.S.<br>Meng, Q.<br>Petrelli, M.<br>Qu, X.<br>Song, R.<br>Szeto, W.Y.<br>Wei, M.<br>Yan, Y.<br>Jara-Díaz, S. | 8 |
| Estrada, M.<br>Gao, Z.<br>Gschwender, A.<br>Guan, W.<br>Kepaptsoglou, K.<br>Liu, T.<br>Mesa, J.A.<br>Sun, B.<br>Tang, J.<br>Yu, J. | 7 |

**Overall**

| Author | N. of published papers |
|---|---|
| Schonfeld, P. | 29 |
| Ceder, A. | 24 |
| Chien, S.I.J.<br>Wirasinghe, S.C. | 19 |
| Muñoz, J.C. | 16 |
| Cassidy, M.J.<br>Gkiotsalitis, K. | 13 |
| Daganzo, C.F.<br>Liu, Z. | 12 |
| Estrada, M.<br>Fan, W.<br>Gschwender, A.<br>Gu, W.<br>Jara-Díaz, S. | 11 |
| Laporte, G. | 10 |
| Liu, T.<br>Lo, H.K. | 9 |
| Badia, H.<br>Cats, O.<br>Guan, W.<br>Ibeas, A.<br>Madanat, S.<br>Mahmassani, H.S.<br>Meng, Q.<br>Ouyang, Y.<br>Petrelli, M. | 8 |
| | Qu, X.<br>Quadrifoglio, L.<br>Robusté, F.<br>Song, R.<br>Szeto, W.Y.<br>Wei, M.<br>Yan, Y. |

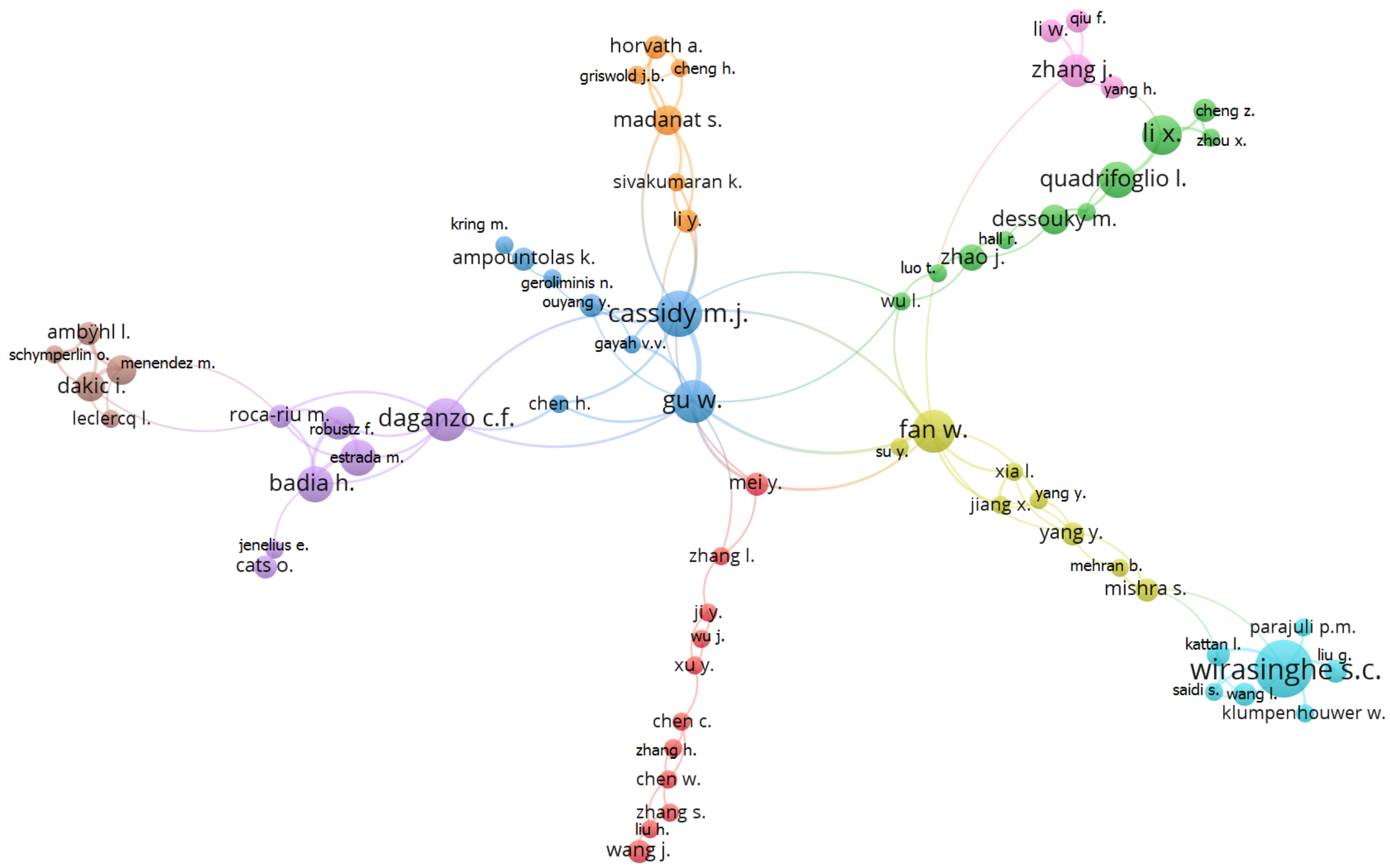

**Figure 3.** Map of co-authorship network related to published papers in the analytical studies (with at least 2 papers).

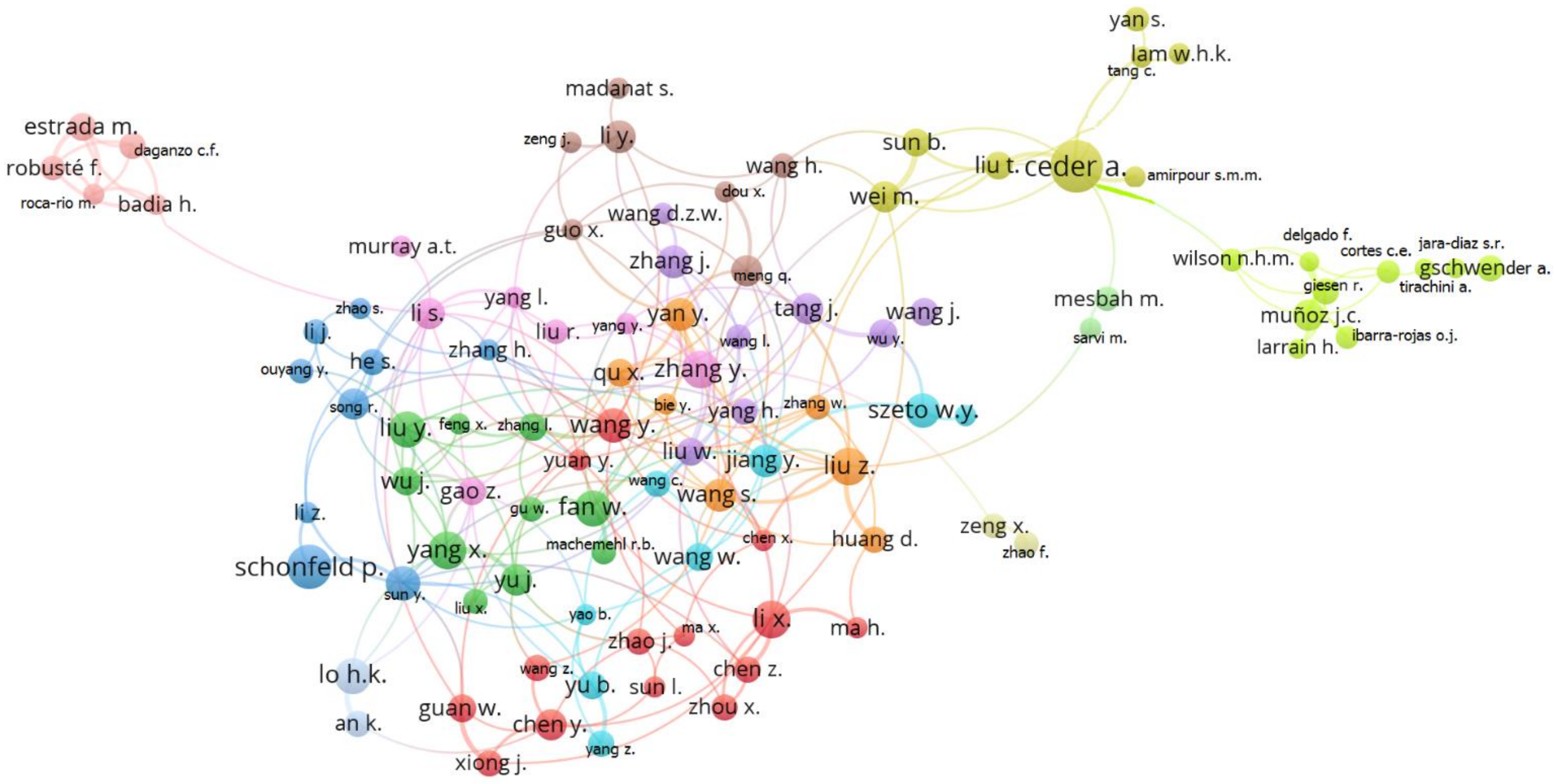

**Figure 4.** Map of co-authorship network related to published papers in the mathematical programming (with at least 4 papers).

The most cited papers in any area are one of the best references for researchers who are interested in that area. These papers are mostly the pioneers in their class of the problem. Also, analyzing the citations of these papers will provide a clear insight into the ongoing trend in the related class of the problem. Especially, for beginners, it would be very helpful to read the main references in their area of interest. Table 5 shows the most cited papers in each category. According to this table, "Some issues relating to the optimal design of bus routes", written by Gordon F. Newell, one of the most famous scientists in transportation studies, is the most cited paper in the first category. For the second category, the most cited paper is "Optimal strategies: a new assignment model for transit networks", written by Spiess and Florian. This paper has been cited more than 1000 times, which indicates that this study is a pioneer in this category.

**Table 5.** Most cited papers in each category (Based on Google Scholar)

| Analytical models | | Optimization and Algorithms | | Overall | |
|---|---|---|---|---|---|
| **Article** | **Citations** | **Article** | **Citations** | **Article** | **Citations** |
| Newell (1979) | 330 | Spiess and Florian (1989) | 1043 | Spiess and Florian (1989) | 1043 |
| Daganzo and Pilachowski (2011) | 291 | Ceder and Wilson (1986) | 781 | Ceder and Wilson (1986) | 781 |
| Newell (1971) | 273 | Pattnaik et al. (1998) | 451 | Pattnaik et al. (1998) | 451 |
| Hickman (2001) | 268 | Ceder et al. (2001) | 393 | Ceder et al. (2001) | 393 |
| Xuan et al. (2011) | 245 | Mandl (1980) | 359 | Bertolini (1999) | 389 |
| Chang and Schonfeld (1991a) | 226 | Baaj and Mahmassani (1995) | 347 | Mandl (1980) | 359 |
| Wirasinghe and Ghoneim (1981) | 214 | Fan and Machemehl (2006a) | 340 | Baaj and Mahmassani (1995) | 347 |
| Kuah and Perl (1988) | 203 | Baaj and Mahmassani (1991) | 335 | Fan and Machemehl (2006a) | 340 |
| Daganzo (1978) | 185 | Eberlein et al. (2001) | 311 | Baaj and Mahmassani (1991) | 335 |
| Jara-Díaz and Gschwender (2003) | 180 | Daganzo (2010) | 289 | Newell (1979) | 330 |
| Vuchic and Newell (1968) | 176 | Bielli et al. (2002) | 273 | Eberlein et al. (2001) | 311 |
| Zhao et al. (2006) | 173 | Xuan et al. (2011) | 245 | Daganzo and Pilachowski (2011) | 291 |
| Chien and Schonfeld (1998) | 173 | Chakroborty (2003) | 245 | Daganzo (2010) | 289 |
| Dubois et al. (1979) | 170 | Liu et al. (2013) | 233 | Bielli et al. (2002) | 273 |
| Li and Quadrifoglio (2010) | 166 | Delgado et al. (2012) | 232 | Newell (1971) | 273 |
| de Palma and Lindsey (2001) | 166 | Tom and Mohan (2003) | 224 | Hickman (2001) | 268 |
| Quadrifoglio and Li (2009) | 160 | Szeto and Wu (2011) | 223 | Xuan et al. (2011) | 245 |
| | 159 | Zhao and Zeng (2008) | 222 | Chakroborty (2003) | 245 |
| Wu and Hounsell (1998) | | Ibarra-Rojas and Rios-Solis (2012) | 221 | Liu et al. (2013) | 233 |
| Geroliminis et al. (2014) | | Fan and Machemehl (2006b) | 219 | Delgado et al. (2012) | 232 |
| Newell (1974) | | | | | |

PBTNDP&OP is a vast field for research, and different types of problems have been studied in this field. It is important to know what these problems/sub-categories are and which problems have attracted more attention. To this end, Table 6, Figure 5, and Figure 6 show the frequency of the keywords across the literature in each category, the map of the co-occurrence of the keywords in the first category (with at least 3 co-occurrences), and the map of the co-occurrence of the keywords in the second category (with at least 5 co-occurrences), respectively. Looking at these figures, for each category, the major problems related to the PBTNDP&OP that have been investigated in the literature can be identified to some extent. To this end, different clusters of the keywords presented in Figure 5 and Figure 6 must be analyzed. For example, based on Figure 5, "Adaptive control", "Bus holding", "Bus bunching", "Transit operation", and "Dynamic holding" have co-occurred in the literature; hence, bus bunching problems are one of the major subjects investigated in the literature. We can see that in the first category, transit network design, crowding and congestion management, dispatching/scheduling policy, reliability and travel time, feeder transit services, and fixed/flexible transit services are some of the hot sub-categories among the researchers. Most of the published papers will fall into one of these sub-categories. The same sub-categories can also be seen for the second category. In addition, it seems dynamic and bi-level programming for PBTNDP&OP, routing, rapid transit planning, and exclusive lane assigning are other popular sub-problems in this category.

**Table 6.** Frequency of the keywords across the literature in each category.

| **Analytical models** | | **Optimization and Algorithms** | | **Overall** | |
|---|---|---|---|---|---|
| Keywords | Frequency | Keywords | Frequency | Keywords | Frequency |
| Public Transport | 81 | Optimization | 321 | Optimization | 364 |
| Buses | 71 | Buses | 286 | Buses | 342 |
| Optimization | 60 | Bus Transportation | 217 | Bus Transportation | 258 |
| Bus Transport | 59 | Mass Transportation | 162 | Public Transport | 224 |
| Bus Transportation | 55 | Genetic Algorithms | 158 | Mass Transportation | 201 |
| Transportation | 44 | Public Transport | 158 | Bus Transport | 201 |
| Mass Transportation | 43 | Bus Transport | 150 | Travel Time | 175 |
| Transportation Planning | 33 | Travel Time | 146 | Transportation | 169 |
| Travel Time | 33 | Transportation | 131 | Genetic Algorithms | 161 |
| Urban Transport | 31 | Transportation Routes | 126 | Urban Transportation | 142 |
| Numerical Model | 24 | Urban Transportation | 122 | Transportation Routes | 138 |
| Urban Transportation | 23 | Integer Programming | 100 | Transportation Planning | 124 |
| Costs | 23 | Scheduling | 100 | Public Transportation | 108 |
| Public Transportation | 21 | Transportation Planning | 99 | Scheduling | 105 |
| Transportation Routes | 21 | Genetic Algorithm | 91 | Integer Programming | 100 |

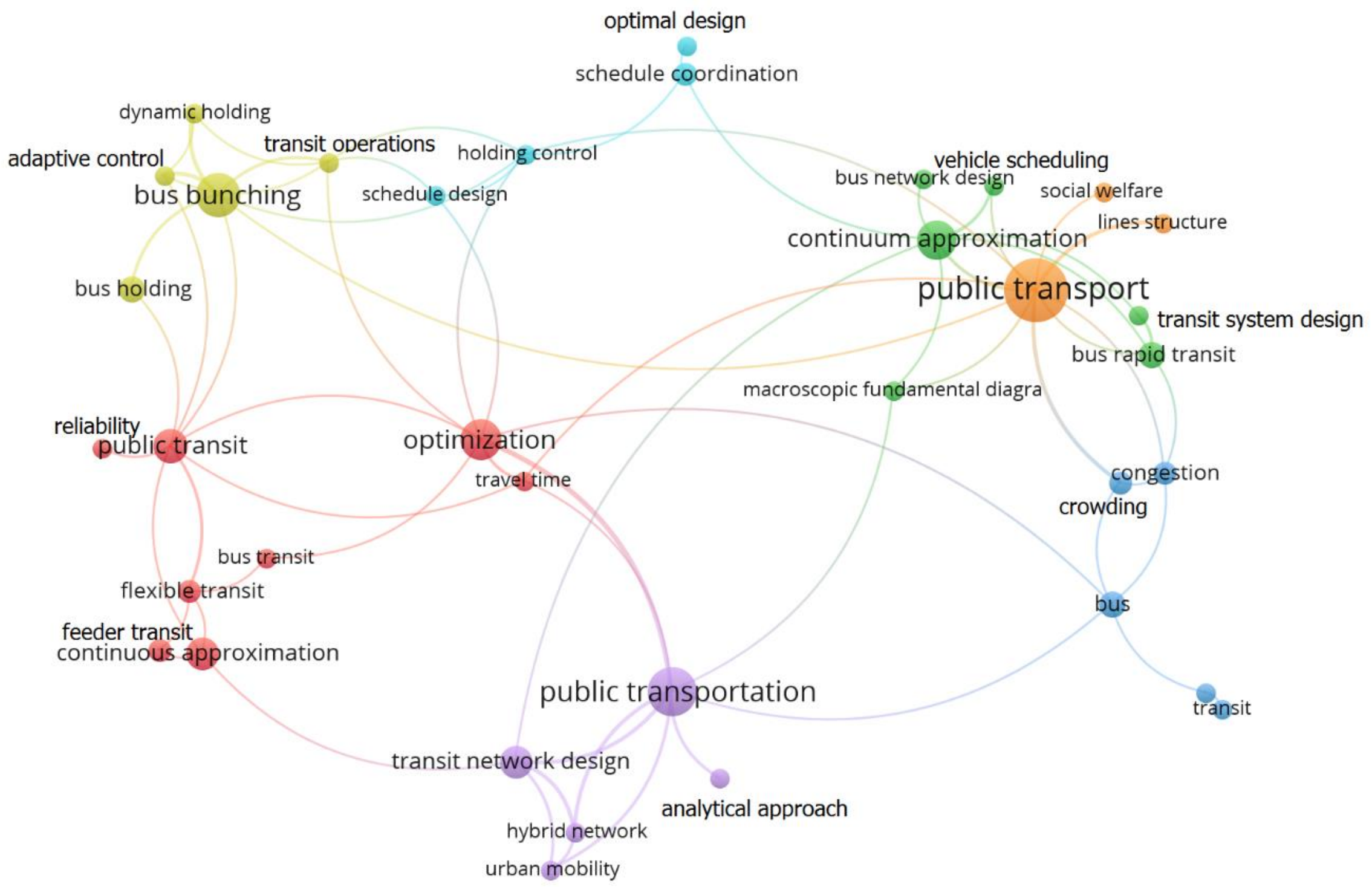

**Figure 5.** Map of most co-occurrence keywords related to published papers in the analytical studies (with at least 3 co-occurrences).

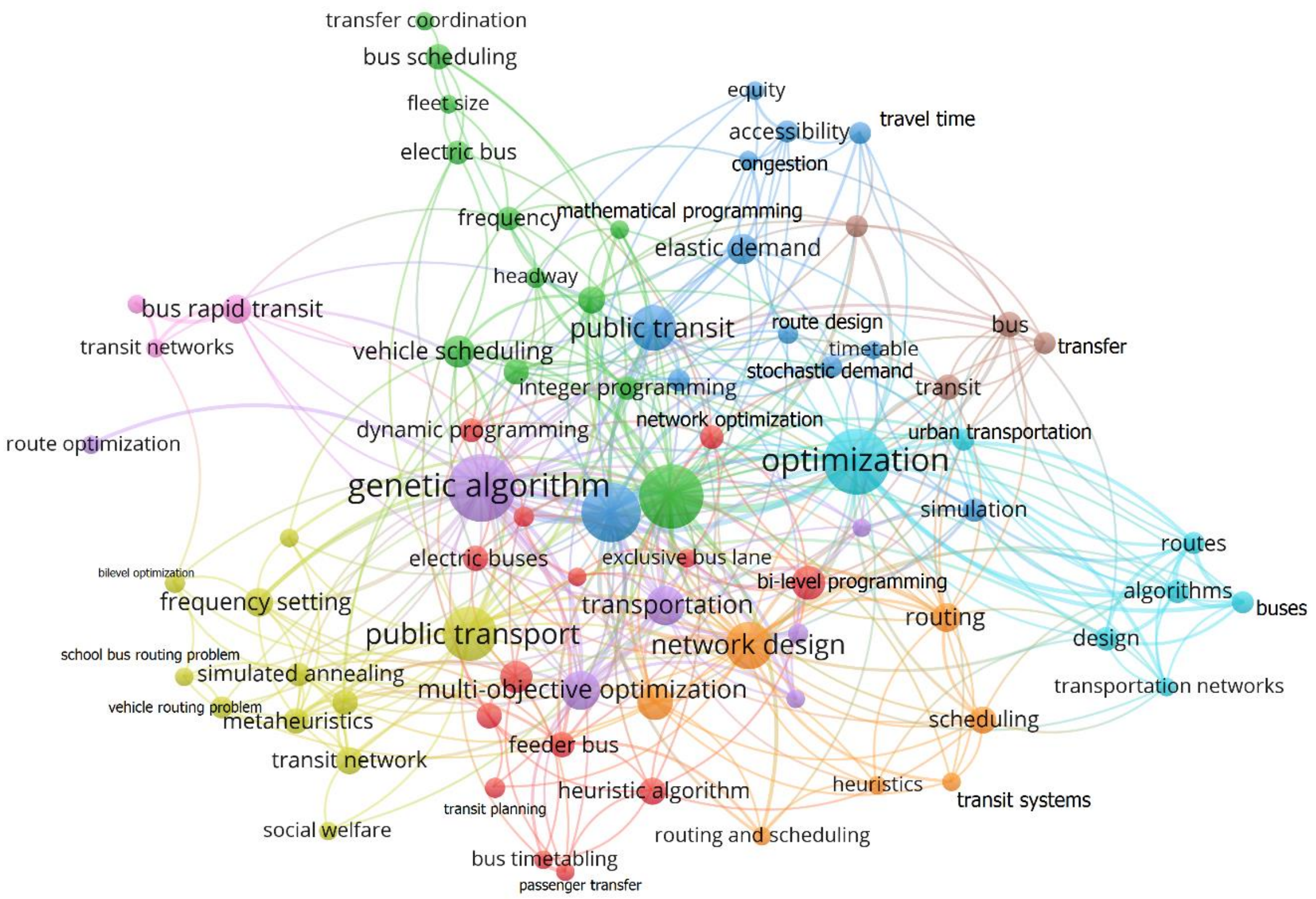

**Figure 6.** Map of most co-occurrence keywords related to published papers in mathematical programming and algorithms
(with at least 5 co-occurrences).

When it comes to reviewing literature statistically, the Critical Path Method (CPM) will be the most interesting, efficient, and helpful method in order to visualize the literature.

CPM visualizes the literature on a diagram which is understandable and trackable for everyone. In particular, a CPM shows (Kaffash and Marra, 2017; Kejžar et al., 2010):

- The origin and the destination of a topic and its subtopics
- The development paths form the origins to the destinations within a topic and its subtopics
- The oldest published paper in a topic and a sub-topic
- How the different sub-topics have been developed during the time,
- Which topics are still attractive for researchers,
- Which topics are getting old with fewer interests in conducting new studies
- When a new sub-topic has appeared in the literature, and how its development is going
- What the recent hot spots are
- Which papers are the pioneer papers in each year

Although CPM can present all of this information just with a figure, it would be very challenging to prepare an accurate CPM as it demands lots of time to conduct a precise search in the literature and find major sub-categories in the development path of a field. Some software such as Pajek can be used to draw this diagram, or it can be drawn manually. Some criteria must be considered in order to draw a CPM for a topic. First, the lowest level in the CPM for a topic must show the oldest identified paper(s) in that topic, and the highest level must show the identified papers in the study period's final year. Second, in each level, only published papers in a specific year must be shown. For example, all published papers at level 5 must be published in the same year. Third, only pioneer papers from each year must be considered. Different measures can be considered to choose a paper as a pioneering paper, such as the number of citations compared to the other papers published in the same year, authors' reputation, the reputation of the journal publishing the paper, and/or the novelty of the subject. However, the main criterion is the number of citations. In the case of having two or more papers with almost the same citation numbers, other criteria can be considered. Fourth, each paper must be connected to at least one other paper (ideally only one paper) at the latest level, if it is not the first paper in an independent sob-category. If a paper has cited another pioneering paper in the latest level and/or has contributed to that sub-topic as a pioneering paper in that level and in that sub-category, that paper must be connected to the latest paper in its sub-category with an arrow. Fifth, in the CPM, each branch is showing a sub-category. If at a level, a new sub-topic has been introduced to the literature by a paper, a new branch must be added to the CPM showing this fact. If the newly introduced sub-topic is somehow related to one of the existing sub-categories, the origin of the new branch can be connected to that existing branch; otherwise, an independent branch must be added to the CPM. Sixth, although it has been suggested to just report one paper for the lower levels, for the top levels (recent years), a group of papers for each level can be reported. To read more about the CPM in reviewing literature see Kejžar et al. (2010), Emrouznejad and Marra (2014), and Maltseva and Batagelj (2019).

Figure 7 and Figure 8 show the CPM for the first and second categories, respectively. According to Figure 7, Vuchic and Newell (1968) was the oldest paper among the identified papers for the first category, the applications of analytical approaches in PBTNDP&OP. This paper was on intersection spacing for rapid transit systems in order to minimize travel time. However, many different topics have emerged later. For analytical studies, 13 main sub-categories were identified from Figure 7, namely: 1) Macroscopic Fundamental Diagram and traffic assignment, 2) Routing, 3) Feeder transit planning, 4) Dispatching policy, 5) Transit network design, 6) Bus holding strategies, 7) Fleet size and capacity of the vehicle, 8) Travel time reliability, 9) Bus priority, 10) Rapid transit and shuttle transit system planning, 11) Flexible and fixed transit services, 12) Emerging technologies in public bus transit systems, 13) Sustainability. This figure also shows which study was the first one in its channel. For example, Wirasinghe et al. (1977) was one of the first papers that considered the feeder bus planning problem. Also, Figure 7 shows the frequency of publishing papers in each channel during the time. For example, Wu and Hounsell (1998) applied analytical approaches to study the bus priority problems; however, as it can be seen, there wasn't any study in this sub-category until 2015 and another one in 2020. The literature, existing research gaps, and possible extensions related to each channel will be discussed in more detail in sections 4.

The CPM diagram has been also developed for the second category shown in Figure 8. Based on this diagram, 15 major subcategories have been identified for PBTNDP&OP and mathematical modeling and algorithms, including 1) Accessibility and coverage, 2) Feeder transit network design, 3) Service pricing and fare management, 4) Fleet management, 5) Dispatching policy, 6) Stop spacing, 7) Mode split and traffic assignment, 8) Network design and Dispatching policy, 9) Network design and routing, 10) Green and sustainable public bus network design, 11) Flexible and fixed transit service, 12) Travel time and Reliability, 13) Disaster and disruption management, 14) Bus priority, and 15) Emerging technologies and electric buses. For example, it can be seen that Newton and Thomas (1969) is the oldest identified paper in this category. Also, subcategories 5, 8 and 9 have been the most popular categories among researchers from the beginning. However, subcategory 15 is getting extremely popular in recent years.

Analyzing the CPMs of two categories can reveal interesting facts. For example, number of studies that have applied mathematical programming to study emerging technologies in PBTSs and sustainability is much more than the number of analytical studies on these problems. We can see that there are four subcategories that include a good number of published studies in mathematical programming, but have not been very popular among analytical studies. These subcategories are “Accessibility and coverage”, “Service pricing and fare management”, “Flexible and fixed transit service”, and “Disaster and disruption management”. In particular, there is almost no study that have applied analytical approaches to investigate Disaster and disruption management in PBTSs.

**Note.** *The literature, existing research gaps, and possible extensions related to each sub-category in the applications of mathematical programming in PBTNDP&OP will be discussed in Part II of this literature review paper.*

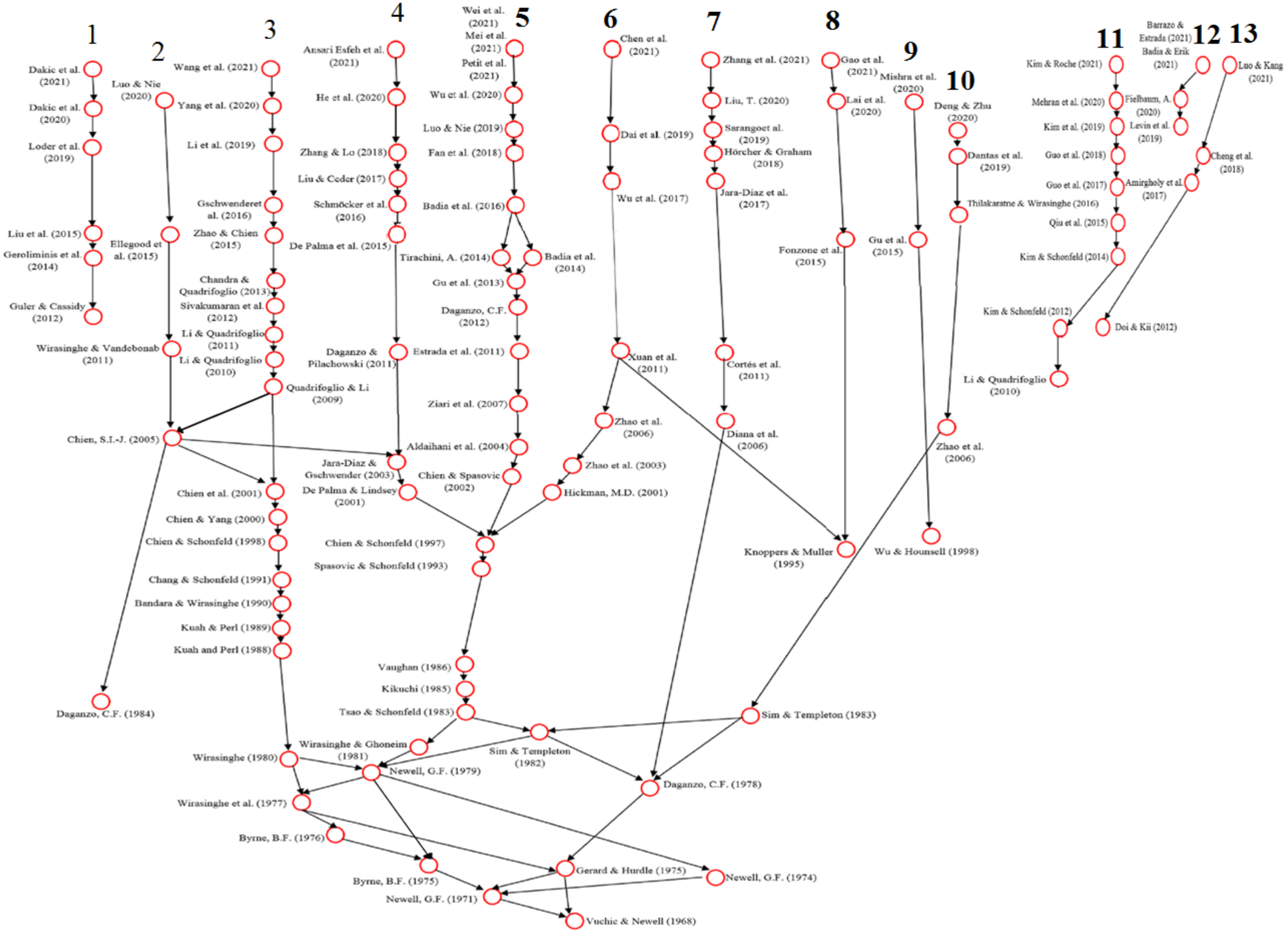

**Figure 7.** CPM of development analytical models applications in PBTNDP&OPs.

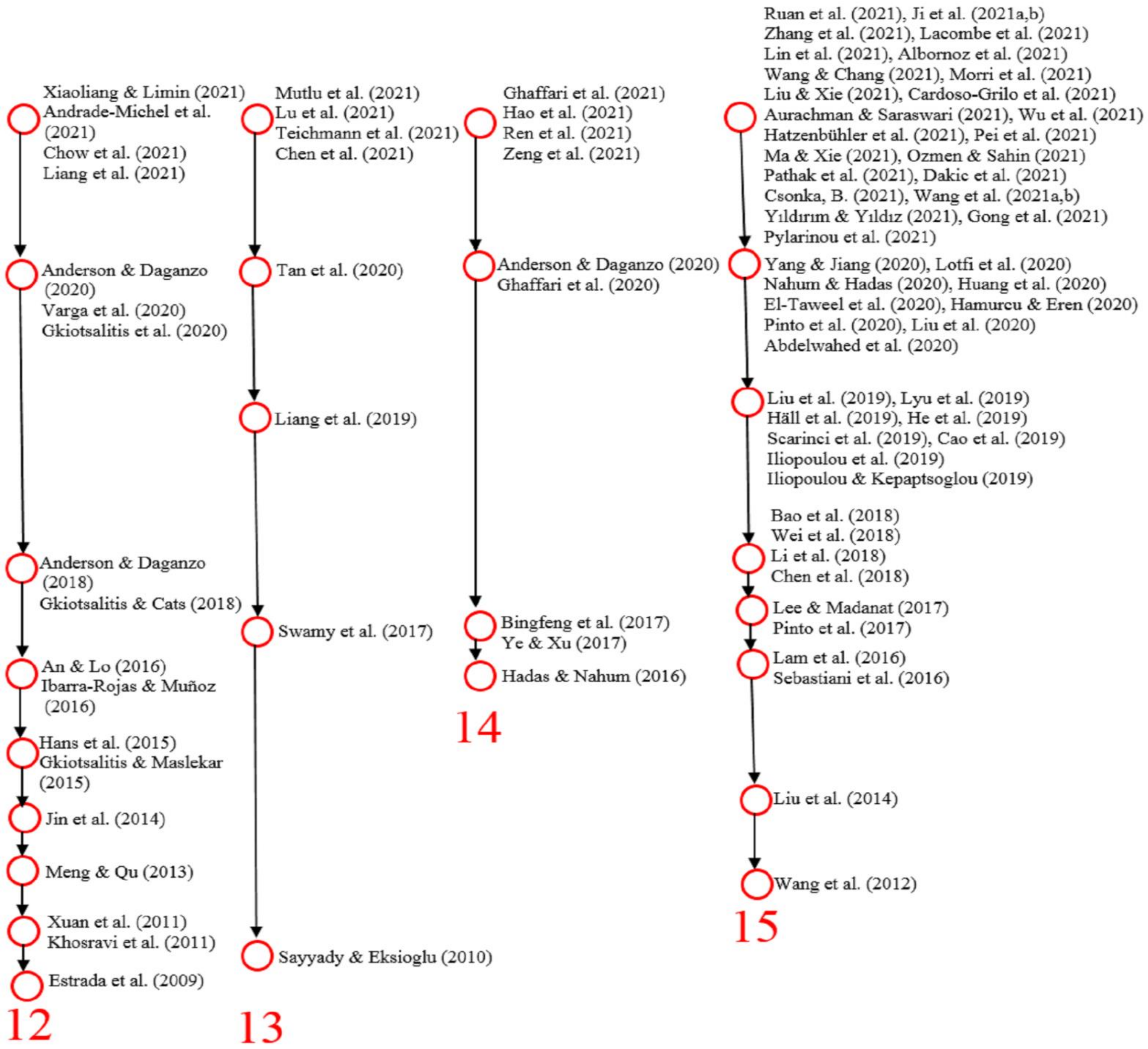

(B) emerging subcategories from 2009

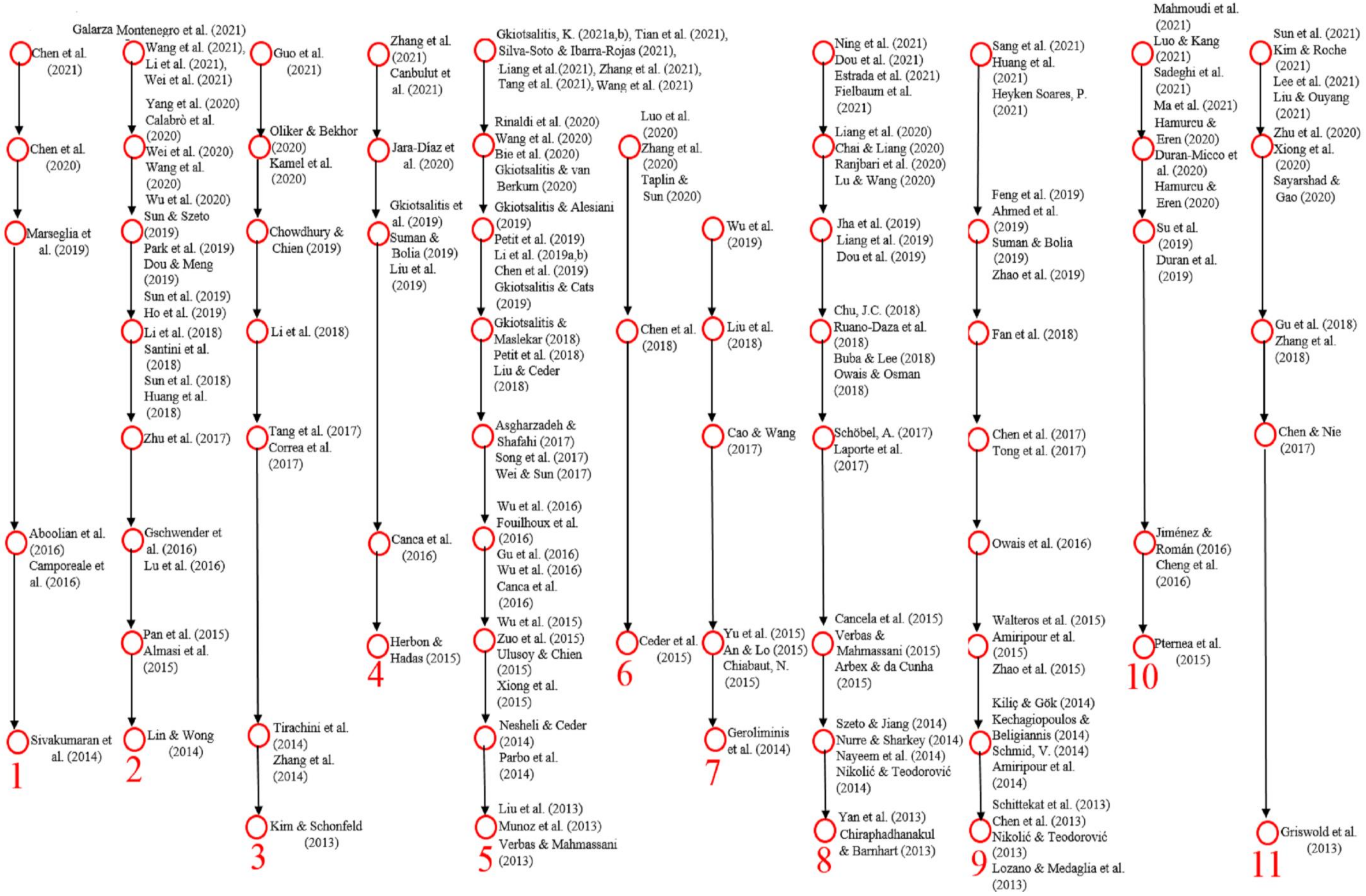

(A) continued: subcategories 1-11 from 2013 to 2021

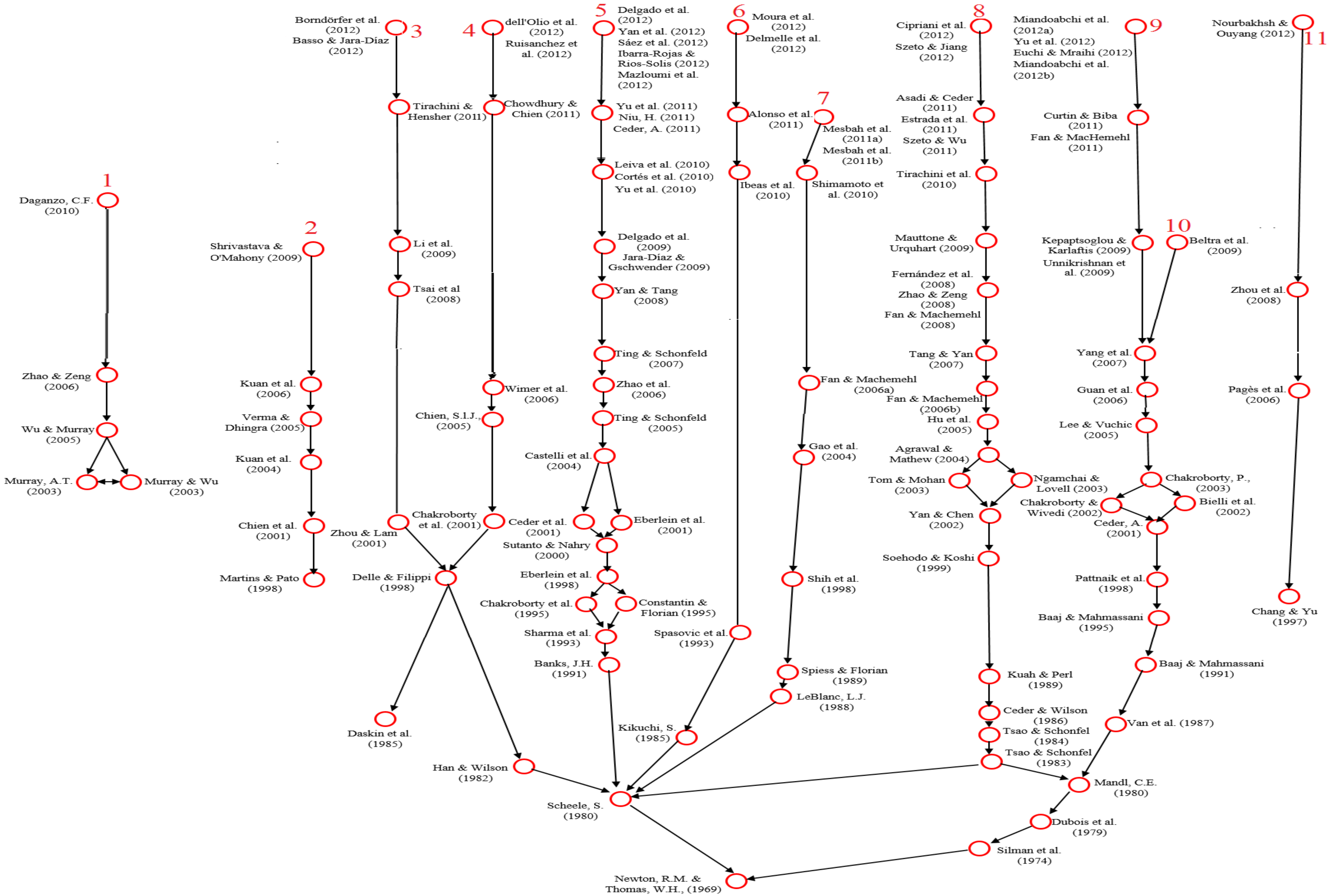

(A) Subcategories 1-11 until 2012

**Figure 8.** CPM of mathematical modeling and algorithms applications in PBTNDP&OPs.

## 4. Comprehensive Literature Review

In this section, the literature related to analytical approaches will be reviewed based on the identified CPM presented in Figure 7.

According to Figure 7, 13 major sub-categories have been identified in the related literature to the applications of analytical methods in PBTNDP&OP. As it can be seen from this figure, at the very beginning, there wasn't any clear direction in the conducted studies, and it has taken a while to see these major sub-categories in the literature.

Vuchic and Newell (1968) is the oldest identified study in this category. In this study, they focused on the intersection spacing problem in rapid public transit systems. They considered a line-haul and a many-to-one transit service and developed an analytical model for this problem with the objective of minimizing total passenger travel time. Considering a time dependent smooth function for users' arrival rate, Newell (1971) applied analytical methods to study the problem of vehicle dispatching for a transit route. They tried to find a dispatching policy that leads to the minimum total waiting time for the users. As a key finding, they showed that when there isn't a capacity constraint, the optimal headway will be approximately proportional to the arrival rate of passengers. Considering a shuttle transit service that serves users with two vehicles and assuming a constant user arrival rate, Newell (1974) used an analytical method to investigate the bus pairing problem. They tried to find a control strategy that leads to the minimum average users' waiting time. For a radial transit network, Byrne (1975) proposed an analytical approach in order to find optimal line spacing and headway for the transit services. The objective of the proposed model was to minimize the summation of users' travel costs and the agency's operation costs. They showed that in the optimal headway there will be a balance between users' and agency's costs. Also, optimal line spacing will be proportional to the population density and circumferential access. Clarens and Hurdle (1975) used continuum approximation approaches to study the optimal operation strategy problem for a commuter bus transit service. They developed a cost function of the zone size and used analytical methods to find optimal values. The results indicated that the zone size will be dependent on the occupied capacity of the vehicle, making it proportional to the square or cube root of the users' destination density.

Considering a network including parallel transit lines in a rectangular city, Byrne (1976) proposed an analytical approach for finding lines lengths, positions, and headways. The objective was to minimize the summation of user travel time and the agency's operating costs. They showed that in the case of uniform population distribution, headway will depend on average headway, variance of the headway, speed, and fleet size. Daganzo (1978) proposed an analytical model to estimate average users' waiting time and riding time in a many-to-many public transit service. First, they used a deterministic model to study different dispatching algorithms. Then, they developed a stochastic approach to obtain more realistic results. Using analytical methods, Newell (1979) studied the problem of optimal headway and line spacing for a rectangular transit network that led to minimizing the summation of waiting costs, access costs, transfer costs and operation costs. They considered two different transit network structures and developed an objective function for each structure. Taking a non-uniform demand pattern that varies smoothly over the space into account, Wirasinghe and Ghoneim (1981) proposed an analytical approach to find the optimal stop spacing on a local bus route that provides a many-to-many transit service. They considered a cost function as the summation of operation costs, waiting and in-vehicle costs, and used continuum approximations in their analysis. Sim and Templeton (1982) considered a transit service with a terminal where users' arrival pattern at the terminal follow a Poisson process. Also, they assumed that to control waiting lines, the agency could dispatch different types of vehicles. Having these assumptions, they applied analytical methods to determine the optimal fleet size and dispatching policy. Tsao and Schonfeld (1983) developed analytical approaches to obtain near to optimal values for headways and zone size for a many-to-one/one-to-many public transit service. To this end, they used calculus or a quasi-Newton algorithm to minimize the summation of operation, waiting, and in-vehicle costs.

As it was already mentioned, after the early published papers, 13 different sub-categories are emerged later. Following, the existing literature in each sub-category has been reviewed. Please note that the related literature to the applications of analytical approaches and mathematical programming in "Emerging technologies in public bus transit systems" and "Sustainability in PBTNDP&OP" will not be reviewed in this paper, as the related papers to these concept has been already published and are available in Mahmoudi et al. (2024) and Mahmoudi et al. (2025).

### 4.1. Dispatching policy

Without any doubt, dispatching policy is one of the most important and oldest problems in the mass transit systems, especially in public bus transit services. Setting an appropriate service frequency and a reliable timetable will lead to lower costs and higher satisfaction for both users and agencies. Dispatching policies are highly dependent on the topology of the network, fleet size, line and bus spacing, and demand patterns. On the other hand, the design of the network, fleet size, line and bus spacing, demand pattern, and other variables are also dependent on the designed dispatching policy.

These interdependencies make the dispatching problem more interesting but complicated. There are a significant number of studies that have applied analytical methods to investigate different problems related to dispatching policies in public bus transit services. de Palma and Lindsey (2001) considered a problem where there is a single transit line, the fleet size is known, and there are different groups of passengers based on their travel time and the cost of the delay. Based on these assumptions, they proposed analytical models in order to find the dispatching policy that minimizes the total schedule delay costs for the passengers. They considered two different possible scenarios for vehicle scheduling. First, a linear scenario, in which assumed that preferred travel times follow a uniform distribution over part of the day and it is not possible to reschedule the travel between days. Second, a circle scenario, in which it is assumed that preferred travel times follow a uniform distribution over the full 24 hours and it is possible to reschedule the travel between days. Finally, using the developed analytical model, they obtained closed form solutions for the decision variables in each scenario and compared the results. Jara-Díaz and Gschwender (2003) conducted a literature review on studies that had focused on the problem of dispatching policies and fleet size, considering a parametric demand and using microeconomic models, and tried to understand all the developments done in this area so far. After doing a comprehensive review considering a transit corridor with a fixed demand in a single period, they proposed a model based on Jansson's model to capture the impacts of vehicle capacity and crowding on the agency’s costs and time value, respectively. Finally, they tried to develop a general model to optimize bus operations. Criticizing the shedule- and headway-based control strategies in dealing with the bus bunching problem in public transit services, Daganzo and Pilachowski (2011) proposed a cooperative strategy based on the idea of sharing information between leading buses and the ones that are following them. To this end, they assumed that there is a closed loop served by a known number of buses and including a known number of stops, where the distance between all consecutive stops is a constant value, the demand is distributed uniformly along the loop, and the arrival time of the users is spatially homogeneous. Having these assumptions, they developed an adaptive control strategy to modify bus cruising speed based on the shared real-time information. The study showed that under the strategies obtained by the proposed model, the transit service would perform much better than the one under the results of the previous models.

Analyzing crowding effects on public transit operations planning has always been an interesting subject for researchers. Using analytical models, De Palma et al. (2015) tried to investigate the impacts of crowding on vehicle capacity, fare pricing, and frequency setting. They considered a discomfort cost/comfort value for the users and obtained the equilibrium values for the decision variables in different service types (many-to-one, one-to-many, and multiple-lines) and different assumptions for the arrival time and allocation of seats. Schmöcker et al. (2016) studied the problem of bus bunching, assuming that the network includes overtaking and common lines (see Figure 9). They developed some analytical approaches in order to find the optimal departure times of a group of buses when a delay has occurred for a leading bus at a stop. The possibility or impossibility of overtaking at the stops led to different models. When overtaking was impossible, they assumed that users would first try the front bus. In another case, authors claimed that users would show an equilibrium behavior in choosing buses while trying to optimize their waiting times. To show how network design can affect the regularity of the services, authors applied the proposed approach to a case study and showed that having a common line would have negative impacts on the regularity when overtaking was not possible, while it would be a positive factor for the transit system if overtaking was allowed.

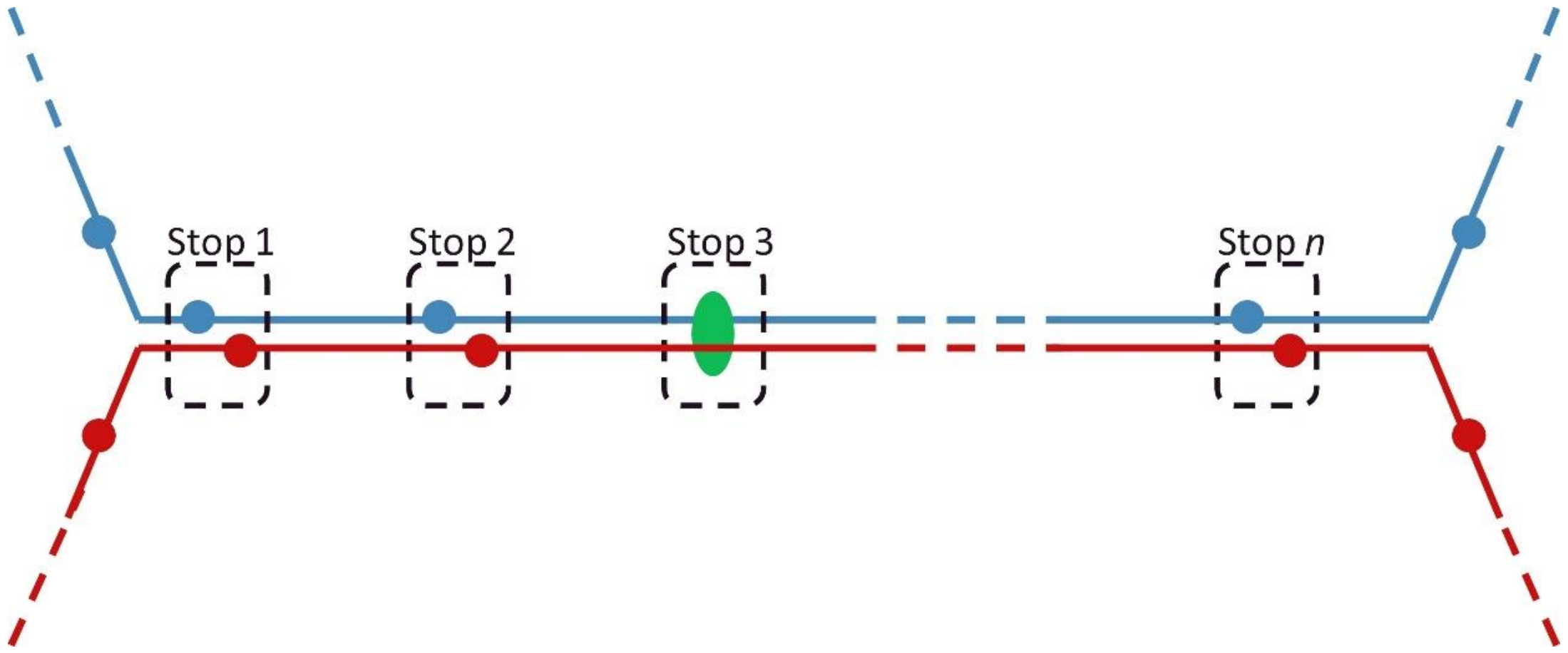

**Figure 9.** A network including overtaking and common lines: two bus lines with some common line stops (Schmöcker et al., 2016)

Zhang and Lo (2018) conducted another study related to the bus bunching problem, taking both deterministic and stochastic running times into account. Using analytical methods, they proposed a control method to obtain the optimal control parameters according to the number of buses on the route. The basic idea for their model was to see both upstream and downstream of a control point to find the differences between the headways. As the key results, the study showed that when there is a deterministic travel time assumption, the headways will self-equalize, while for the stochastic travel time assumption, it was shown that the control strategy would limit the headway variance to a specific amount. He et al. (2020) studied the problem of bus bunching for a public bus transit line that serves users with a high frequency. It was also assumed that buses do not follow any predefined schedule. Having these assumptions, they proposed dynamic analytical approaches to find the best holding strategies based on the difference between a target headway and an expected current headway for a bus. The proposed model was able to deal with a situation where there were multiple control points on the line. However, the study indicated that increasing the number of control points on a line wouldn't always lead to an improvement in the performance of the transit network. When it comes to PBTNDP&OP, Ansari Esfeh et al. (2021) claimed that the exact type of passenger and transit service must be considered in the modeling process. Also, they mentioned that considering waiting time only at the origin and ignoring the waiting time at the destination would lead to unrealistic/inaccurate results. Reviewing the existing literature, they showed that most of the studies do not consider the various types of passengers and transit services in their studies, which calls into question the validity of their findings. They categorized passengers into two groups: those who plan their trips and those who do not, and service types into four categories including schedule-based, frequency-based, high-frequency, and low-frequency services, then suggested an approach to find the average waiting time.

Table 7 represents the main features of the reviewed papers in this subsection. According to this table, most of the reviewed papers have considered a uniform demand pattern and have found closed form solutions for the decision variables. However, in most these studies, the optimizations have been done for an idealized single route transit network. Min bus bunching, min total schedule delay cost, and minimization of the sum of access costs, waiting costs, in-vehicle costs, and operation costs are the most commonly used objective functions.

**Table 7.** Main features of the reviewed papers in subcategory 1.

| Reference | Objective Function | Transit mode | Network structure | Demand Pattern | Decision Variables | Solution Approach |
|---|---|---|---|---|---|---|
| de Palma and Lindsey (2001) | Min. Total schedule delay cost | General | Single line | Uniform distribution | Headway/ Frequency | Closed form solution |
| Jara-Díaz and Gschwender (2003) | Min. summation of access cost, waiting cost, in-vehicle cost and operation cost | General | Single corridor | Uniform distribution | Headway/ Frequency, Fleet size & vehicle size | Closed form solution |
| Daganzo and Pilachowski (2011)* | Min. bus bunching | Bus | A closed loop | Uniform distribution | Bus speed | Closed form solution |
| De Palma et al. (2015) | Min. Crowding cost and Schedule delay cost | General | many-to-one, many-to-one, and multiple-lines | Uniform distribution | Headway/ Frequency, Fare & vehicle size | Closed form solution |
| Schmöcker et al. (2016) | Min. bus bunching | Bus | Two lines | Uniform distribution | Departure times | Closed form solution - Iterative algorithm |
| Zhang and Lo (2018) | Min. bus bunching / Headway variance | Bus | A closed loop | - | Bus control parameters | Closed form solution |
| He et al. (2020) | Min. bus bunching | Bus | A closed loop | A time dependent demand function | Bus control parameters | Closed form solution |
| Ansari Esfeh et al. (2021) | - | General | - | - | Headway & Waiting time | Closed form solution |

*: The paper has considered an emerging technology

### 4.2. Feeder transit system planning

Bus transit services can serve passengers both as an independent public transit service or as a feeder system to take passengers to another mass transit system such as LRT, BRT, Metro, etc. Especially, in the big and populated cities where there is at least a rapid public transit system, feeder transit system planning will be an important problem for the transit system operators. From the very beginning, there was significant interest in applying analytical methods to study problems related to feeder transit planning. Although it is an old problem in this area, there are still a significant number of studies that have focused on modeling this problem using analytical models.

Wirasinghe et al. (1977) conducted one of the first studies related to feeder transit system design. They assumed that travel to the central business (CBD) district of a city was served by a bi-modal public transit system, including rail transit services and public bus services. It is also assumed that buses can be directly used to reach the CBD or that passengers can take feeder buses to a train station and then continue their trip with the train. Considering an idealized metropolitan region as shown in Figure 10 and assuming that passenger demand is in a deterministic form and varies slightly with location, they proposed an analytical model to obtain the optimal parameters for this transit service that led to the minimum total cost function. They defined the total cost function as the summation of the users' and agency's costs. It was also assumed that the railway network is a radial network and there is a radiocentric regional highway network in the CBD district. Having these assumptions and using continuum approximation, for the considered network, they obtained railway station spacing, feeder service boundary, and train dispatching policy. Another similar study has been done by Wirasinghe (1980). He considered a bi-modal public transit system where a rail system serves a many-to-one transit demand during peak time and is being fed by a public bus system as shown in Figure 11. As it is clear from this figure, a rectangular highway network was assumed where the rail line is parallel to one of the axes. Applying analytical methods and graphical inferences, he proposed an approach to obtain near to optimal values of the parameters, such as dispatching policies on each side of the rail line and density of the bus routes on each side.

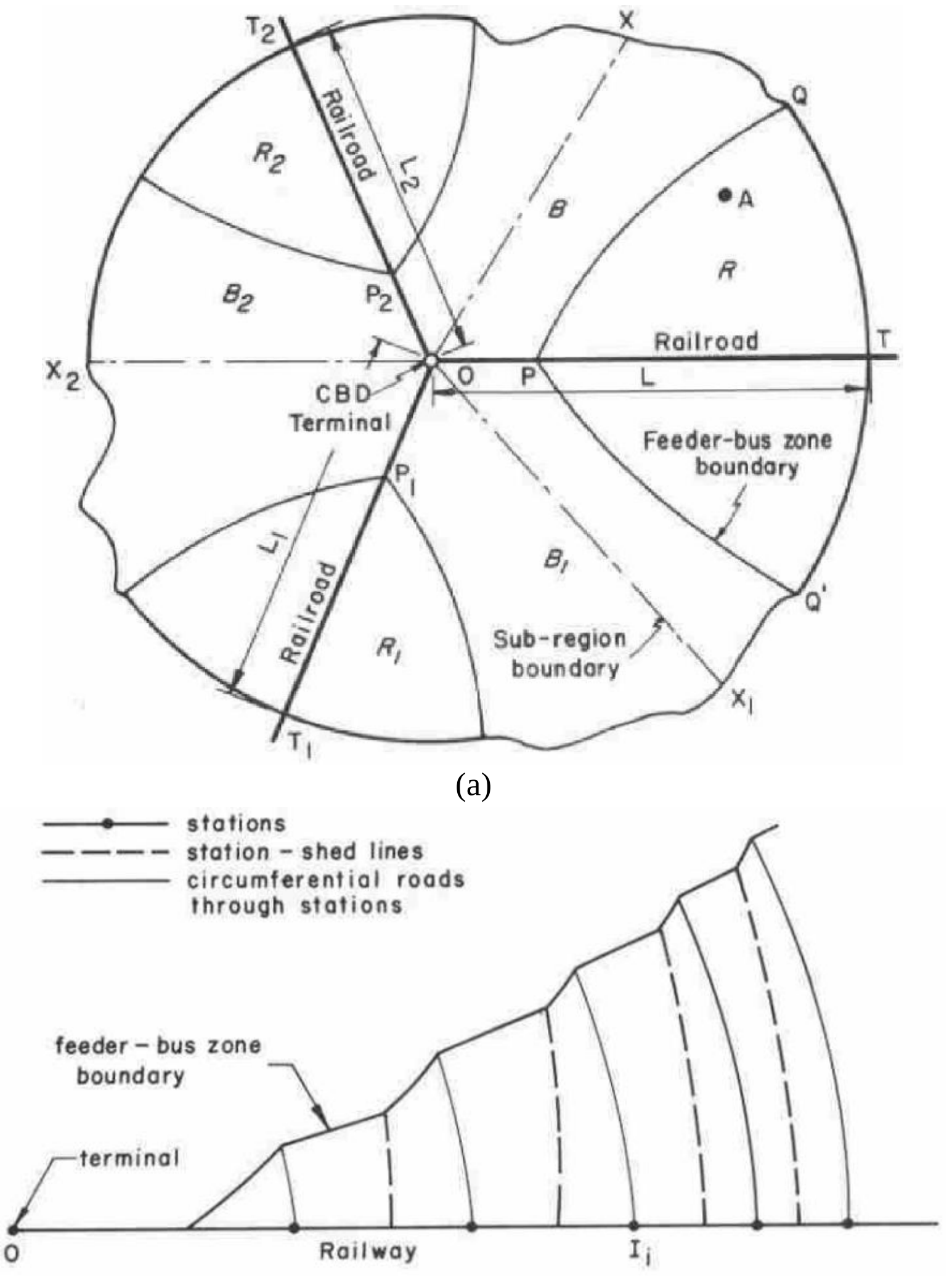

(a)

(b)

**Figure 10.** Idealized metropolitan region and feeder transit system considered in Wirasinghe et al. (1977)

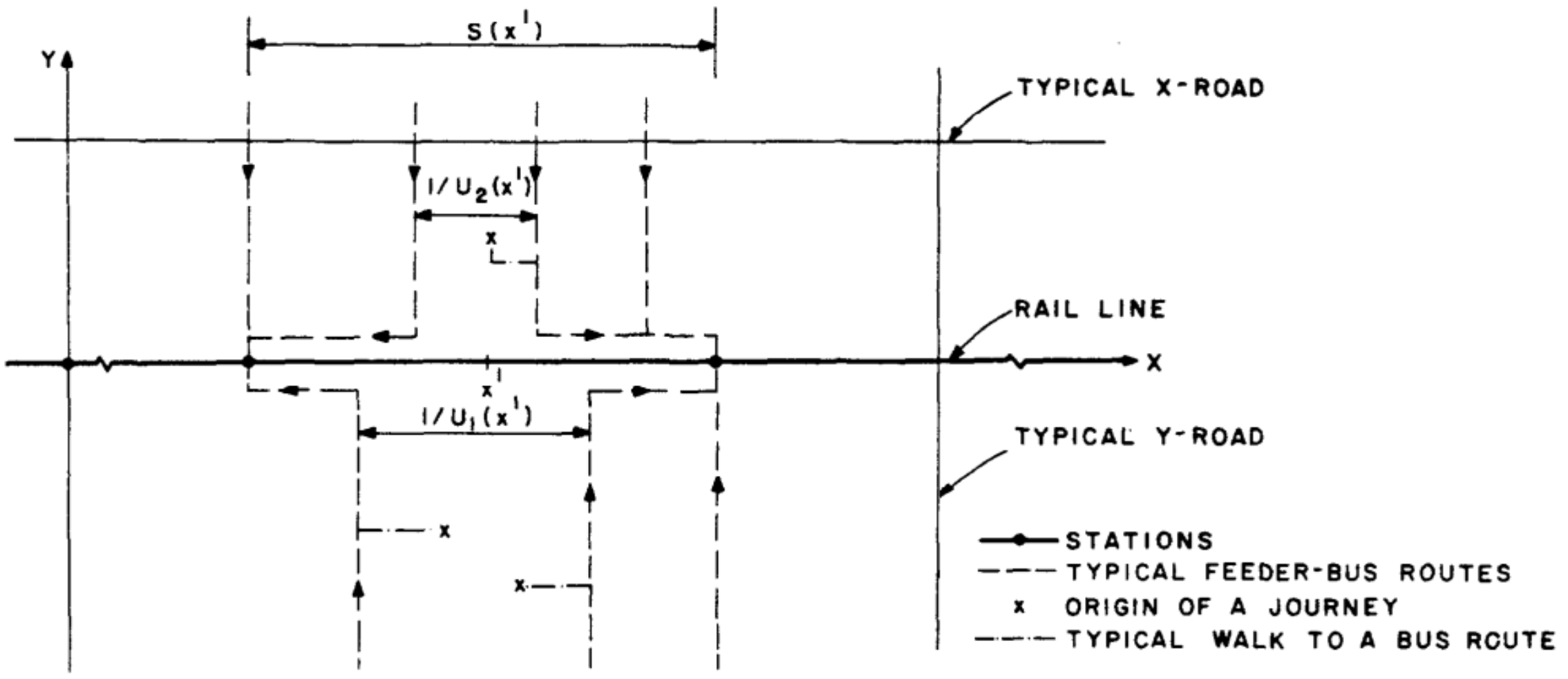

**Figure 11.** Integrated feeder bus and rail network for a many-to-one transit service considered in Wirasinghe (1980)

Similar to the previous studies, Kuah and Perl (1988) considered an integrated feeder bus and rail transit network as shown in Figure 12 and applied analytical models to obtain the optimal network features, including line and stop spacing, and frequency setting. As a contribution, instead of considering these variables separately and using the sequential approach, they combined all three variables and used an analytical approach to optimize them. For bus-stop spacing, three different scenarios were considered: (a) it is a constant value for the entire network, (b) it is constant along any specific line, and (c) it could take any value for any consecutive stops in the entire network. As a key finding, they showed that stop spacing is not significantly sensitive to changes in the values of other parameters of the transit network, such as route location and demand pattern. Kuah and Perl (1989) assumed that there was a rail transit system serving a city, and the problem they considered was to design an efficient feeder-bus transit network to feed this system under many-to-one and many-to-many demand patterns (see Figure 13). To this end, they proposed a mathematical optimization model to solve the network design problem and applied analytical methods to find the optimal frequency of public bus transit service. Because of the complexity of the proposed integrated method, they suggested a heuristic algorithm that was able to find acceptable solutions under variable demands. Chang and Schonfeld (1991b) proposed analytical approaches to make a comparison between the performance of the fixed-route and flexible-route subscription public bus transit services. They assumed that both services, as feeder services, were transferring users to a specific transportation node on the network, such as a railway station or a terminal. They considered the average cost per trip as the performance comparison criterion, and taking the total cost function as the summation of the agency and passengers' costs, they found closed form solutions for the optimal vehicle size and service zone size. They also modified the proposed models for two specific cases. For the first case, they assumed that travel demand and costs are time dependent. In the second case, they assumed that the maximum load factor could be changed according to existing stochastic demand. The result showed that the flexible route transit system can be the preferable system when the service zone is small, ride time, operation and vehicle costs are low, and transit speed, waiting time, and accessibility are high.

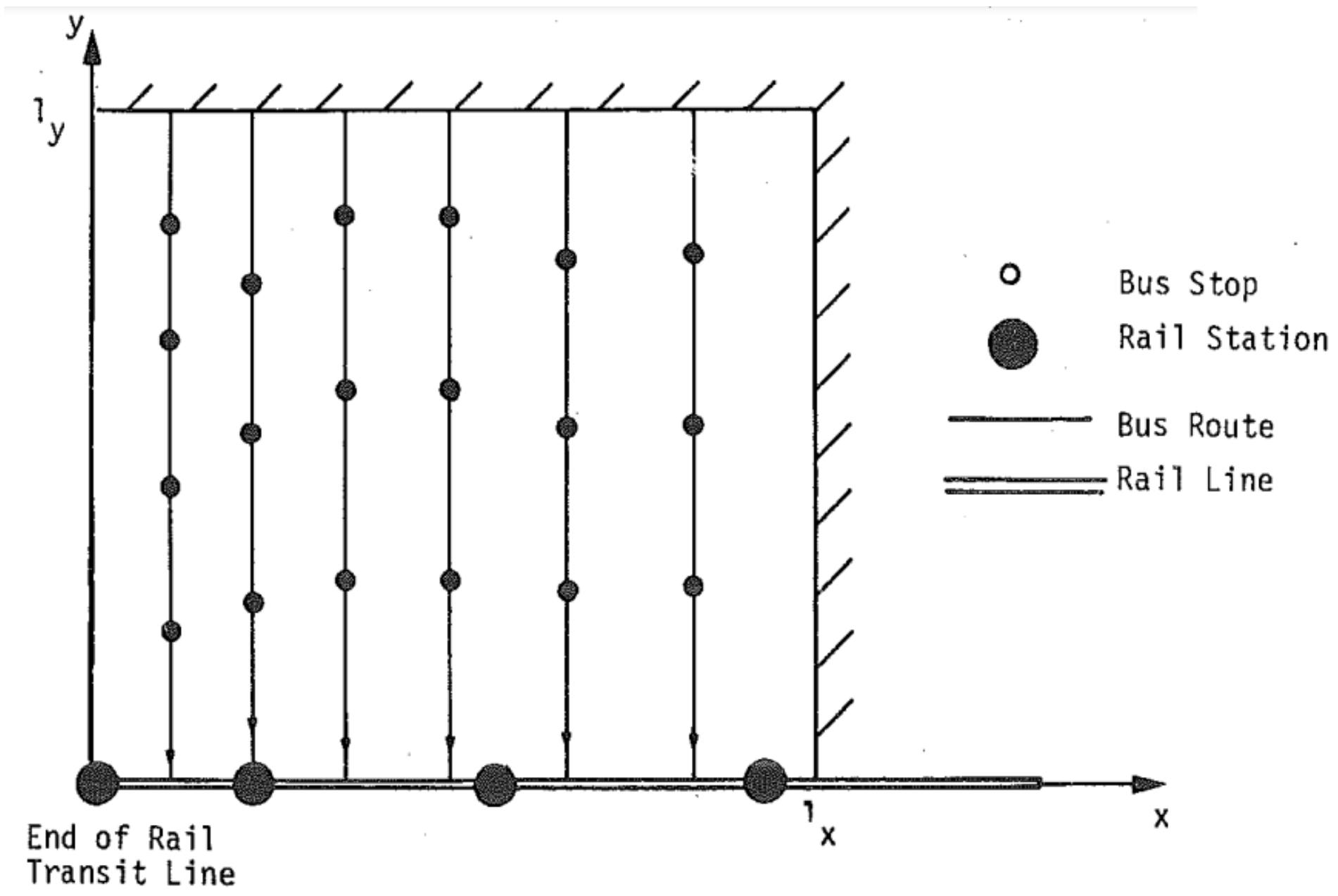

**Figure 12.** Bi-modal public transit network structure considered in Kuah and Perl (1988)

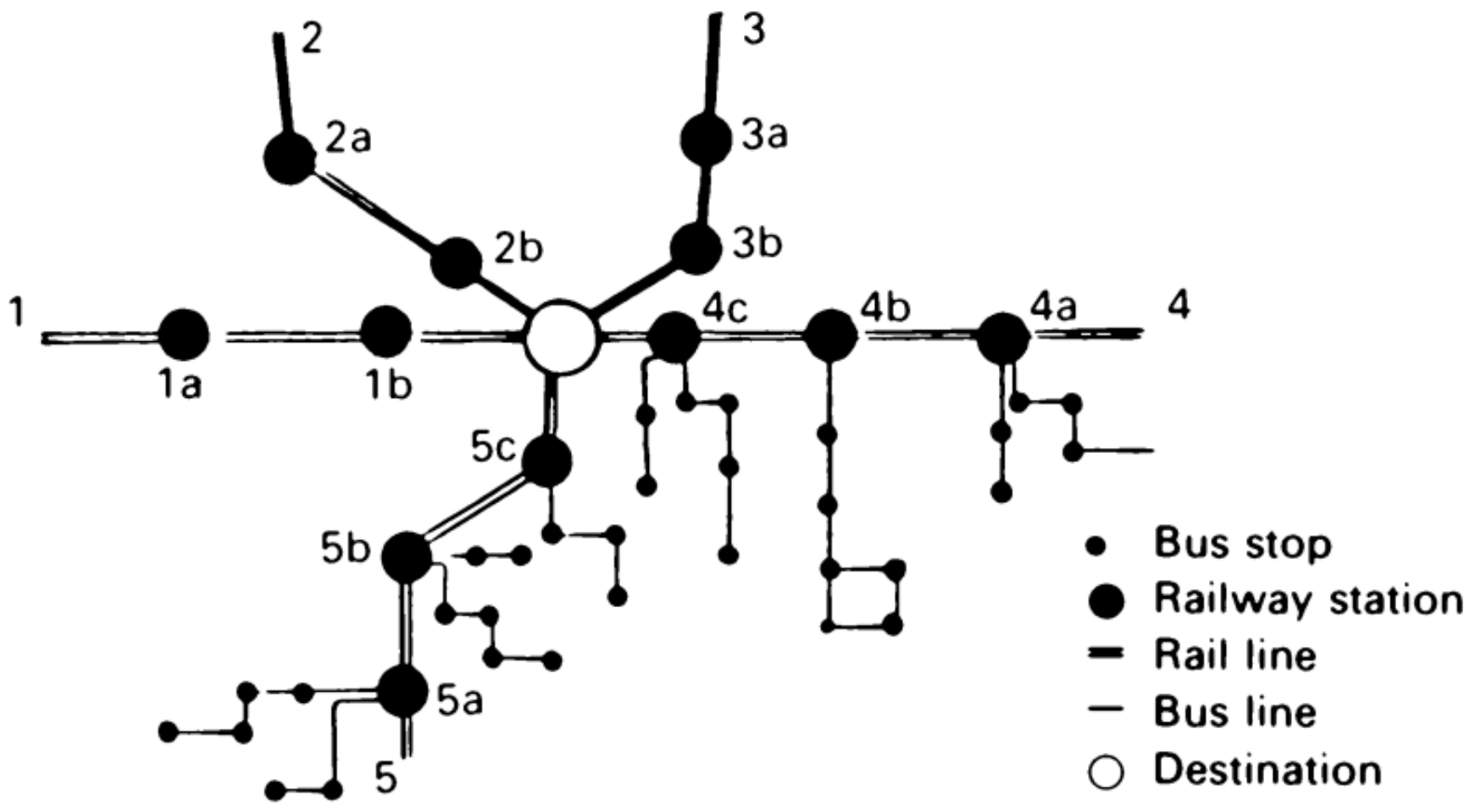

(a)

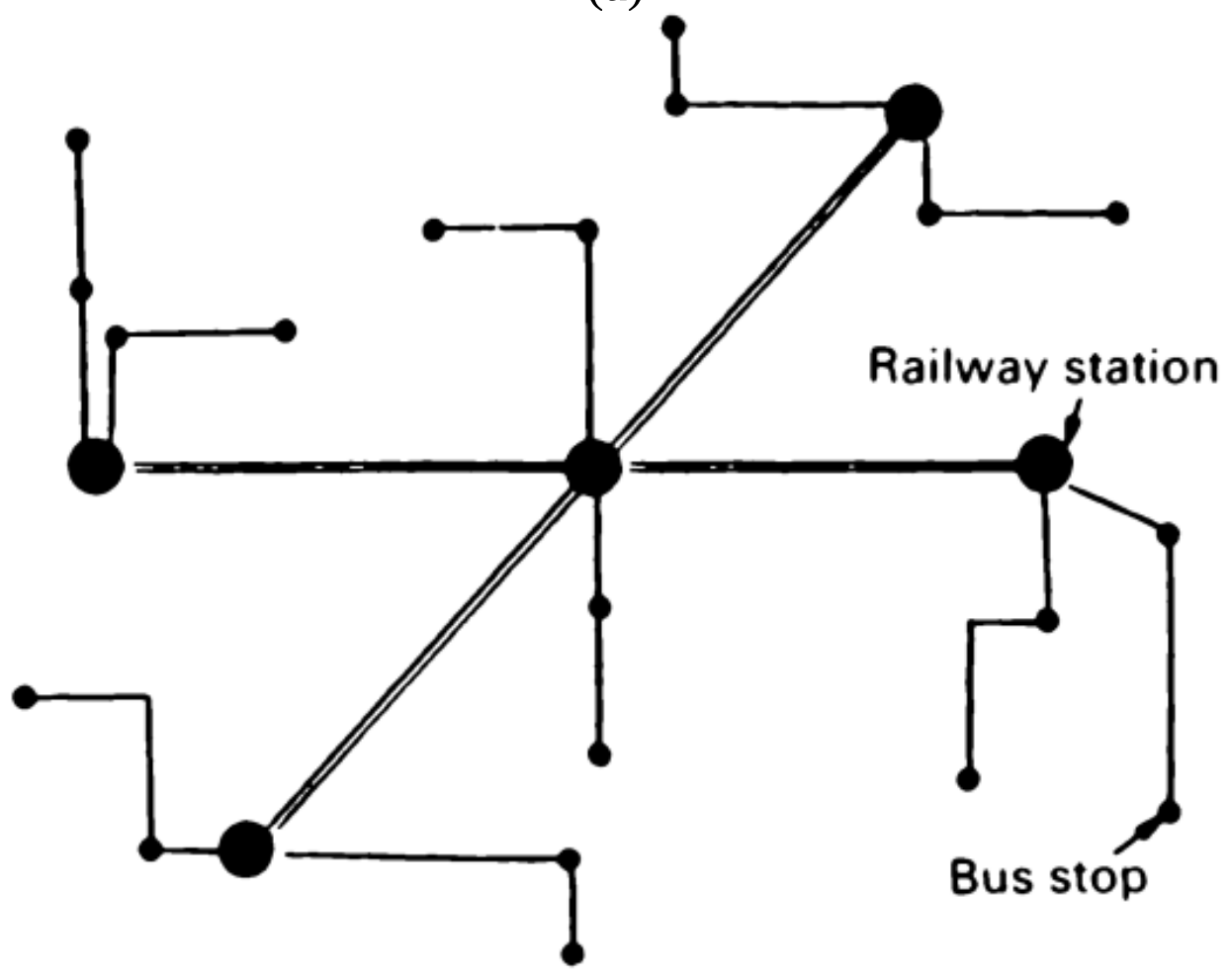

(b)

**Figure 13.** A public transit network including a rail transit service and a feeder bus network serving many-to-one (a) and many-to-many (b) demand patterns considered by Kuah and Perl (1989)

Considering irregular discrete distributions for demand and a bi-modal public transit network as shown in Figure 12, Chien and Schonfeld (1998) presented an analytical model with a capability of jointly obtaining optimal features for both rail and feeder bus transit services. The objective of the proposed model was to minimize the total cost function that was the summation of users' and agency's costs in the considered integrated transit service. The decision variables were railway length, rail stations and bus stop locations, bus line spacing, and finally bus dispatching policy. An iterative algorithm and a computer program were also developed to solve the proposed model. The results showed that the algorithm is able to find many near to optimal solutions and there is a slight difference between the performance of the network under each obtained near to optimal solution. Hence, the locations of obtained lines can be shifted to the nearest existing street without significant increases in the objective function. Chien and Yang (2000) considered an irregular street network with a many-to-one travel pattern (see Figure 15) and developed an analytical model to design an optimal feeder bus transit network by taking intersection delays into account. Using a cost function and an irregular discrete distribution for travel demand similar to the previous study, the proposed model was used to obtain optimal bus route locations and dispatching policies on each route. They also developed computer programs to find the global optimal solutions to the routing problems. Using an analytical model considering an average total cost function (summation of the users' and agency's costs) and a probabilistic demand pattern, Chien et al. (2001) tried to make a comparison between fixed-route and flexible feeder bus systems in feeding another mass transit system, as shown in Figure 16. The objective of the proposed model was to obtain optimal vehicle size and route spacing for the fixed-route system and optimal vehicle size and service zone for the flexible-route one. From the proposed model, a threshold value of demand was obtained to find which system was performing better at any given demand value, from the average cost perspective.

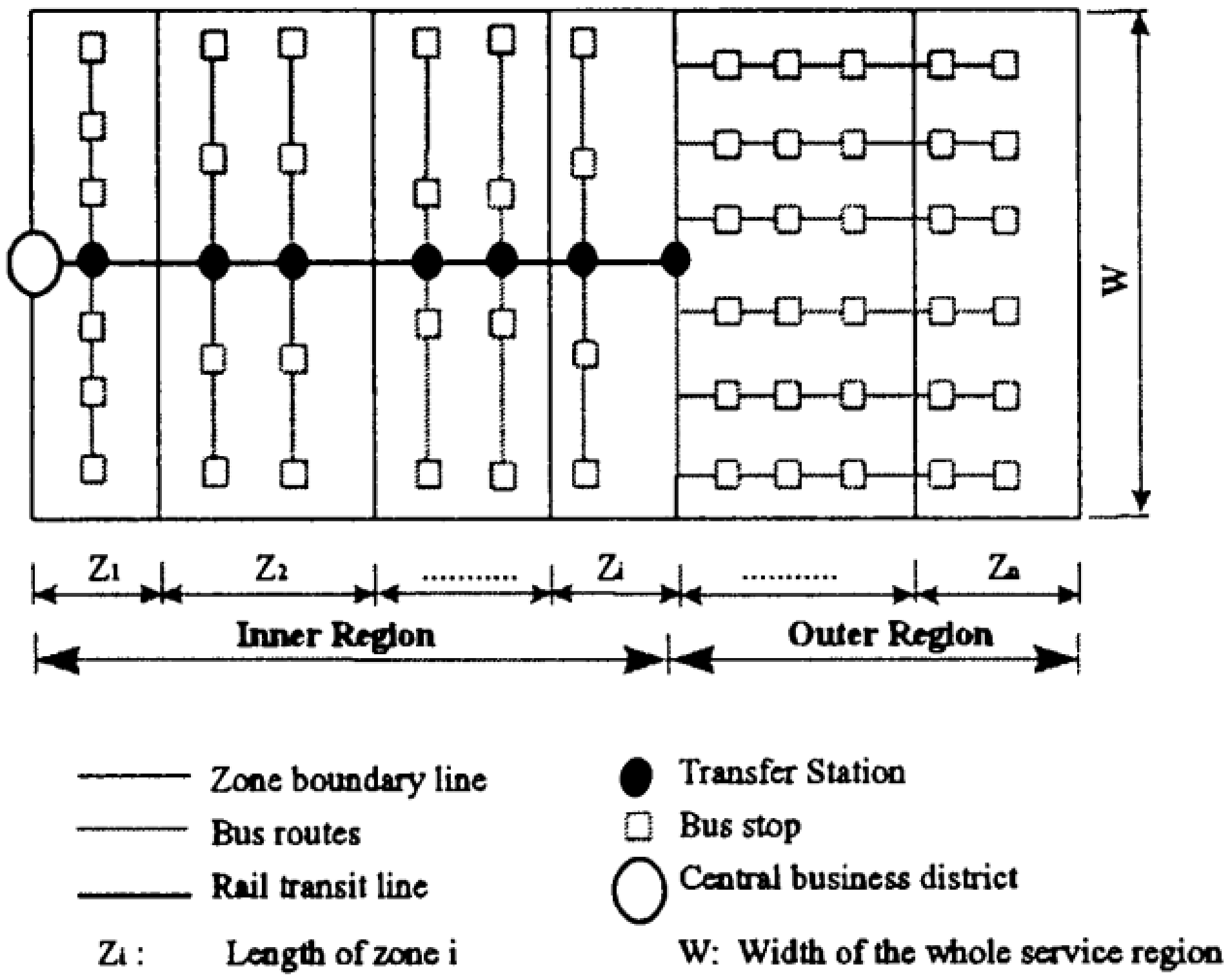

**Figure 14.** The bi-modal transit network considered in Chien and Schonfeld (1998)

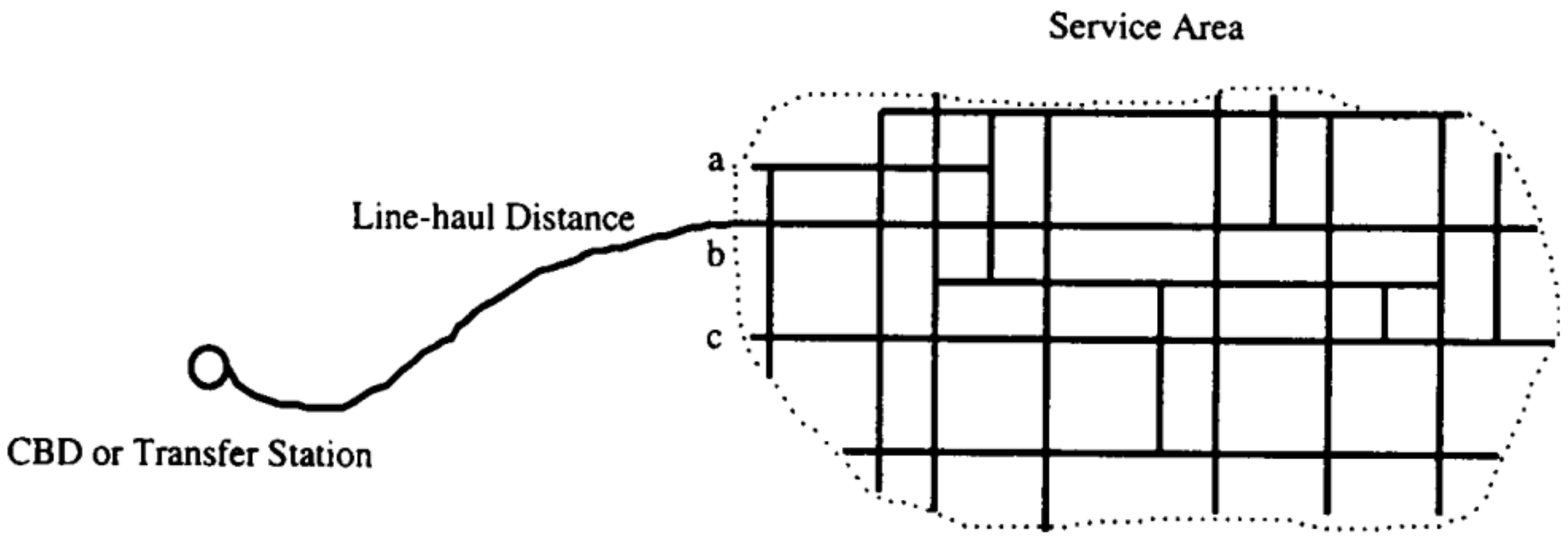

**Figure 15.** Irregular network considered in Chien and Yang (2000)

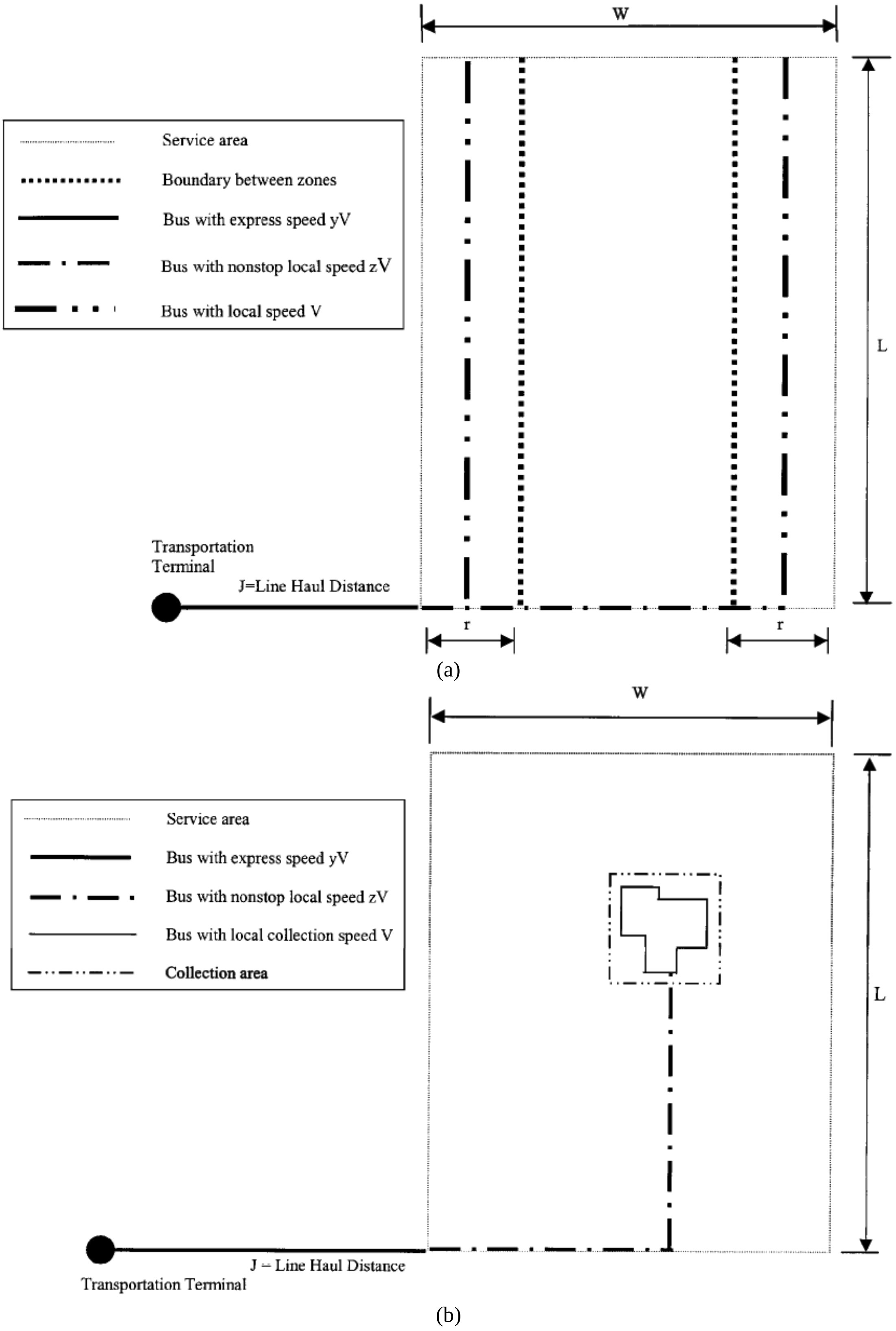

**Figure 16.** Conventional (a) and subscription bus (b) systems considered by Chien et al. (2001)

Li and Quadrifoglio (2010) developed an approach to assist transportation system managers in determining which feeder transit service (flexible or fixed) was more suitable for a network using analytical models and simulation methods. The proposed approach was able to determine if/when it is required to change from one service to another. To develop this approach, Li and Quadrifoglio (2010) considered a weighted total cost function including walking, waiting, and riding costs. As a key finding, they showed that demand will highly affect the switching point between fixed and flexible transit services in the network. Li and Quadrifoglio (2011) investigated the problem of feeder transit zone design where two vehicles are serving each zone, demand is high, and the length of the service area is significantly long. To this end, they developed an analytical model with the objective of making a balance between service quality and agency costs. They used continuous approximation principles in the modeling process and obtained closed-form solutions for the zone numbers. Sivakumaran et al. (2012) considered a mass transit system that provides a many-to-one transit service and is being fed by another public transit network and studied the problem of vehicle schedule coordination. First, they assumed that the vehicle schedule for the trunk service is known and given, and showed that in this case, coordinating the feeder system with the trunk system will lead to reduced passenger costs outweighing the increase in the operation cost of the feeder system. Second, they assumed both trunk and feeder service frequencies as the decision variables and showed that in this case, in the optimal dispatching strategy, both passengers and operation costs can be reduced. Chandra and Quadrifoglio (2013) first defined service quality as the inverse of the summation of waiting and in-vehicle times in a demand responsive feeder transit system. Then, they developed an analytical model to obtain the length of the terminal-to-terminal cycle for the feeder system, which maximizes the service quality. The parameters related to the topology of the network and the demand pattern were the inputs of the proposed model.

Zhao and Chien (2015) considered a many-to-one/one-to-many demand pattern from the residential area to the CBD and vice versa during the peak times and proposed an analytical model to obtain optimum stop spacing and service frequency for the feeder system. As a contribution, they considered the variance of inter-arrival time and its effects on the users' travel time and buses' dwell time in the modeling process. Their results showed that compared to existing studies, taking the variance of inter-arrival time impacts into account would lead to a longer stop spacing, while controlling this variance could lead to a shorter stop spacing and better accessibility for the users. Considering a simple network structure and different line structures as shown in Figure 17, Gschwender et al. (2016) used analytical methods to make a compression between a feeder-trunk system and three direct line systems to find which structure is more cost-efficient. The methods used to obtain optimal fleet and vehicle sizes that minimize the summation of users and the agency's cost.

The study indicated that the results of the comparison are highly dependent on the demand volume and pattern, trunk line length, and transfer cost.

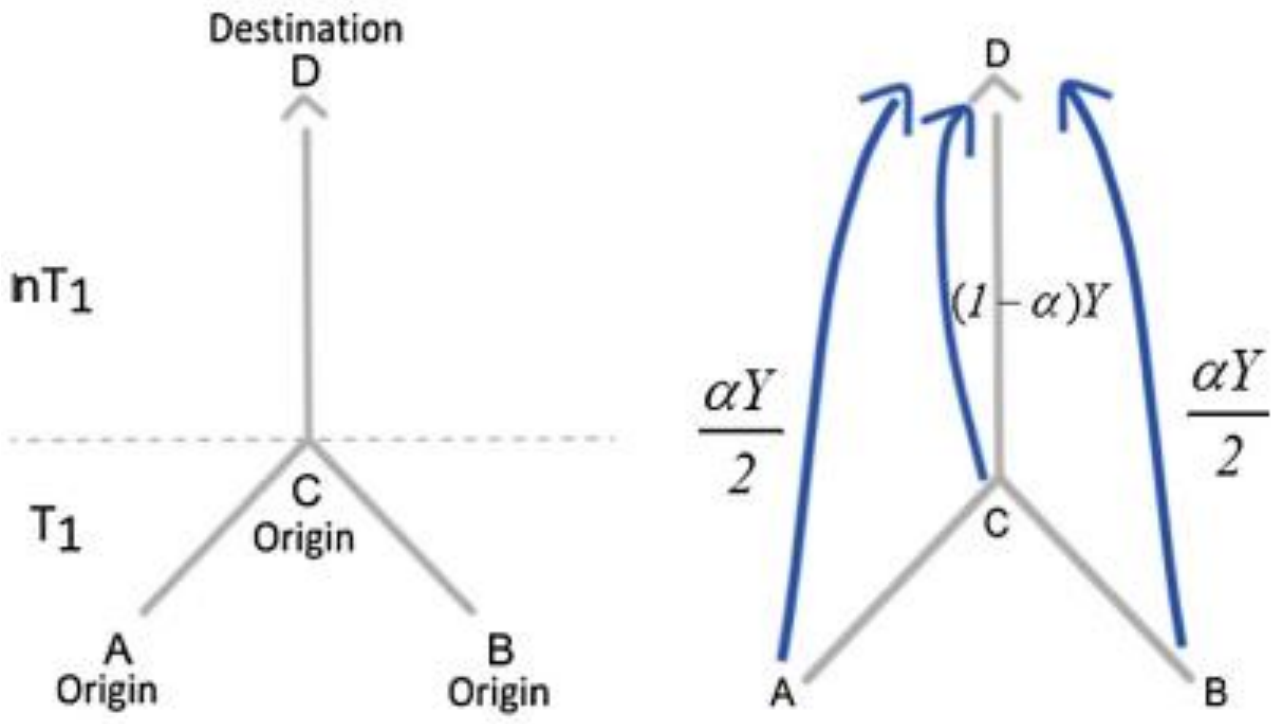

(a)

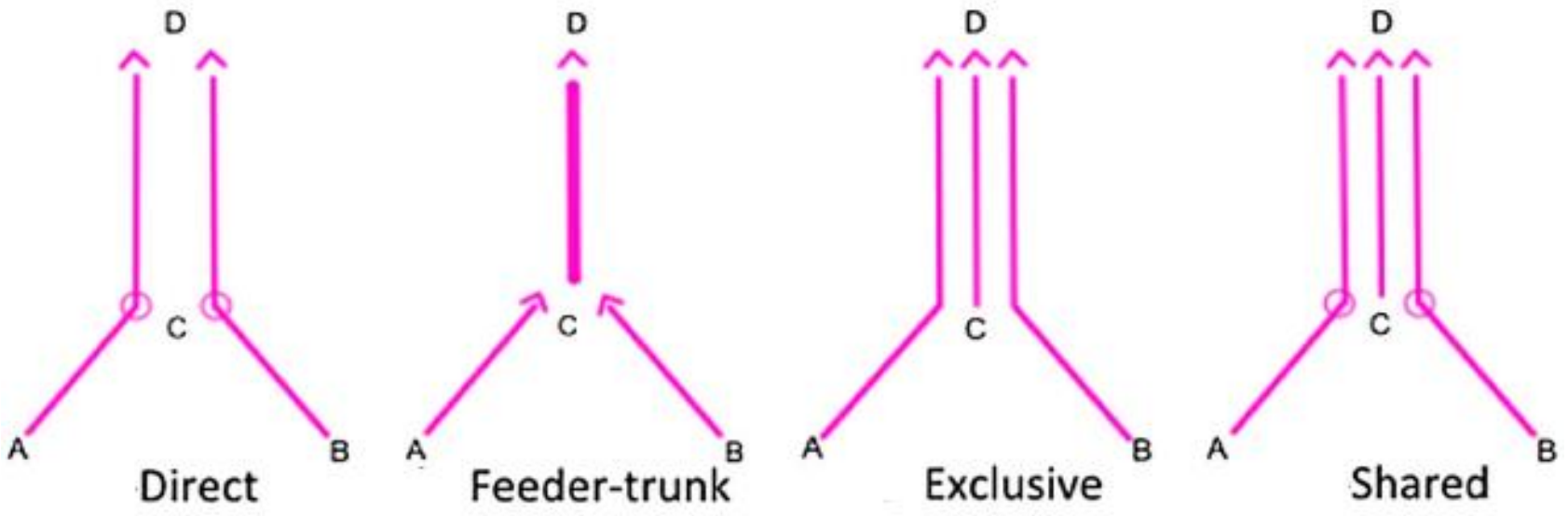

(b)

**Figure 17.** Network structure and flow pattern (a) and line structures (b) considered by Gschwender et al. (2016)

Taking the Purple Line Corridor in Maryland as the case study, Li et al. (2020) proposed an approximation model to deal with the problem of feeder bus network design along with an existing trunk service. They considered the Purple Line Corridor in Maryland in an idealized form (see Figure 18), then, applied the proposed model to obtain optimal density and headway of feeder services with respect to the bus capacity constraints. To address the problem of schedule coordination in a bi-modal transit network where a feeder system services a truck transit system, Yang et al. (2020) proposed a mathematical approach to minimize the total cost function by considering a spatially heterogeneous demand pattern for the service area. Figure 19 shows the hypothetical network structure considered in this study. First, they integrated discrete optimization and the continuous approximation approach to find the optimal uncoordinated schemes and a feasible layout. Then, they added trunk station locations, feeder transit routes' locations and lengths, and the passenger flow distribution in the next stage of designing the coordination. To solve the proposed model, an integrated Genetic Algorithm (GA), analytical approach, and optimization algorithm are proposed. Solving the model, they showed that when there is a lower-level heterogeneous demand, the best strategy will be the coordination one.

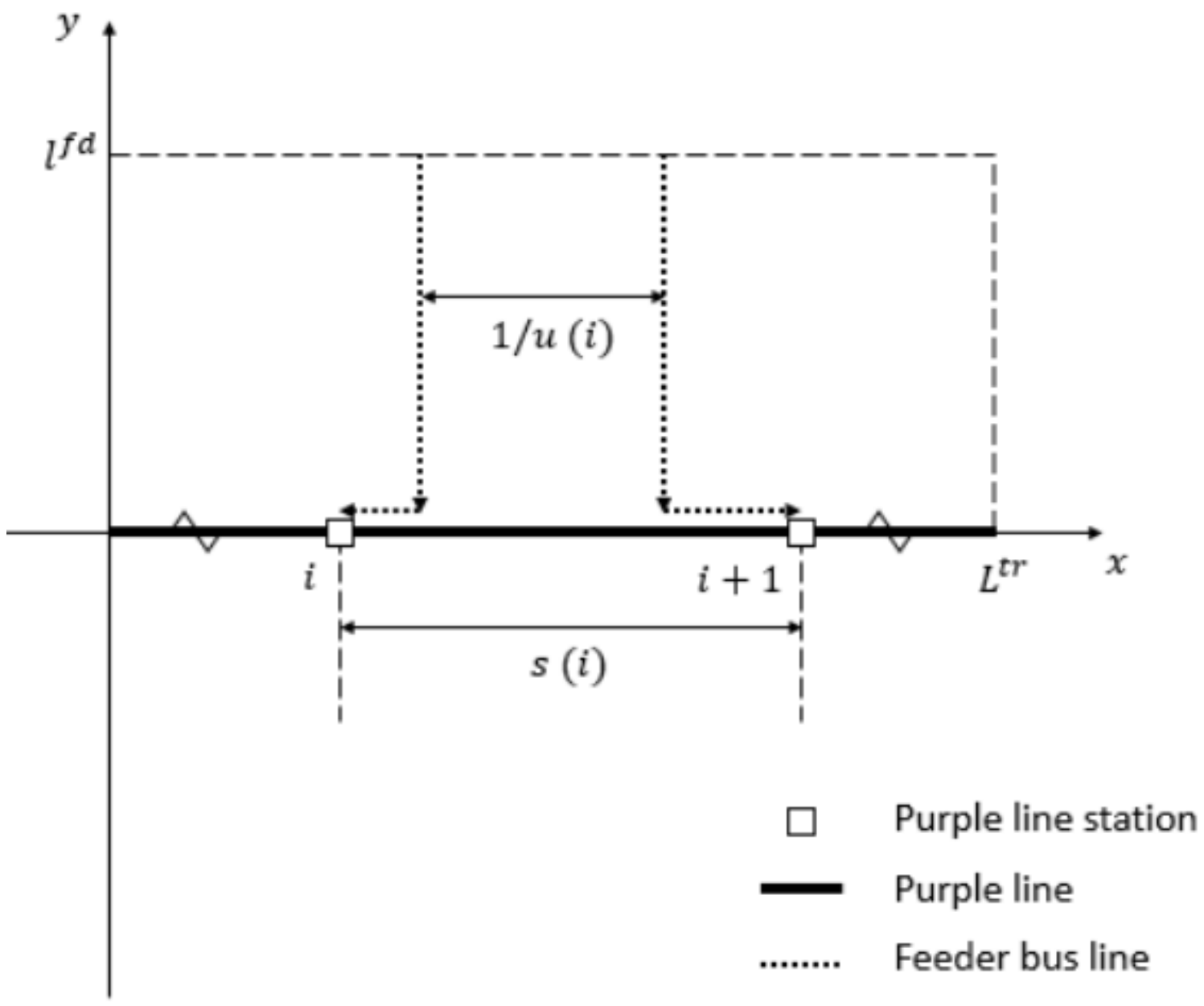

**Figure 18.** The idealized network structure along Purple Line Corridor in Maryland in Li et al. (2020)

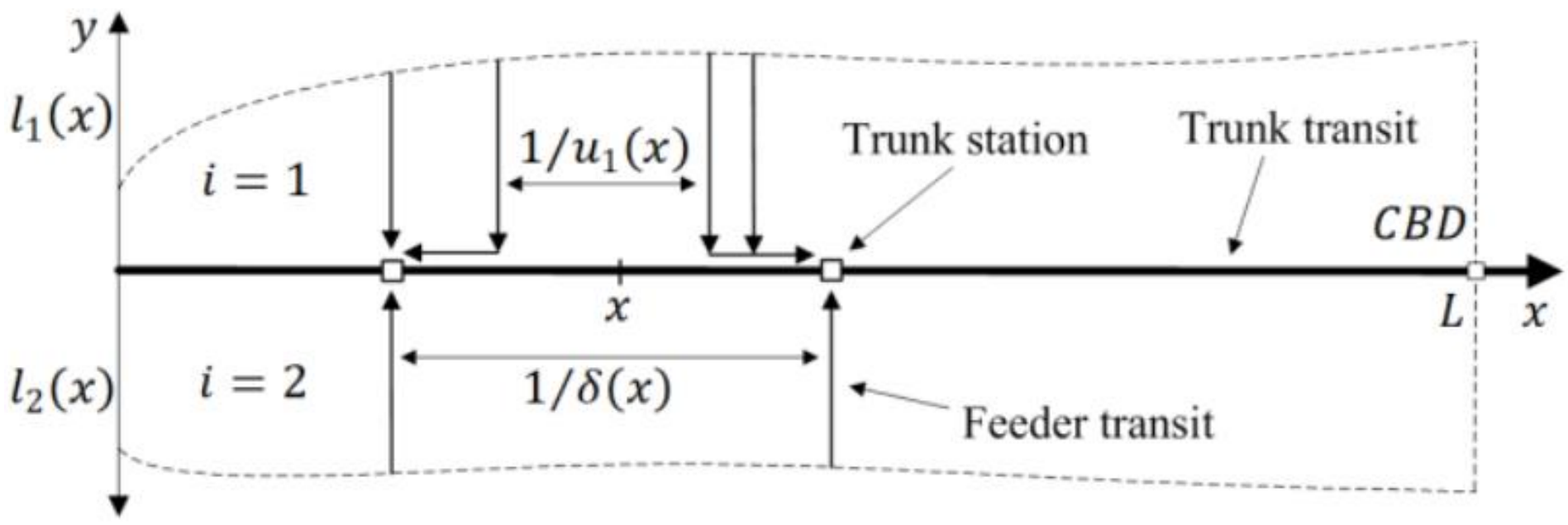

**Figure 19.** Hypothetical network structure considered in Yang et al. (2020)

Table 8 represents the main features of the reviewed papers in this subsection. According to this table, minimization of the sum of access costs, waiting costs, in-vehicle costs, and operations costs is the most commonly used objective function in the literature. Most of the reviewed papers considered an idealized rectangular network. A uniform demand pattern and a location-dependent demand function are the two most considered demand patterns by the reviewed papers. In addition, service frequency, stop location, route density, vehicle and fleet sizes, and service zone size are the main decision variables in this sub-category. Based on Table 8, most of the studies were able to find closed form solutions for the decision variables using analytical approaches.

**Table 8.** Main features of the reviewed papers in subcategory 2.

| Reference | Objective Function | Transit mode | Network structure | Demand Pattern | Decision Variables | Solution Approach |
|---|---|---|---|---|---|---|
| Wirasinghe et al. (1977) | Min. sum of users' and agency's cost | Bus & Rail | Radiocentric & Radial | Deterministic Demand varying slightly with location | Railway stations spacing, feeder service boundary & train dispatching policy | Closed form solution |
| Wirasinghe (1980) | Min. sum of users' and agency's cost | Bus & Rail | Rectangular highway & linear rail network | Deterministic many-to-one demand varying slightly with location | Headway/ Frequency, Bus routes density | Closed form solution |
| Kuah and Perl (1988) | Min. sum of operation, waiting, in-vehicle & access costs | Bus & Rail | Linear | Continuous inelastic location dependent function | Headway/ Frequency, Line spacing & Stop spacing | Closed form solution |
| Kuah and Perl (1989) | Min. sum of operation, waiting, in-vehicle costs | Bus & Rail | Non-idealized | Location dependent deterministic demand | Headway/ Frequency, Stop and station spacing | Closed form solution - Heuristic algorithm |
| Chang and Schonfeld (1991b) | Min. sum of operation, waiting, in-vehicle & access costs | Bus & Trunk | Rectangular | Time dependent & Stochastic Demand | Vehicle size, Route spacing & Service zone size | Closed form solution |
| Chien and Schonfeld (1998) | Min. sum of operation, waiting, in-vehicle & access costs | Bus & Rail | Rectangular | Uniform distribution within a zone- Many-to-many | Railway length, stations/stops locations, bus line spacing & Headway/ Frequency | Closed form solution |
| Chien and Yang (2000) | Min. sum of operation, waiting, in-vehicle & access costs | Bus & Trunk | Irregular network divided into rectangular zones | Discrete uniform distribution within a zone - Many-to-one | Headway/ Frequency, Bus route spacing | Closed form solution - Heuristic algorithm: Near to Optimal |
| Chien et al. (2001) | Min. sum of operation, waiting, in-vehicle & access costs | Bus & Trunk | Rectangular | Time dependent probabilistic demand | Vehicle size, Route spacing & Service zone size | - |
| Li and Quadrifoglio (2010) | Walking, waiting, and in-vehicle costs | Bus & Trunk | Rectangular | Poisson distribution - Uniform distribution | - | Closed form solution - Simulation |
| Li and Quadrifoglio (2011) | Min. sum of operation, waiting, in-vehicle & access costs | Bus & Trunk | Rectangular | Uniform distribution | Number of zones & Number of bus stations, | Closed form solution |
| Sivakumaran et al. (2012) | Min. sum of operation, waiting, | General feeder-trunk | Rectangular-parallel feeder lines | Continuous, time-independent | Headway/ frequency, feeder | Closed form solution |
| | transfer & access costs | | | density function & Uniform distribution | route spacing & stop spacing | |
| Chandra and Quadrifoglio (2013) | Max. service quality | General feeder-trunk | Rectangular | Spatially and temporally uniformly distributed (Poisson process) | Cycle length | Near to optimal - Simulation |
| Zhao and Chien (2015) | Min. sum of operation, waiting, in-vehicle & access costs | Bus & Trunk | Linear | Uniform distribution, many-to-one/one-to-many | Headway/ frequency & stop spacing | Near to optimal – Iterative algorithm |
| Gschwender et al. (2016) | Sum of operation, waiting, transfer & in-vehicle costs | Bus & Trunk | - | Deterministic demand | Fleet & Vehicle size | Numerical |
| Li et al. (2020) | Min. sum of operation, waiting, in-vehicle & access costs | Bus & Trunk | Rectangular | Pre-determined: Real case | Bus stop density & Headway/ frequency | Optimization Software |
| Yang et al. (2020) | Min. sum of line infrastructure, station, fleet, operation, waiting, in-vehicle & access costs | Bus & Trunk | Single corridor | Known, uniformly distributed continuous density function, Spatially heterogeneous, Many-to-one | Station spacing, line spacing & Headway/ frequency | Closed form & GA |

### 4.3. Transit network design

The most important part of any public transit system design is the network design part, where the objective is to design an optimal topology for the transit network. The network design problem is always assumed as a highly complicated problem, both in the modeling process and the solving process. Although most of the existing studies related to PBNDPs have applied mathematical programming and algorithms to deal with this problem, there are a number of studies that have used analytical models to investigate PBNDPs. The PBNDP can include different problems such as line spacing, stop spacing, etc.

Considering a single fixed-route public transit service, Kikuchi (1985) used analytical approaches to analyze the relationship between the number of stops and service frequency. While the fleet size was a known parameter, the number of stops and service frequency were the decision variables of the proposed model. As a key finding, the results showed that for any given fleet size, a unique solution can be obtained for the decision variables. Vaughan (1986) proposed an analytical model for jointly optimizing the route spacing and headway setting in a radial bus network (see Figure 20) with a many-to-many type of public transit service. They considered fleet size as a given parameter and assumed buses would always travel at a constant speed. The results showed that the optimal values of both of the decision variables are proportional to the inverse of the cube root of the ratio of the number of joining passengers to the number of passengers who are leaving the route.

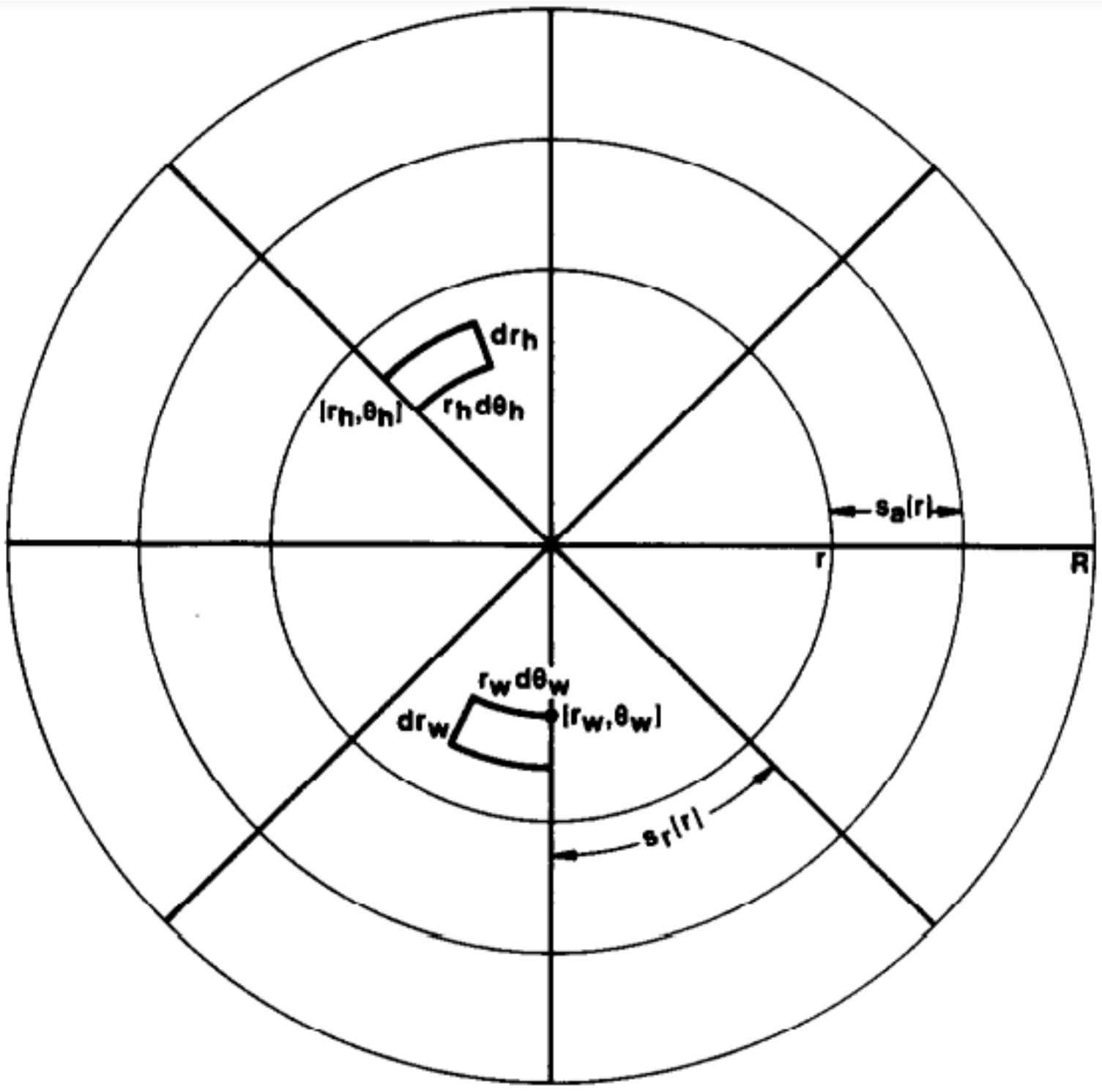

**Figure 20.** Radial bus network considered by Vaughan (1986)

Spasovic and Schonfeld (1993) considered a radial transit network between the suburbs and the CBD with a many-to-one transit service for a uniformly distributed demand along with the network as shown in Figure 21. Then, they proposed an analytical approach to find the optimal route length, route spacing, stop spacing, and service frequency for this network. The objective function was the summation of users' and the agency's costs. They also integrated the proposed model with an algorithm to solve the same problem for rectangular and wedge-shaped corridors when there is a vehicle capacity limitation. For this case, they assumed uniform and linearly decreasing demand patterns. As the main finding, they showed that the minimum total cost will be achieved when there is a balance between agency cost, access cost, and waiting cost. As a contribution to the literature, Chien and Schonfeld (1997) proposed an analytical model to design public transit systems for non-idealized geographical characteristics and demand distributions. The developed model was used to design an optimal grid public transit network with heterogeneous many-to-many demand and supply features by minimizing the summation of users' and agency costs. The decision variables were route and stop/station spacing, and service frequency. They also considered vehicle capacity and different routes costs for each zone in the optimization process. Almost the same problem, with a slight difference, has been investigated by Chien and Spasovic (2002). However, in that study, the optimal route and stop/station spacing, service frequency, and service fare were the decision variables while the objective was maximizing the summation of the agency's profit and social welfare.

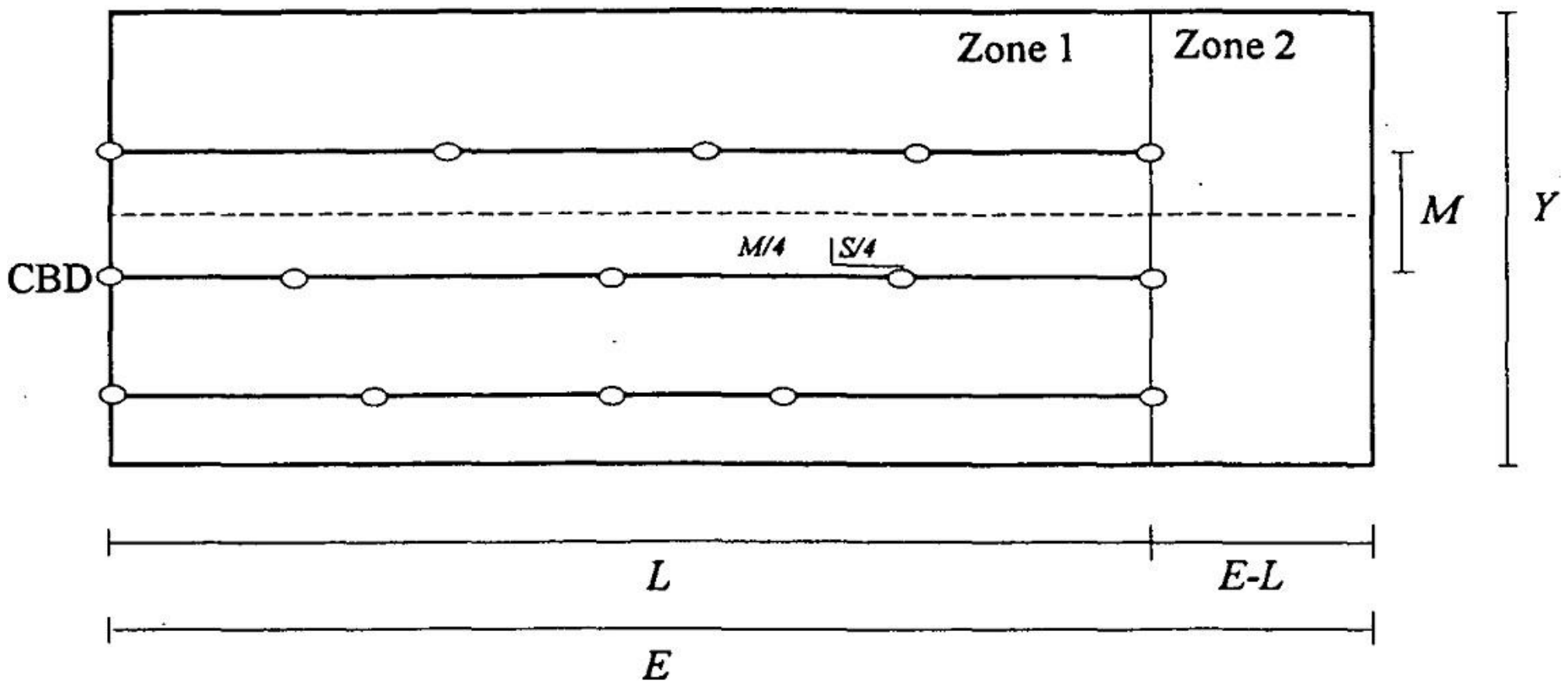

**Figure 21.** Transit network considered by Spasovic and Schonfeld (1993)

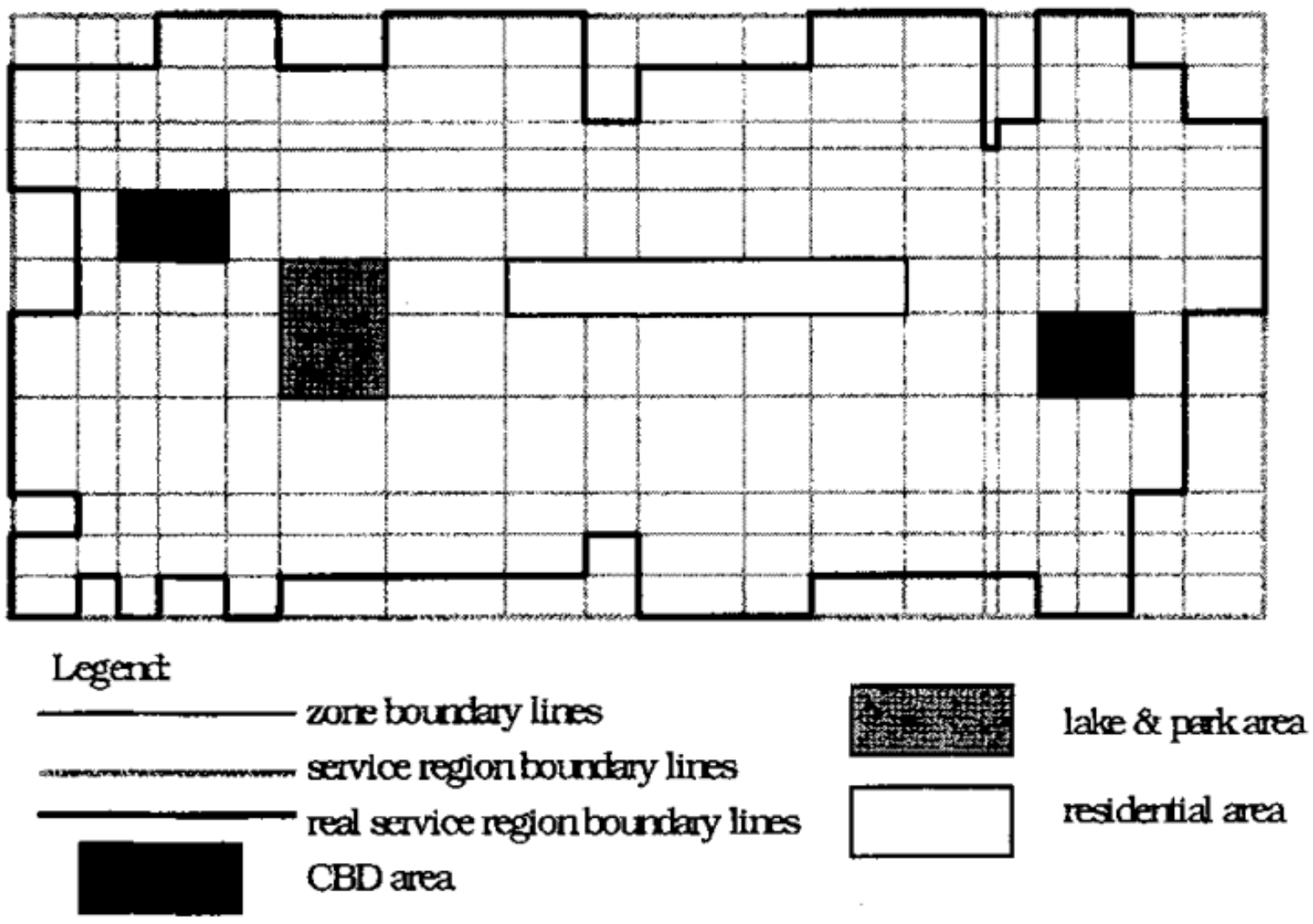

**Figure 22.** Irregular network structure considered by Chien and Schonfeld (1997)

Aldaihani et al. (2004) proposed an analytical framework to design a grid hybrid fixed- and flexible-route transit network as shown in Figure 23. They assumed that if the destination is in another zone, the flexible service will take users to a fixed-route stop, but if it is in the same zone, the flexible service will simply transfer users to their destination. Having this assumption, the objective of the proposed model was to find the optimal number of zones that minimizes the summation of users' costs, flexible-route service vehicle costs and fixed-route service costs. Ziari et al. (2007) proposed an analytical approach in order to study the problem of station/stop location problem with the objective of maximizing the accessibility of the users. Compared to the previous models, the proposed approach needed fewer variables to find optimal locations. Estrada et al. (2011) integrated an analytical approach and mathematical programming to propose a model for designing an efficient public transit network and applied it to the case of Barcelona, Spain, where the objective was to minimize the summation of the different cost types for both the agency and users. Line spacing, stop spacing, and service frequency were the main decision variables of the proposed model. Also, the study proposed a method to idealize the topology of a city to design a rectangular transit network. They considered the obtained optimal values for the idealized network as a design target when it comes to design a plan for the real network.

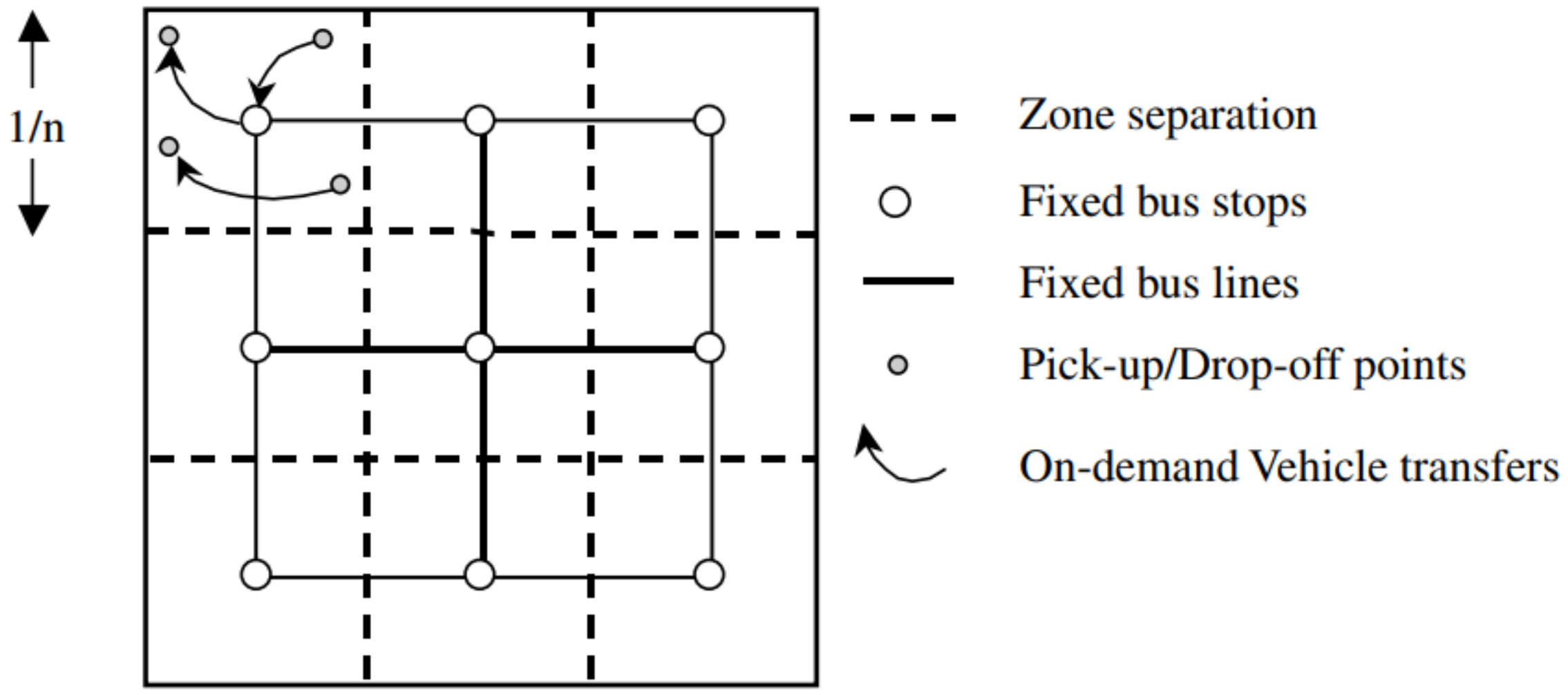

**Figure 23.** Hybrid network structure considered by Aldaihani et al. (2004)

Combining mathematical modeling and analytical approaches, Daganzo (2012) proposed a general method to design a public infrastructure system that serves multiclass users. The objective of the proposed model was to maximize the welfare of society. The model was used to obtain the service prices, the demand, and the general design of the system. Near to intersection stops can lead to queues in the network. Badia et al. (2014) proposed an analytical model for optimal transit network design in a city with a radial/circular network pattern in the core region and a hub and spoke pattern in the periphery area (see Figure 24). The objective was to find the optimal design by finding the optimal values for core area size, service frequency, line and stope spacings that minimize total passengers and operation costs. It was claimed that the proposed model is a useful approach to design a competitive public transportation network to serve passengers with high speeds and low headways. The study also analyzed the robustness of the proposed model against changes in the parameters. Tirachini (2014) used an analytical method to estimate the required number of stops along a public bus transit route as shown in Figure 25. To this end, first they started by reviewing the existing literature on bus stop location problems. Then using a case study, they showed that a Poisson model will lead to unrealistic results in estimating the probability of stopping and the number of bus stops for an on-demand service. Considering data of city of Sydney, to overcome this limitation of the Poisson model, they empirically derived two alternative models. In addition to the probability of having a stop in low demand zones, they also investigated the relationships between stop size, vehicle speed, stop spacing, and passengers in the zones with high demand. On the other hand, for the fixed-route system, they showed that vehicle speed, headway and dwell time, and demand pattern will highly affect the results of the bus location problem. Badia et al. (2016) used analytical approaches in order to conduct a comprehensive comparison between four possible public transit network structures (see Figure 26) to see which one fits better with a specific city. To this end, they considered idealized patterns for mobility with different concentration degrees, including very-concentrated cities, intermediate-concentrated cities, and decentralized cities, and showed that the selected network structure will be highly dependent on transport technology and geographical patterns of demand and mobility. Fan et al. (2018) considered a public transit network including an express transit service and local transit service that together served passengers' requests to travel in a city in different directions, as shown in Figure 27. They then proposed continuum approximations to design a grid bi-modal public transportation network. The objective was to find an optimal design that leads to the minimum total users' and operation costs. As one of the main findings, they showed that the best strategy will be jointly designing local-express transit network. Because the existing models, which solve this problem separately, won't consider all the advantages of bi-modal transit services.

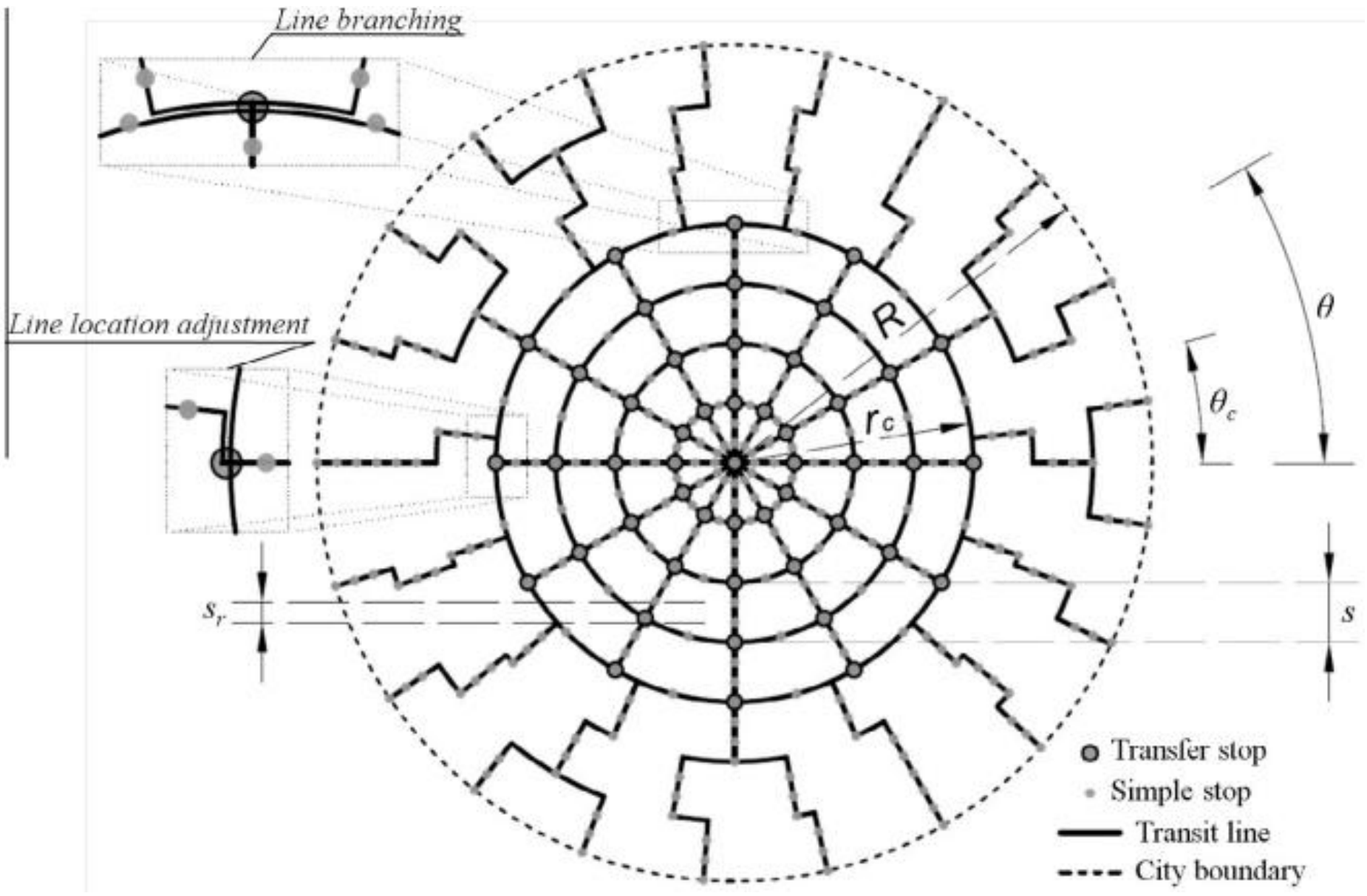

**Figure 24.** Network structure considered by Badia et al. (2014)

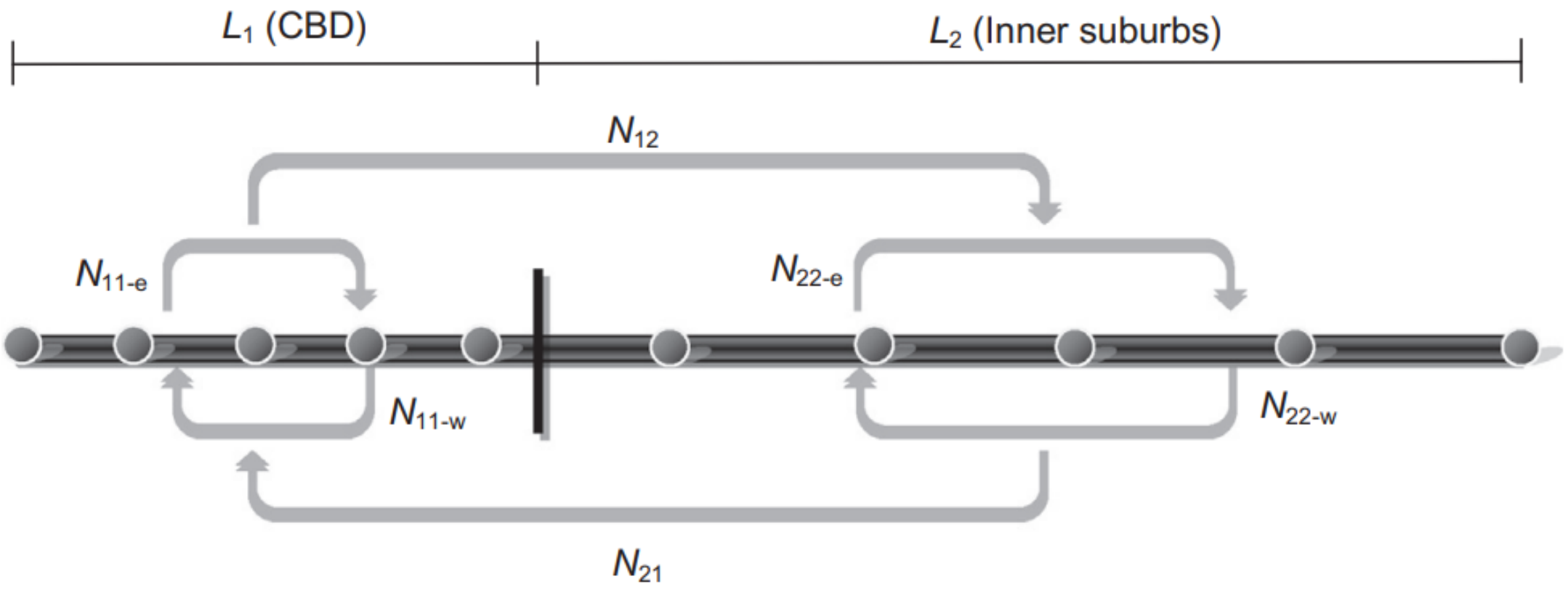

**Figure 25.** Considered configuration for a bus route and surrounding demand in Tirachini (2014)

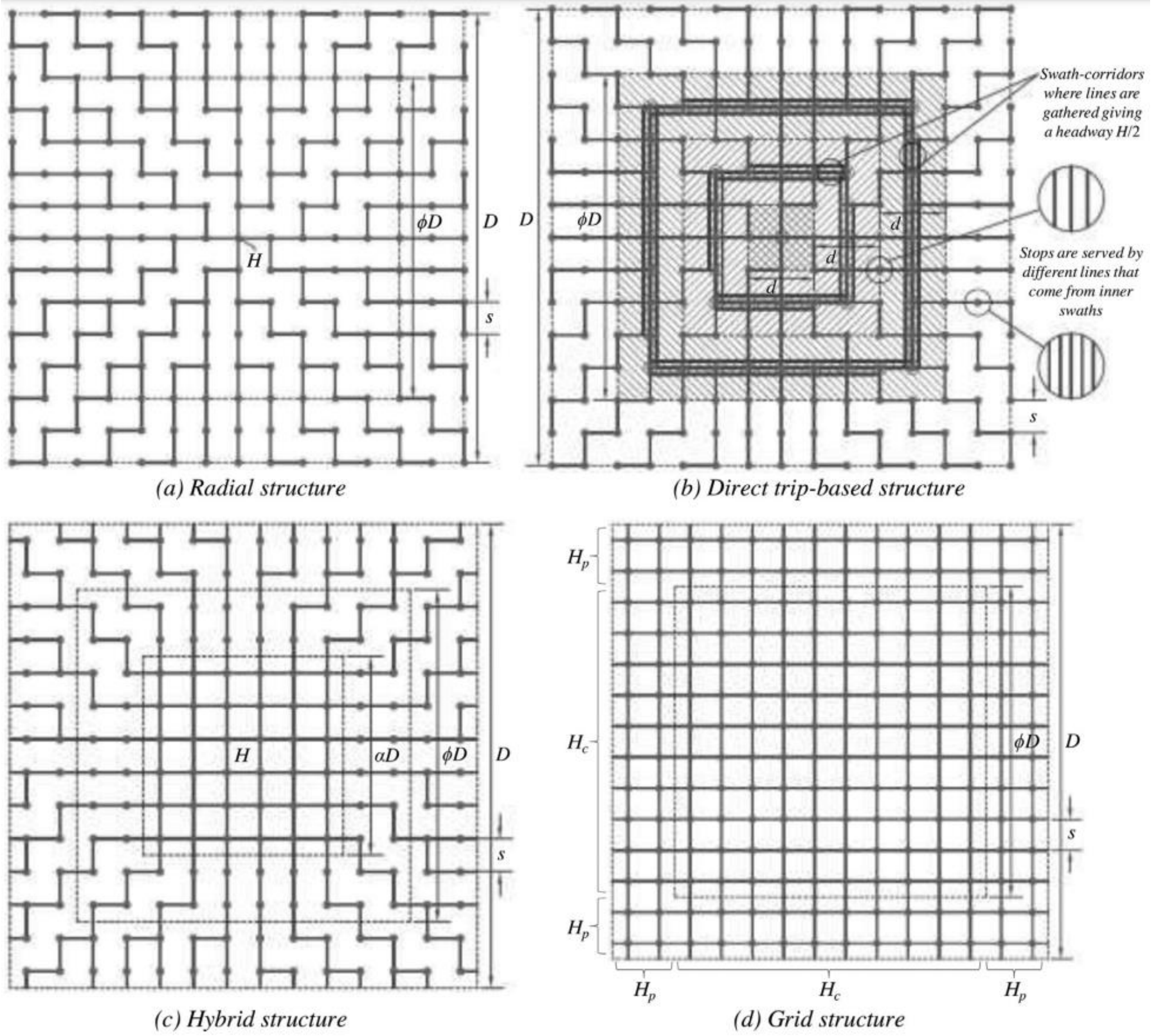

**Figure 26.** Four different structures considered in Badia et al. (2016)

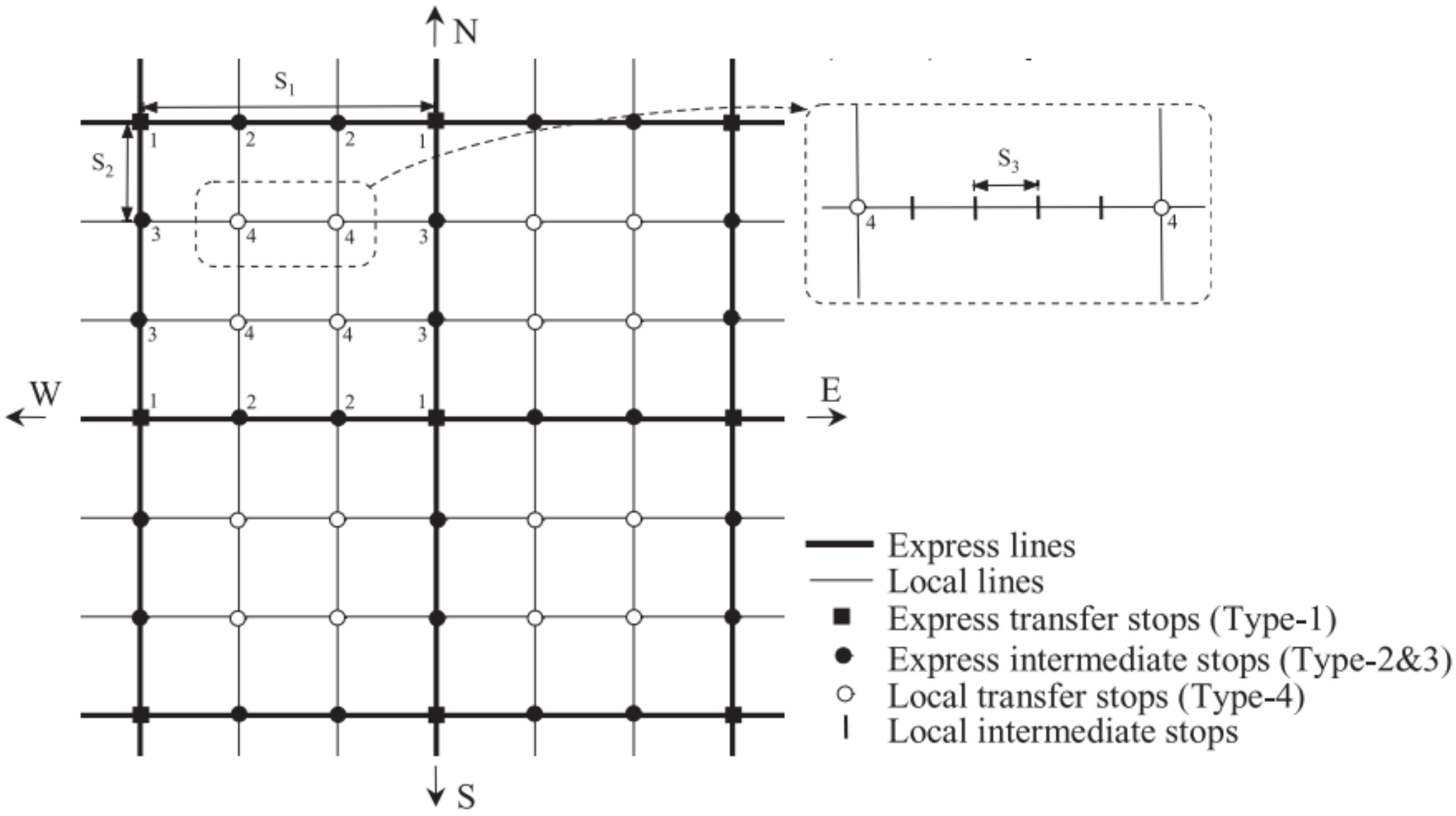

**Figure 27.** Considered hybrid express-local public transit network by Fan et al. (2018)

Considering a monocentric city (see Figure 28) and assuming that the majority of transit demand has focused on the central area, Luo and Nie (2019) proposed an approach to study the paired-line hybrid transit network design problem. They separated the city into two major areas: a central area, and the surrounding suburbs, and it was assumed that the distribution of the trips ending inside each area had a homogeneous pattern. Approximation methods have been developed to estimate cost functions under an exponentially decreasing demand pattern from the center to the suburban. Wu et al. (2020) integrated analytical methods and mathematical programming to design a public transit system that is being fed by shared bikes. They assumed that passengers could use bikes for both accessing the trunk system or reaching their destination after egressing from the truck service, or they could choose to walk to the stop/station or from the stop to their destination. Taking walking and biking as two feeder services into account, the optimal price for shared bikes has been obtained in order to have the system-optimal traffic assignment in the network. They considered a network structure with square grids for the central area and a radial structure for the suburbs (see Figure 29). The objective of the proposed model was to minimize total travel costs. They also proposed similar models for the case that walking or waking-fixed-route transit service are the alternatives to reach the trunk service. Considering a heterogeneous public transportation network as shown in Figure 30, Petit et al. (2021) proposed a combined analytical and mathematical programming approach to design an efficient urban bus transit network and obtain dynamic service frequency for a city with a heterogeneous demand pattern. They used continuum approximations to determine frequency setting and route spacing and used dynamic programming to solve the route and mode choice, and user equilibrium traffic assignment problems. For this study, they assumed that assigning exclusive lanes to the buses would reduce the reserved capacity for the private vehicles, and that as a result, congestion in the city would be highly dependent on the design of the transit network. Mei et al. (2021) considered an AB-type public transit service, as shown in Figure 31, and assumed that demand would slightly change along the route. Then, they proposed an analytical model to design an optimal skip-stop transit system. A heuristic algorithm was suggested to solve the proposed model. Solving different numerical examples, they showed that in most cases, the skip-stop service will have more efficient performance compared to the all-stop service.

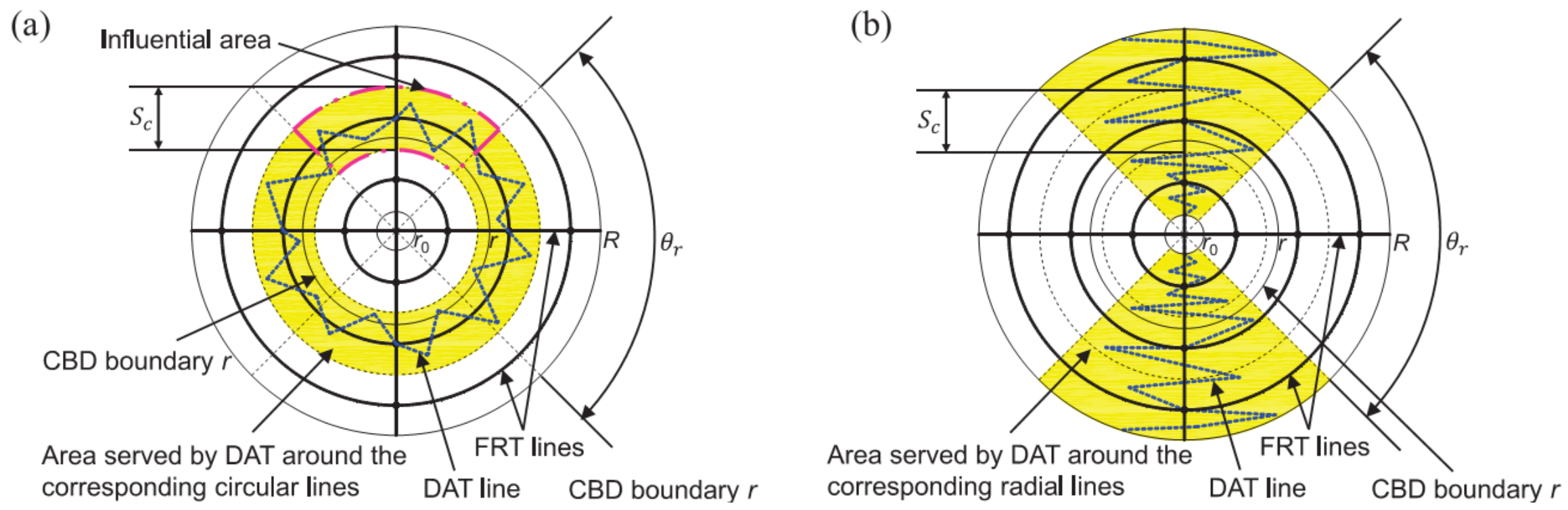

**Figure 28.** C-model (a) and R-model (b) hybrid public transit services considered Luo and Nie (2019)

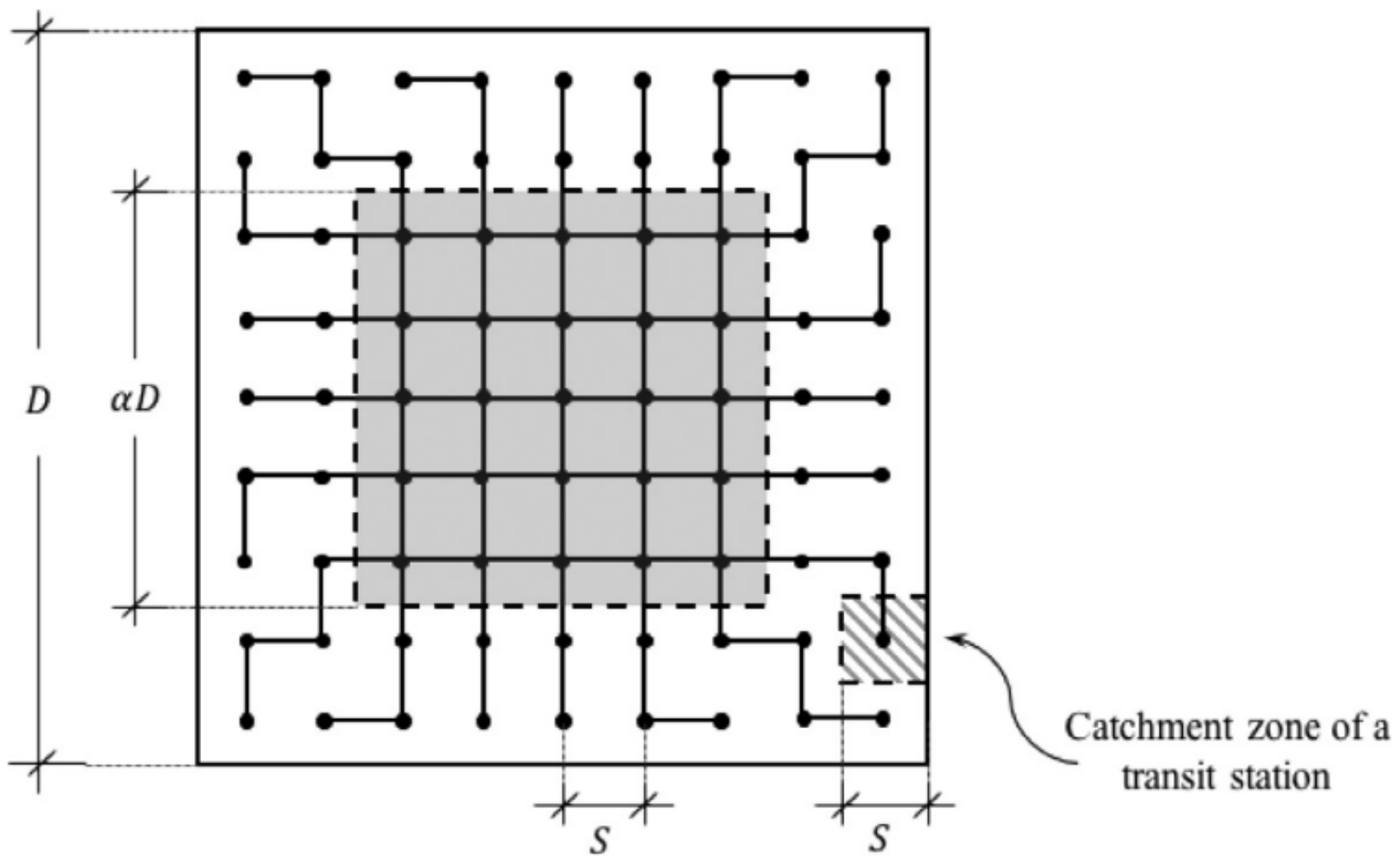

**Figure 29.** Considered network in Wu et al. (2020)

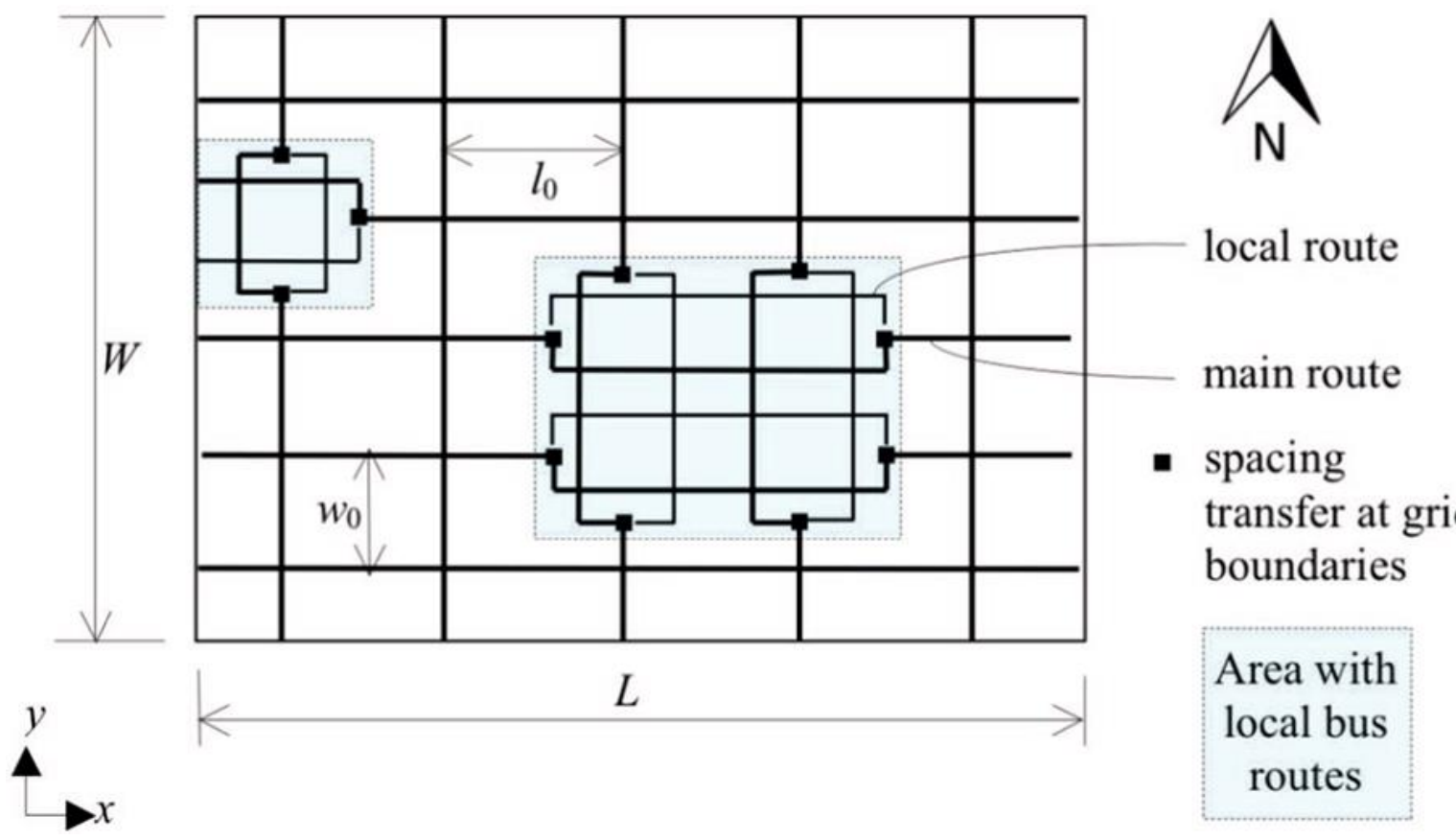

**Figure 30.** Heterogeneous public transportation network considered by Petit et al. (2021)

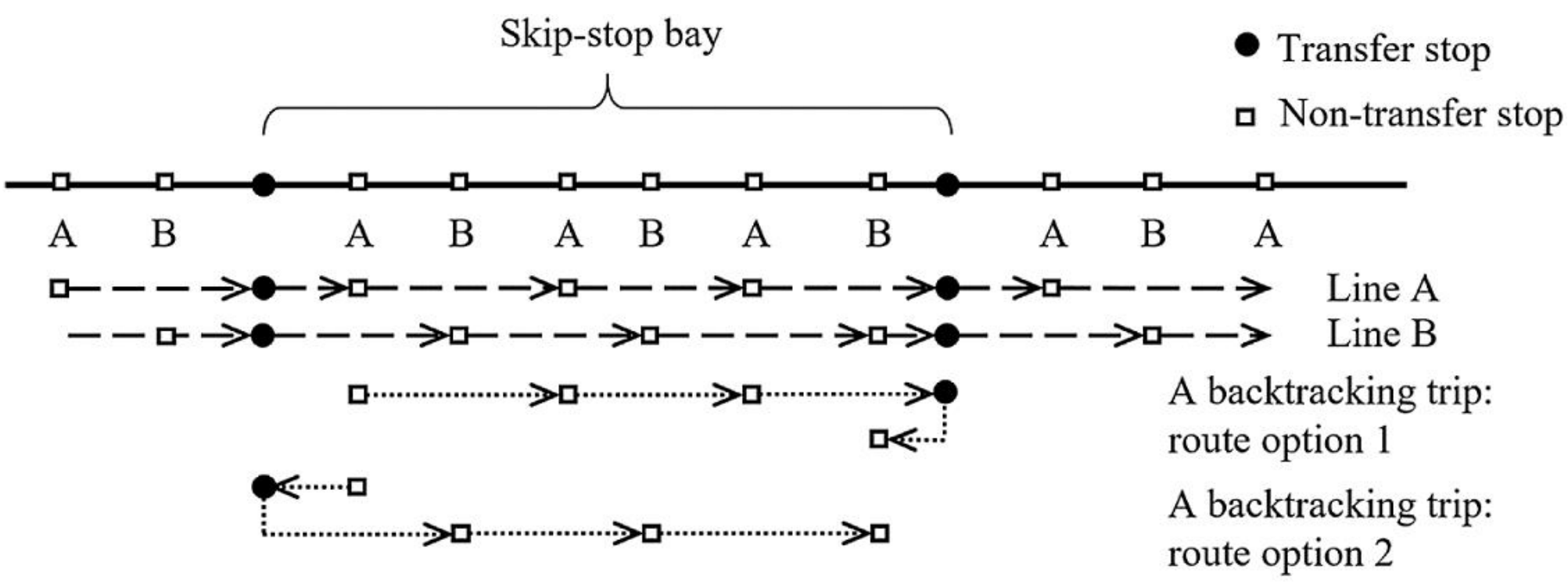

**Figure 31.** AB-type public transit service considered by Mei et al. (2021)

Table 9 represents the main features of the reviewed papers in this subsection. According to this table, the most commonly considered cost factors in this subcategory are access costs, waiting costs, in-vehicle costs, transfer costs, operation costs, and infrastructure costs. While some papers have focused on the rectangular grid networks, others have investigated the transit network design problems for the radial and ring networks. Although the uniform demand pattern is the most considered pattern in the literature, some researchers have conducted their studies based on random and location-dependent demand density functions. Headway, route spacing, stop spacing, route length, and service zone size are the main decision variables in most of the studies in this subcategory. As can be seen from the table, most studies were not able to find closed-form solutions for their decision variables. This is because of the high complexity of the problem of designing the topology of transit networks.

**Table 9.** Main features of the reviewed papers in subcategory 3.

| Reference | Objective Function | Transit mode | Network structure | Demand Pattern | Decision Variables | Solution Approach |
|---|---|---|---|---|---|---|
| Kikuchi (1985) | Min. sum access, waiting & in-vehicle times | General | Single fixed route | Uniform distribution | Headway/ Frequency & Number of stops | Numerical |
| Vaughan (1986) | Min. sum of walk, wait, in-vehicle & interchange times | Bus | Radial | Many-to-many | Headway/ Frequency & Route spacing | Closed form solution |
| Spasovic and Schonfeld (1993) | Min. sum operation, access & waiting costs | General | Radial, rectangular, wedge-shaped corridors | Uniform distribution, linearly decreasing, Many-to-one | Route length, route spacing, stop spacing & Headway/ frequency | Closed form solution, Algorithm |
| Chien and Schonfeld (1997) | Min. operation cost, infrastructure, transit line, waiting, access & in-vehicle costs | General | Rectangular | Random, uniformly distributed at each station | Route spacing, stop/station spacing & Headway/ frequency | Closed form solution |
| Chien and Spasovic (2002) | Max. total revenue minus operation cost | General | Rectangular | Heterogeneous many-to-many, sensitive to fare and travel time | Route spacing, stop/station spacing & Headway/ frequency | Closed form solution |
| Aldaihani et al. (2004) | Min. users', flexible-route service vehicle & fixed-route service costs | General | Rectangular | Constant rate | Number of zones | Numerical |
| Ziari et al. (2007) | Max. accessibility | Bus | - | - | Stop/station location | Closed form solution |
| Estrada et al. (2011) | Operation, vehicle, infrastructure, access, waiting, riding and transferring | General - Bus | Rectangular | Uniformly and independently distributed | Line spacing, stop spacing & frequency/ headway | Optimization software |
| Daganzo (2012) | Max. Social welfare | General | - | Fixed – elastic demand | Price, Demand & general design | Algorithm |
| Gu et al. (2013) | - | Bus | - | - | Stop location & holding strategy | Numerical |
| Badia et al. (2014) | Min. sum operation, infrastructure, fleet, access, waiting, in-vehicle & transfer costs | General | Radial – hub & spoke | Uniform distribution | Core area size, line spacing, stop spacing, headway/ frequency | Grid search |
| Tirachini (2014) | - | Bus | Single route | Uniform Demand – Poisson model | Number of stops | Numerical |
| Badia et al. (2016) | Min. sum infrastructure, operation, access, waiting, in-vehicle & transfer costs | General - Bus | Rectangular - Radial | Different scenarios | Headway/ frequency, stop spacing & swath-corridor spacing | Software |
| Fan et al. (2018) | Min. sum infrastructure, stop, operation, access, waiting, in-vehicle & transfer costs | Bus, BRT, Rail | Rectangular grid network | Uniform distribution | Headway/ frequency & line spacing | Iterative algorithm |
| Luo and Nie (2019) | Min. sum fleet, vehicle operation, access, waiting, in-vehicle & transfer costs | Bus | Radial-and-ring | Exponential Demand | Headway/ frequency & line spacing | Numerical & simulation |
| Wu et al. (2020)* | Min. average agency and users cost | General, Bike | Radial & square grid | Average density function | Headway/ frequency & line spacing | Software |
| Petit et al. (2021) | Operation, access, waiting, in-vehicle & transfer costs | Bus | Rectangular | Spatial and temporal heterogeneous demand | Headway/ frequency, route spacing, route choice, model choice, traffic assignment | Closed form solution & Iterative algorithm |
| Mei et al. (2021) | Min. sum infrastructure, vehicle, operation, access, waiting, in-vehicle & transfer costs | Bus & Rail | Linear | Spatial heterogeneous demand | Headway/ frequency, stop location & route plan | Heuristic algorithm |

*: The paper has considered an emerging technology.

### 4.4. Fleet size, vehicle capacity and crowding

Other critical decisions that managers of the public transportation systems must make are related to fleet size and vehicle capacity, which will be highly impacted by demand patterns and crowding levels in the network. Fleet size and vehicle capacity are two important parameters that significantly affect the performance of the network, users' costs and operation cost. Therefore, there are lots of studies that have applied various methods to study problems related to fleet size and vehicle capacity in public bus transit services.

Diana et al. (2006) developed analytical approaches in order to investigate the fleet size problem for a flexible-route transit service when there is a maximum allowed detour and maximum waiting time at each stop considered for this service. The demand distribution and service level were the only requirements of the developed probabilistic model, and it was shown that this model performs better than simulation-based approaches when there is no access to detailed data. Cortés et al. (2011) applied analytical approaches to suggest an efficient strategy for obtaining the optimal values of service headways for low- and high-demand areas, and vehicle capacities for a single public transit line. To this end, they integrated short turning and deadheading in the modeling process. As a key finding, they showed that the short turning strategy will have a bigger impact on the total cost compared to deadheading. Considering a single line transit service, Jara-Díaz et al. (2017) applied analytical equations to investigate the fleet size and bus capacity problems by minimizing total users' and agency's costs. It was assumed that travel demand and travel lengths would be different for peak and off-peak times. The results showed that considering both peak and off-peak times would lead to different optimal values for the decision variables compared to the case where only peak time was considered. When both periods have been considered, a larger fleet of smaller buses will serve passengers compared to the peak hour only case. In addition, passengers would be served with smaller headways during peak hours and larger headways during off-peak times compared to headways obtained by solving the model for each period separately. Also, the optimal vehicle capacity would be a value between the optimal values obtained for peak and off-peak periods separately. Considering a heterogeneous demand pattern which varies by time and location, Hörcher and Graham (2018) formulated the multi-period public transport supply problem analytically. They defined the passengers' cost function as the sum of waiting and crowding costs, and assumed that this cost function is a heterogeneous function in order to produce more realistic results. Sarango et al. (2019) claimed that service quality in urban transit networks is directly related to the crowding level in the vehicle, and then proposed a model to estimate the in-vehicle crowding level. The model was developed based on actual data and its results were used to improve the service quality. They also investigated the possibility of applying analytical approaches to estimate in-vehicle crowding levels in different passenger traffic situations. Criticizing the existing models for obtaining fleet size for public transit systems, Liu (2020) proposed a model based on a deficit function that was able to overcome the limitations of previous models. As another contribution, they used continuum approximations to make the proposed model differentiable at its significant points. Considering a time-dependent demand pattern, Zhang et al. (2021) proposed a two-stage dynamic heuristic method to study the problem of fleet management. The suggested method first finds the optimal schedule for vehicle dispatching, then, obtains the optimal schedule for vehicle purchasing and retiring.

Table 10 represents the main features of the reviewed papers in this subsection. According to this table, the most commonly considered cost factors in this subcategory are waiting costs, in-vehicle costs, crowding costs, and operation costs. Most of the papers have considered a single transit line and tried to optimize headway, vehicle size, and fleet size by considering a realistic demand pattern. Except for one of the reviewed papers, others were not able to find closed-form solutions for their decision variables.

**Table 10.** Main features of the reviewed papers in subcategory 4.

| Reference | Objective Function | Transit mode | Network structure | Demand Pattern | Decision Variables | Solution Approach |
|---|---|---|---|---|---|---|
| Diana et al. (2006) | - | General | Square area | Discrete distribution | Fleet size | Numerical-software |
| Cortés et al. (2011) | Min. sum of operation, waiting & in-vehicle costs | Bus | Single public transit line | Spatially disaggregated demand information | Headway/ frequency & vehicle capacity | Closed form solution |
| Jara-Díaz et al. (2017) | Min. sum of operation, waiting & in-vehicle costs | Bus | Single public transit line | Given parameter, different for peak and off-peak periods | Headway/ frequency, vehicle capacity & Fleet size | Numerical |
| Hörcher and Graham (2018) | Min. sum of operation, waiting, in-vehicle & crowding | General | Single line | Exogenous parameter & fluctuating elastic demand | vehicle capacity | Numerical & simulation |
| Sarango et al. (2019) | - | General | Single line | Real data | Crowding level | - |
| Liu (2020) | - | General | - | - | Fleet size | Numerical |
| Zhang et al. (2021) | Min. sum of fleet, operation and maintenance & salvage value | General | - | Time varying demand | Number of vehicle purchased in each time, type | Heuristic algorithm – Tabu search |

| | | |
|---|---|---|
| | | of each vehicle, operation and retirement schedules of each vehicle |

### 4.5. Routing

The vehicle routing problem has always been considered an important but highly complicated problem in the public transportation industry. A few studies have applied analytical methods to deal with routing problems in PBTNDPs and operations planning. Daganzo (1984) conducted one of the first studies in this area. They applied analytical approaches to make a comparison between the performance of the fixed-route bus transit systems and a dial-a-ride transit service. They considered a simple routing strategy in the network in order to make the comparison process easier. I-Jy Chien (2005) integrated analytical approaches and numerical techniques to design an optimal bus feeder transit service that served a shuttle transit system. They assumed that the shuttle system carried passengers between two major hubs. The proposed methods were used to obtain optimal dispatching policy, vehicle capacity, and route choice for the feeder bus system. Considering a conventional fixed-route public bus transit system, Wirasinghe and Vandebona (2011) developed an analytical method to solve the problem of the path chosen for public bus services. They claimed that despite having different planning time frames, it was not possible to choose an efficient transit route without considering other parameters such as fleet size and headways. The total cost function was considered as the summation of the users' access cost, the agency's operation cost, waiting cost and traveling cost. Also, they took both grid and non-grid transit network structures into account. The proposed model started with an initial routing solution and then updated the solution step by step. Ellegood et al. (2015) focused on the routing problem for school buses. They developed continuous approximation models to deal with the school bus routing problem, assuming that a mixed load strategy was allowed in the network. A mixed load strategy means a bus route will serve multiple schools at the same time. As a key finding, they showed that when the service zone is large enough, schools are close to each other and stops are serving multiple schools' students, a mixed load strategy will be the optimal strategy.

The passenger route choice problem is another important subject, where the problem is to understand how (based which criteria) passengers choose their route when they are using public transit services. Luo and Nie (2020) applied some sketched design models to investigate route choice behavior in public transit systems. A continuous approximation approach was used to analyze three competitive public transportation systems and to find how the route choice behavior of the passengers affected the operations of the transit system. In addition, for each transit system, they proposed some analytical models to obtain the aggregate share of all routes based on the users' utility.

Table 11 represents the main features of the reviewed papers in this subsection. According to this table, different network structures have been considered in the literature to solve the routing problem using analytical approaches. However, most of them have considered operation costs, waiting costs, access costs, and in-vehicle costs as the cost factors in their cost function. Location-dependent demand functions are considered in the reviewed paper and except for one of them, others were not able to find closed form solutions. The Routing Problem is a high-complexity problem related to the PBTNDP&OP.

**Table 11.** Main features of the reviewed papers in subcategory 5.

| Reference | Objective Function | Transit mode | Network structure | Demand Pattern | Decision Variables | Solution Approach |
|---|---|---|---|---|---|---|
| Daganzo (1984) | Sum. operation, waiting, access & in-vehicle costs | Bus | General | Exogenous parameter | Routing strategy | Approximation |
| I-Jy Chien (2005) | Min. sum operation, waiting, access & in-vehicle costs | Bus-rail | Single feeder route | Exogenous parameter | Headway/ frequency, vehicle size & routing strategy | Closed form – Iterative algorithm |
| Wirasinghe and Vandebona (2011) | Min. sum operation, waiting, access & in-vehicle costs | Bus | Grid & non-grid | Average passenger demand per day, different for each generator | Routing strategy | Iterative algorithm – software |
| Ellegood et al. (2015) | Expected distance | Bus | Non-idealized | - | Routing strategy | Numerical |
| Luo and Nie (2020) | Min. sum operation, infrastructure, fleet, access, waiting, in-vehicle & transfer costs | Bus | Radial, rectangular-grid, hybrid | Spatially homogeneous demand | Routing strategy | Numerical-Optimization software |

### 4.6. Rapid transit and shuttle transit system planning

Bus Rapid Transit (BRT) and shuttle transit systems are some of the most popular public bus transit modes in most populated cities. There are a few studies that have used analytical methods to investigate different problems related to these systems. Sim and Templeton (1983) used an analytical model to study the shuttle system operation planning problem by obtaining optimal fleet size and service frequency. They considered a shuttle transit system where the passengers are served using vehicles with unlimited capacity. The vehicles transfer passengers from a major hub to different destinations and take them back to the major hub. They also assumed that trip times are independent and follow an exponential distribution pattern, the users' arrival times follow a Poisson process, and the length of the users' queue and their waiting time can be calculated in accordance with steady-state conditions. Thilakaratne and Wirasinghe (2016) compared the performances of the traditional bus transit systems and BRT systems in order to find the advantages of transitioning from the first one to the second one. To this end, they considered a corridor that a long classic bus route serves the users, a BRT system with a limited number of stops that serves passengers along a part of the corridor, and a classic bus system that also serves passengers along the whole corridor. Considering a total cost function including waiting and riding time, bus operation costs, and capital costs, they proposed a model to obtain the optimal values of the service frequency for both services, and BRT terminals and stop locations. Dantas et al. (2019) used analytical equations to conduct a comprehensive sensitivity analysis on the performance of the BRT services with respect to changes in the critical parameters. The goal was to understand how the system works and how its performance can be improved, especially by analyzing factors impacting average system size, average queue length, and average time of queue. As a key result, the study showed that service headway will have the most significant effect on all aspects of the system. Deng and Zhu (2020) developed an analytical model to design a hybrid feeder-express bus transit service, where the feeder system transfers passengers to the express bus stops and the express bus transfers passengers between two terminals at a higher speed compared to the conventional local bus services. They applied the proposed approach to compare the passengers' cost in regular bus services and express bus services, and showed that this cost in the express system is significantly lower than the cost in the regular service.

Table 12 represents the main features of the reviewed papers in this subsection. According to this table, a single transit corridor and a general transit network have been considered by different authors. Except Thilakaratne and Wirasinghe (2016), other studies were not able to find closed form solutions for their decision variables. Instead, they have solved numerical examples to investigate the performance of the proposed models. Operation costs, capital costs, and waiting costs are the main considered cost factors in this section.

**Table 12.** Main features of the reviewed papers in subcategory 6.

| Reference | Objective Function | Transit mode | Network structure | Demand Pattern | Decision Variables | Solution Approach |
|---|---|---|---|---|---|---|
| Sim and Templeton (1983) | Min. sum of operation, fleet & waiting costs | Bus | General | Poisson process | Headway/ frequency & fleet size | Numerical |
| Thilakaratne and Wirasinghe (2016) | Min. sum of operation, capital, in-vehicle & waiting costs | Bus | Single corridor | Continuous smooth function | Headway/ frequency, terminals & stop location | Closed form solution |
| Dantas et al. (2019) | - | Bus | - | - | Average system size, average queue length, & average time of queue | Numerical |
| Deng and Zhu (2020) | Sum of in-vehicle & waiting costs | Bus | Single corridor | Uniform distribution – different for each zone | - | Numerical |

## 4.7. Travel time reliability

Service and travel time reliability in public transit is one of the most important factors for both users and agencies. Unreliable public transit services will make the users disappointed and encourage them to use more private transit modes. In any case, unreliable transit services will lead to higher users' and agency costs. Because of this importance, there are a significant number of studies that have focused on the service and travel time reliability in public transportation systems. However, only a few studies have applied analytical methods to study this problem in public bus transit systems.

Knoppers and Muller (1995) was one of the first papers that used analytical methods to study the travel time fluctuation in public transit. They proposed a model to optimize the transfer time for the users by considering the uncertainty in the actual users' arrival time. To obtain the optimal transfer time, the study minimized the users' expected waiting time. The results indicated that the optimal transfer time can be obtained only in the case of having users' arrival times fluctuate within specific time slots. Fonzone et al. (2015) studied the bus bunching problem with respect to the passengers' reactions to the system's reliability. To this end, they used a logit model to model the users' choice for their arrival time at the boarding points by considering the possibility of missing a bus, then, obtained the expected waiting time of the users and the dwell time of the consecutive buses. They showed that the insufficient boarding rate and non-uniform arrival patterns are the main reasons for the bus bunching. The findings suggested considering the impacts of non-uniform arrival patterns on bus bunching when it comes to determining service control measures in public bus transit service. To design a reliable bus transit service, Lai et al. (2020) developed a model for planning a resilient schedule, considering delays caused by priority signals in intersections or/and passengers boarding in the stops, and travel times on each link in the transit corridor. The objective of the proposed model was maximizing service resiliency according to the realistic volume of the passenger flow, and the results showed this target will be achieved when there is a balance between users' and the agency's cost. Gao et al. (2021) applied analytical models to analyze the reliability of public bus transit services during the COVID-19 pandemic in China. They tried to study the problem of having reliable service while preventing the spread of the virus through public transportation services. Taking the constraints related to the buses into account, they developed a timetable design model, and then proposed a reliability-based method to analyze the timetable systems. Finally, to obtain the failure rate of a timetable, the study presented an event tree based model. Failure was defined as the failure to reach the downstream stop on-time, and the failure rate was obtained based on the upper and lower bounds of failure probability.

Table 13 represents the main features of the reviewed papers in this subsection. According to this table, service frequency, transfer time, expected waiting time, and dwell time are some of the decision variables considered in this section. However, only Lai et al. (2020) were able to find an optimal closed form solution in this sub-category. It was predictable, because most of the reliability related problems are time-dependent and dynamic problems. Hence, it would be a great job to find closed-form solutions for these types of problems using analytical approaches. As can be seen from this table, in-stop crowding, in-vehicle crowding, and traffic are important factors that have not been considered by any study.

**Table 13.** Main features of the reviewed papers in subcategory 7.

| Reference | Objective Function | Transit mode | Network structure | Demand Pattern | Decision Variables | Solution Approach |
|---|---|---|---|---|---|---|
| Knoppers and Muller (1995) | Min. transfer waiting time | General | General | - | Transfer time | - |
| Fonzone et al. (2015) | - | Bus | General | Deterministic time dependent & time and location dependent probability distribution function | Expected waiting time & dwell time of the consecutive buses | Numerical-algorithm |
| Lai et al. (2020)* | Min. sum of operation, waiting transfer, delay & missed connection costs | Bus | Single corridor | Uniform demand different for each | Headway/ frequency | Closed form solution |
| Gao et al. (2021) | - | Bus | Single route & General | - | Headway/ frequency & failure rate | Numerical |

*: The paper has considered an emerging technology.

### 4.8. Bus queue, bus and transit signal priority

Giving priority to public buses, especially in congested cities and at signalized junctions, has been accepted as an efficient strategy in order to make the public transportation services more productive. Hence, another application of analytical approaches in PBTNDP&OP is bus and transit signal priority. Although the topic itself is a hot topic in the public transit planning field, only a few articles have applied analytical methods to study this problem. Wu and Hounsell (1998) was the oldest identified study in this sub-category. Considering three different pre-signal categories, they tried to provide analytical approaches to deal with the problem of designing and operating pre-signals. Pre-signals are being installed close to the end of the bus lanes in order to give priority to buses with access to the downstream junction. They first categorized pre-signals with respect to their different requirements for signaling, operational features, the types of effects that they had on the network's capacity, and delays in the travel time. Then, for each category, they developed analytical models in order to perform pre-implementation evaluation of pre-signals, by taking the design, capacity, and signal adjustment problems into account. The proposed approach was able to estimate the delay or the saved time for both non-priority and priority traffic and to predict traffic queue relocation at pre-signal areas. In addition, the bus advance area is designed using analytical models. Considering curbside bus stops, Gu et al. (2015) applied analytical approaches to present steady-state queueing models in stops. They assumed that stops belonged to multiple bus routes and that they were giving services independent of the traffic signals and other stops. The objective of the proposed model was to predict the maximum possible number of buses that can service a stop at the same time without going beyond the considered average bus delay time. Obtaining the number of bus berths for each stop under a known bus demand and average delay time was another application of the proposed model. Taking both the benefits at major streets and delays on side streets into account, Mishra et al. (2020) developed an analytical approach to making a passenger-based priority for public buses in the city. The main idea of their approach was to make a balance between the benefits and delays. They considered three factors in order to assign priority to the scheduled buses using transit signals: adherence to the schedule, waiting time for the passengers at the stops, and passenger occupancy.

Table 14 represents the main features of the reviewed papers in this subsection.

**Table 14.** Main features of the reviewed papers in subcategory 8.

| Reference | Objective Function | Transit mode | Network structure | Demand Pattern | Decision Variables | Solution Approach |
|---|---|---|---|---|---|---|
| Wu and Hounsell (1998)* | - | Bus | Single route – multiple lanes | - | Bus advance area, traffic queue, signal times, delay & saved times | Closed form solution |
| Gu et al. (2015) | - | Bus | Multiple routes | - | Maximum number of serviceable buses at a stop, bus berths for each stops & delays | Closed form solution - approximation |
| Mishra et al. (2020)* | - | Bus | Single route | - | Priority assigning & minimum number of bus passengers | Simulation |

*: The paper has considered an emerging technology.

### 4.9. Bus holding strategies

Holding strategies can be considered as a dispatching policy; however, since the number of the identified articles that have applied analytical methods to study holding policies was large enough, these articles are categorized in a different sub-category titled "bus holding strategies".

Considering a stochastic transit service model and by developing an analytical approach, Hickman (2001) tried to find the best bus holding time at a control stop. The idea of taking the stochasticity effects into account in the holding process was to make the problem close to real world conditions. Both gradient or line search techniques were easily applicable to solve the proposed holding model for a bus as a convex quadratic model. Zhao et al. (2003) developed an approach for obtaining optimal holding policies for each bus at each stop. To this end, they considered both the bus and the stop as agents assuming that they are able to communicate/negotiate with each other in real-time to determine holding time based on a total cost function. Finally, to verify the validity of the proposed algorithm, they made a comparison between the suggested real-time scenarios, on-schedule, and even-headway scenarios. In order to provide reliable public transit services, Zhao et al. (2006) proposed a model to obtain optimal holding time for the buses when it comes to obtaining headways for dispatching the buses. The objective of the model was to minimize the expected waiting cost of the users, and the results confirmed that the services would be reliable if the managers considered the slack time in the scheduling process. To estimate the travel/delay times for a transit route with multiple routes, they used approximation methods and confirmed the results by comparing them with simulation results. To overcome the limitations of headway-based dynamic holding algorithms, Xuan et al. (2011) developed new dynamic holding models based on the actual and scheduled arrival time differences. The proposed model was able to obtain the closed form solution just in the case of having a one parameter model. Dai et al. (2019) integrated analytical approaches, game theory, and simulation methods to design a framework for headway-based bus holding. The method included a dynamic algorithm to evaluate the control points, rank them, and suggest the best one. They considered the effects of each control point separately and the collective impacts of all control points, bus travel times, and changes in the demands in the modeling process. Using analytical equations and simulation approaches, Chen et al. (2021) developed a method for solving the bus holding problems. They considered the gap between actual waiting time and perceived waiting time of the users in the expected waiting cost function and tried to find a solution to deal with the bus holding problem for bunching buses and speed adjustment problem for lagging buses. For the first problem, they applied a threshold method, while for the second one, a Markovian decision model was used. As the key findings, the study showed that (a) a rigid holding control can lead to more regular headways, (b) increased perceived users' waiting time, the demand patterns and vehicle occupancy will impact the results, and (c) suitable speed adjustment won't negatively affect the average cruising speed in a trip.

Table 15 represents the main features of the reviewed papers in this subsection. According to this table, minimization of the expected waiting time and slack time are two commonly used objective functions in this subcategory. Compared to other subcategories, in this subcategory, most studies have considered more realistic demand patterns and were able to find closed form solutions. "Control Point" and "slack time" are the two main decision variables in this section.

**Table 15.** Main features of the reviewed papers in subcategory 9.

| Reference | Objective Function | Transit mode | Network structure | Demand Pattern | Decision Variables | Solution Approach |
|---|---|---|---|---|---|---|
| Hickman (2001) | Min. total expected waiting time | General & Bus | Single bus route | Poisson process (different for each stop) & binomial distribution | Load on vehicle, running time, waiting time & headway | Numerical & Programming |
| Zhao et al. (2003) | Min. sum of off-bus and on-bus waiting costs | Bus | Single bus loop | Spatial and temporal heterogeneous demand | Holding strategy & dispatching policy | Closed form solution & algorithm |
| Zhao et al. (2006) | Min. passengers' expected waiting time | Bus | Single bus loop | A probability density function | Slack time | Closed form solution, approximation & algorithm |
| Xuan et al. (2011) | Min. slack times | Bus | General | Stationary (different for each stop) | Control coefficient & control point | Closed form solution & algorithm, software |
| Dai et al. (2019) | - | Bus | General | Dynamic demand-stop-level time-dependent function-real GPS/smart card data | control point | Dynamic solution & simulation |
| Chen et al. (2021) | - | Bus | Single bus route | Stationary (different for each stop) | Holding strategy & speed adjustment | Algorithm & Simulation |

### 4.10. Emerging technologies in public bus transit systems

For this part, please see Mahmoudi et al. (2024).

### 4.11. Macroscopic Fundamental Diagram and traffic assignment

All types of issues related to the Macroscopic Fundamental Diagram (MFD) and Traffic Assignment (TA) have always been among the most popular but highly challenging issues in transportation studies. There are some studies that have applied analytical tools to study these problems, considering bus transit services as one of the transit modes in the network. Guler and Cassidy (2012) investigated the delay problem in bottlenecks when buses and other vehicles share the lanes in a bottleneck. The main idea for their model was to force vehicles to take shared lanes in the bottlenecks to make more available space in the unshared lanes until adding a new vehicle causes a delay in the bus travel time. Using analytical methods, they showed that the lane sharing strategy can increase the capacity of the bottleneck compared to the available capacity in the exclusive lane assigning strategy. For a bi-modal urban transportation network where public buses share the same roads with other cars, Geroliminis et al. (2014) conducted a study to figure out how it is possible to evaluate the performance of the network using aggregated relationships and how the interactions between buses and cars affect the quality of the network's performance. They also analyzed the impacts of the bus stops on the network's performance. They proposed an analytical model to estimate passenger flow and used simulation results to develop a 3D MFD related to the vehicle flow in the considered bi-modal network and a 3D MFD for the passengers. Using a clustering approach, they divided the network into a set of zones. Then they demonstrated that the level of interactions and their effects vary from zone to zone. Taking the same bi-modal network, Liu et al. (2015) developed analytical methods in order to modify the cell transition model to consider both modes to better describe the moving bottlenecks caused by having buses in the network. They considered capacity reduction and speed differentials in their model and applied the proposed approach to solve the mixed traffic assignment problem with the objective of minimizing total users' cost in both classes. Another similar study has been conducted by Loder et al. (2019). They proposed an approach to develop a functional MFD for multiclass networks where both cars and buses serve users in the network. Their model had two components: first, a method to restrict the number of travels for any given vehicle number in each class. They used topological and operational specifications of the network to estimate the parameters of this method. Second, based on the topological features of the network and real traffic data, a smoothing parameter was estimated to understand how vehicles interactions can decrease the number of travels.

Similar to previous studies, Dakic et al. (2020) considered bi-class traffic flow in a transit network and tried to develop a passenger MFD. As a contribution, they considered the stochasticity in the moving bus bottleneck and modified variational theory by taking the effects of this stochasticity into account. To this end, they developed a probabilistic variational theory graph, and also considered the correlation between different bus arrival times, the dynamic passengers' mode choice behavior, and the impacts of the different possible traffic conditions on the performance of the public transit service. Their analysis showed that the proposed model will be highly helpful when the public bus transit system provides short and variable headway services. Dakic et al. (2021) analyzed the impacts of users' behavior and trip length on different cost items in the network, the impacts of the topology of the public bus transit network on the private vehicle traffic flow, the impacts of the traffic conditions on the mode choice function of the users, and the interactions among multiple transit classes. In particular, they used 3D MFD to take the interactions and speed of each class into account in the modeling process. Then, they used analytical methods to solve the PBTNDP for free-flow and congested network cases, considering more effective factors. As the main finding, they showed that considering uniform trip patterns would lead to incorrect fleet size, resulting in crowding at the bus stops at peak times. To overcome this problem, they suggested considering actual trip length distribution.

Table 16 represents the main features of the reviewed papers in this subsection.

**Table 16.** Main features of the reviewed papers in subcategory 11.

| Reference | Objective Function | Transit mode | Network structure | Demand Pattern | Decision Variables | Solution Approach |
|---|---|---|---|---|---|---|
| Guler and Cassidy (2012) | | Bus & cars | General | - | - | Simulation |
| Geroliminis et al. (2014) | - | Bus & cars | Real network | Trapezoidal shape & real data | Passenger flow | Optimization, closed form, algorithm & simulation |
| Liu et al. (2015) | Min. of total system travel time | Bus & cars | - | Parameter | Number of vehicles in A cell, inflow and outflow of a cell, flow capacity, the ratio of the free-flow speed and backward propagation speed, & maximum number of vehicles held in a cell | Optimization software |
| Loder et al. (2019) | Min. total production | Bus & cars | - | - | Accumulation of cars, Accumulation of buses, Total travel production, Average speed in the network & Speed of cars and buses | Closed form solution, real cases & simulation |
| Dakic et al. (2020) | Min. walking distance or min. number of transfers | Bus & cars | Rectangular | Parameter (different for each corridor) | Stop spacing, line spacing, headway/ frequency & fractions of dedicated bus lines | Algorithm |
| Dakic et al. (2021) | - | Bus & cars | - | Uniform distribution & real data | Shortest path & type of bus stop | Algorithm |

### 4.12. Flexible and fixed transit services

The comparison between on-demand, flexible and fixed transit services and operations planning for each type of transit system has always been an interesting challenge for researchers (Winter et al., 2018). Li and Quadrifoglio (2010) used analytical models and simulation methods to develop an approach to assist the managers of the transportation systems to find which feeder transit service (flexible or fixed) is more fitting for the network. The proposed approach was able to determine if it was necessary to change from one service to another, and if "yes" when it should be done. To develop this approach, Li and Quadrifoglio (2010) considered a weighted total cost function including walking, waiting, and riding costs. As a key finding, they showed that demand will significantly affect the switching point between fixed and flexible transit services in the network. Using analytical methods, Kim and Schonfeld (2012) conducted a study to find the effectiveness of purely fixed-routes/schedules, purely flexible and variable types of public transit systems in serving different levels of travel demand in a city. They defined a variable transit type as a combination of fixed and flexible ones, where high demands are being served by fixed transit services and flexible systems serve passengers in the off-peak periods. They developed a new model for operations planning for the variable type of the service. However, for the pure modes, they used existing models in the literature and compared the results. Kim and Schonfeld (2014) investigated the possibility of reducing the total travel cost in the network by integrating fixed and flexible transit services in different transit zones. To this end, they considered a network including one terminal and multiple local zones and developed a stochastic optimization model to capture the variability in both travel time and waiting time effects on the performance of the network. The analytical models were applied to find the optimal values of the key variables including the number of the transit zones, service type for each transit zone, the vehicle, capacity, dispatching policy, fleet size, and slack times. A similar study has been done by Qiu et al. (2015), where they developed a model that helps owners of the system to make a decision about the type of public transit service in a transit zone. They assumed a varying demand in the network can be served by a fixed-route or a flex-route transit service, and presented a service quality function to evaluate the performance of the system. Using analytical approaches and simulation methods, they analyzed the network's operation under different scenarios to find optimal switching points between two transit systems.

Considering a transit zone that must be served with a many-to-one last mile transit service using fixed or flexible routes, public buses, Guo et al. (2017) and Guo et al. (2018) developed a stochastic dynamic analytical model to find the optimal switching points and to plan the operation of shared autonomous public transit fleets. They also tried to understand how the different patterns of demand would affect the switching point in the network. Kim et al. (2019) developed an analytical model to optimize transit zone sizes and dispatching policies for a flexible public bus transit system that serves many-to-one travel patterns. Using the proposed model, they tried to investigate the relationships between dispatching policies and optimal transit zone size. Different cost factors such as bus operating cost, waiting cost, in-vehicle cost, and supplier cost were considered in the modeling process. It was assumed that vehicles have capacity limitations. Mehran et al. (2020) defined the semi-flexible transit system as a system that is a combination of fixed-route and flexible-route bus transit services. They assumed the semi-flexible transit system was serving passengers with a flexible headway in a low demand area, on a fixed public bus route that includes just a few stops. Having these assumptions, Mehran et al. (2020) presented analytical models to calculate the operational costs of different service types, estimate demand patterns, and identify switching points. Kim and Roche (2021) proposed an analytical optimization method to obtain optimal headways and transit zones for a flexible public bus transit service that served areas with a low/mid passenger travel demand. They tried to find closed form solutions for the decision variables, aiming to minimize the total travel cost function. Compared to Kim et al. (2019), as a contribution, Kim and Roche (2021) considered a multi-period demand pattern in their modeling process.

Table 17 represents the main features of the reviewed papers in this subsection. As can be seen from this table, most of the reviewed papers in this section have considered rectangular network structures and a uniform demand pattern. However, some studies also considered time- and location-dependent demand density functions in their modeling process. Vehicle/fleet size, service area/zone, service frequency, and type of service are the most commonly considered decision variables in this sub-category. While some of these studies were able to find closed-form solutions, others have applied simulation techniques and numerical examples to find optimal solutions. Waiting costs, walking costs, in-vehicle costs, capital costs, and operation costs are mostly considered by the reviewed studies as the cost terms in the objective function.

**Table 17.** Main features of the reviewed papers in subcategory 12.

| Reference | Objective Function | Transit mode | Network structure | Demand Pattern | Decision Variables | Solution Approach |
|---|---|---|---|---|---|---|
| Li and Quadrifoglio (2010) | Weighted sum of expected waiting, walking & in-vehicle times | Bus | Rectangular | Poisson distribution for temporal distribution & uniform distribution for spatial distribution | Best feeder transit system | Simulation |
| Kim and Schonfeld (2012) | Min. sum of capital cost, operation cost, waiting, access & in-vehicle costs | Bus | Rectangular and a line-haul route | Uniform distribution | Vehicle size, service area, route spacing & headway/ frequency | Closed form solution |
| Kim and Schonfeld (2014) | Min. sum of operation cost, waiting, access, transfer & in-vehicle costs & min. sum of operation cost, waiting, transfer & in-vehicle costs | Bus | Rectangular and a line-haul route | Uniformly distributed over space: different for each region | Number of zones, vehicle size, fleet size, service area, route spacing, headway/ frequency & slack time | Closed form solution, numerical optimization & GA |
| Qiu et al. (2015) | Sum. of expected walking time, expected waiting time & expected waiting | Bus | Rectangular and a flex-route – real case | Uniform distribution | Best service type | Simulation |
| Guo et al. (2018)* | Min. sum of operation, waiting, access & in-vehicle costs | Bus | Rectangular and a line-haul route | Random time-dependent demand & uniform distribution (many-to-one /one-to-many) | Optimal switching point, route spacing, vehicle size, headway/ frequency & service area | Closed form solution |
| Kim et al. (2019) | Min. sum of operation, waiting & in-vehicle costs | Bus | Rectangular and a line-haul route | Different demand densities (many-to-one /one-to-many) | Zone size & headway/ frequency | Closed form solution |
| Mehran et al. (2020) | Sum of operation & fleet purchase | Bus | General | Estimated O-D matrix | Best service type | - |
| Kim and Roche (2021) | Min. sum of operation, waiting & in-vehicle costs | Bus | A single zone with a line-haul route | Uniform distribution (many-to-many) | Zone size & headway/ frequency | Closed form solution |

*: The paper has considered an emerging technology.

### 4.13. Sustainability

For this part, please see Mahmoudi et al. (2024).

## 5. Possible extensions and future research directions for the applications of analytical approaches in PBTNDP and operations planning

In this section, according to the reviewed literature, some possible extensions related to the applications of analytical methods in PBTNDP&OP have been introduced. In general, subjects that have not been yet analyzed or where there has been little research include the following (for additional discussion for this part see Mahmoudi et al. (2024) and Mahmoudi et al. (2025)):

- **Reliability.** While reliability is one of the most important factors in evaluating the performance of a PBTS from both users' and agency's perspectives, interestingly, the number of published papers that have applied analytical methods to study this subject is very limited. Also, most of the existing studies have only considered a single transit line to model the reliability of PBTSs, while public bus transit networks are much more than a single line, and any unexpected change in any location of the network can affect the performance of the whole network. Overall, researchers must put an effort in this area to model the reliability of the PBTSs considering different factors, such as crowding, unexpected recurrent traffic, accidents, any crisis in the network, etc., and different network and service structures.

  In addition to modeling the reliability of the networks, analytical models can be applied to analyze different possible strategies to enhance the reliability of the network. For example, the possibility of supporting the existing system during peak hours or when there is an unexpected crowding in the system with any fixed- or flexible-route service, shared transit mode, or pubic buses of various capacities, as well as the impacts of this strategy on service reliability, can be analyzed using analytical methods.

Using analytical approaches to design a reliable and coordinated timetable for the whole bus transit network while considering other transit modes in the network is another important problem that a limited number of studies have only focused on. Different aspects of this problem subject can be the main focus of future studies.

- **Feeder services.** So far, many studies have focused on the planning of feeder transit services where public bus service was acting as a feeder system for another mass transit service or was being served by a feeder service. However, especially with the introduction of emerging technologies in private and public transit modes, in addition to newly emerging problems, all previously studied problems must be reanalyzed in order to update the existing feeder-trunk systems. Analyzing the advantages and disadvantages of fixed- and flexible-route feeder transit services that are feeding a trunk service by regular buses, shared vehicles, shared automated vehicles, bikes, e-bikes, connected automated vehicles, e-hailing services, etc., using analytical methods will be a hot research direction in the near future. For example, an interesting problem related to feeder transit is fleet management and vehicle capacity. While users face high crowding costs during peak hours, the agency faces unused capacity costs during off-peak hours. One can suggest using mixed small and big buses to serve passengers in a feeder system. Another can suggest using one of the big or small sized buses but with different operation plans during off-peak and peak hours. However, another interesting scenario could be using modular buses, where the size of the bus can be adjusted based on the demand level. Providing an optimal operation plan for each scenario and comparing the performance of a feeder transit service under different scenarios can be an important research direction.
- **Network structure.** There are a large number of studies that have investigated the network design problem for public bus transit services. The main objective of these studies is to design or redesign the topology of the network, but most studies have considered many simplifying assumptions regarding the topology of the network. In many cases, the existing studies have considered an idealized radial or rectangular/square network. Idealized assumptions were mostly about the distributions of demand and origins/destinations, line spacings, etc. However, in real-world problems, developing models based on these assumptions can lead to unrealistic and inaccurate results. A few researchers have attempted to study the PBNDP using analytical methods without taking the simplifying assumptions into account; however, still, there is much to study. For example, except for a very small number of studies, there isn't any study that has considered a hybrid network structure. In a hybrid network structure, the topology of the network in the different zone areas will be different. For example, a part of the network could be rectangular and another part radial, and the PBTSs will serve passengers in both areas, which means, for instance, the first stop of a bus route will be in the rectangular part of the network and the last stop will be in the radial part. Using analytical methods for this type of problem will be a challenging but interesting research direction. Another suggestion is not to consider the fixed line spacings for the different network structures.

Another limitation of the existing studies in PBTNDP&OP using analytical approaches, is consideration of a specific density function or distribution function for the travel demand in the whole network. Most studies have considered a uniform distribution function, and a few studies have considered other functions. However, in the real world, the demand function in each service area and any time period will be different. For example, zone 1 can follow density function A, and zone 2 can follow density function B, between 8 a.m. and 9 a.m. But, between 9 a.m. and 10 a.m. the demand density function for zones 1 and 2 can be function C and function D, respectively. Hence, it would be recommended to consider time-dependent location-dependent and form varying demand density functions in the modeling process. However, whenever analytical approaches are not capable of managing the complexity of such problems, they could be integrated with other methodologies.

- **Uncertainty.** One of the main limitations in many of the identified studies is considering deterministic form for most of the key parameters in the modeling process. In the real world, for example, it appears that walking time, in-vehicle travel time, waiting time, operations costs, and maintenance costs, which are the main cost factors considered in the majority of the analytical studies, should not be in the deterministic form. Hence, developing analytical models based on the uncertainty of travel times and costs is another possible extension.

- **Simplifying assumptions and idealized structures**. In general, one of the main criticisms of most of the analytical studies is that they consider many simplifying assumptions in the modeling process. These studies mostly prefer to work on idealized structures rather than those close to real-world conditions. Although there are a good number of attempts to overcome this weakness, it is time to conduct more challenging studies by applying analytical models to non-idealized networks and under more realistic assumptions. For example, in most cases, researchers have assumed a uniform distribution of travel demand in the network; however, in the real world, it is hard to have this assumption in any public transit network.

  Another simplifying assumption considered by most of the well-cited studies in this area is about the type of services. Most of the proposed analytical models are developed based on single route, many-to-one and one-to-many transit services, while the real world conditions are different. Urban public transit networks are highly interconnected networks where the performance and condition of a part of the network will affect the performance of the other parts of the network. Therefore, to capture real world conditions instead of a many-to-one/one-to-many transit route, researchers must consider the whole transit network in modeling process. Although a number of studies have applied analytical methods to cover this limitation, many problems can still be defined and modeled in this area.
- **Hybrid methods**. Generally speaking, hybrid methods, when at least two methods have been integrated to investigate a problem, perform more efficiently in many cases than non-hybrid methods. This is because, in most cases, the goal of integrating methods is to compensate for a model's limitation(s) with the benefits of another model. Reviewing the literature shows that there is a very small number of studies that integrated analytical methods with other methods to study the PBNDP&OP, while there are so many advantages to doing it. For example, analytical methods can be integrated with mathematical programming and heuristic/metaheuristic algorithms to deal with the PBNDP&OP. In this case, a part of the problem can be solved by analytical methods, especially when it is necessary to have closed-form solutions for that part, and another part can be solved with mathematical programming. Since in many cases, using mathematical programming to model a problem related to the PBNDP&OP leads to a model with very high complexity, different algorithms will be used to find near to optimal solutions for the proposed model. However, the quality of the obtained solution and the execution time of these algorithms is highly dependent on the quality of the starting point or initial solution. As such possible research contribution can lay in modeling different problems related to the PBNDP&OP with both analytical and mathematical programming methodologies, then using the closed-form solution obtained by the analytical method as the initial solution for the proposed algorithms to solve the model developed using mathematical programming. Also, a specific problem can be solved with both methodologies, and results can be compared to find managerial insights.

Another idea is to integrate analytical methods with data-driven and empirical methods. In particular, when the objective is to propose a dynamic analytical model to study a specific problem related to the PBNDP&OP or a model to make real-time decisions (e.g. real-time headway setting), analytical methods can be combined with data-driven ones to find more realistic results. The results using these methods can also lead to more reliable transit services.

- **Demand patterns and network structures**. All of the studies related to PBTNDP&OP have assumed that there is a given uniform demand pattern or any other density function for demand, then they have applied an analytical model to design the network. A highly challenging but interesting question that has not yet been addressed in the literature is looking at PBTNDP&OP from a different perspective. It means a good research direction will involve applying analytical models to find how the structure of an existing system affects the demand pattern and forms the demand density function in the different areas of the network. This type of information can help decision makers to redesign PBTNs or to re-plan the operation of PBTSs.

## 6. Analytical methods versus mathematical programming

Analyzing the number of published studies in each category shows that although most of the earlier published works between 1968 and 2021 belong to the first category (analytical studies), after 2001 there is a significant increase in the number of studies that have applied mathematical modeling and optimization algorithm to study problems related to the PBTNDP&OPs. By the end of 2021, the total number of published papers in the second category is almost three times the number of publications in the first category. However, there is not necessarily any superiority of one method over the other. Each method has its own advantages and disadvantages. Here there is only a critical comparison between these two popular methodologies.

**Network structure: Idealization versus realism.** Analyzing the existing literature, it can be seen that there are significant differences between the structure of the networks analyzed by these two methodologies. In most of the analytical studies, there are many simplifying assumptions about the topology of the network, size of the network, and the demand pattern in the network, while published studies in the second group have considered more realistic assumptions. Analytical studies mostly have focused on idealized network structures such as a single linear transit line or corridor, rectangular grid network, idealized radial networks, etc. They have also mostly considered uniform demand patterns and small size networks. Although there are a few attempts in the analytical studies to consider more realistic assumptions, in general, it can be claimed that analytical approaches go well with idealized network structures and assumptions. On the other hand, looking at the literature it can be seen that most of the studies belonging to the second group have considered real transit networks with a big variety of demand density functions. Also, the number of studies that have applied multimodal transit network design problems using mathematical programming and optimization algorithm is way more than the number of published studies by analytical methods. While analytical studies mostly analyze only bus transit systems or bus-based feeder-trunk transit service, using mathematical programming and optimization algorithms, there is a good number of identified studies for the second category that have investigated the PBTNDP&OPs considering all possible transit modes in an urban area. Although there is not any mathematical model that can cover everything about a real-world problem accurately and produce 100% precise results, it seems that mathematical programming and optimization algorithms provide more near to real-world solutions compared to analytical studies.

**Modeling.** Analytical methods are powerful tools to develop cost functions, model the relationships between different variables, and convert geometric observations to mathematics. These methods work well with continuous variables and specific cost function forms. For these methods, the art is to identify and comprehensively understand the problem, and model it, trying to capture the real-world conditions as much as possible but at the same time keep the formulation solvable. However, when it comes to the optimization process, mathematical models and optimization algorithms can optimize the process considering a group of constraints with any decision variable type. In addition, mathematical models and optimization algorithms can be used to globally optimize or find near to optimal solutions for the problems with any cost function form. Also, these methods can be used for multiple-level and multi-objective optimization processes such as bi-level programming which is a common decision-making process in transportation studies.

**Computational complexity.** One of the main advantages of the analytical methods compared to the mathematical programming and optimization models is the complexity of the developed models using these methods. In most cases, the closed-form optimal solutions for the developed analytical models can be obtained by hand, using prepared coding packages or commercial software. Although using mathematical programming and optimization algorithms, the systems can be modeled based on more realistic assumptions and with subject to a large number of constraints, since PBTNDP&OPs for real-world transit networks are mostly highly complicated problems (e.g. NP-hard), it is impossible to find exact solutions or even near to optimal solutions for the understudied problem in the real-world networks. Because of this complexity, this method was not that populated method in the 20th century. However, by the increasing availability of advanced computers in the late 1990s and 21st century, most researchers started to use mathematical modeling and optimization algorithms in PBTNDP&OPs. Hence, for the second category, in addition to the modeling, the main art is to efficiently code the proposed model and the associated optimization algorithm(s). An efficient algorithm must be able to find near to optimal solutions with acceptable deviation from the optimal solution and in a reasonable execution time. Hence, while having access to advanced computers and having advanced coding abilities are essential requirements for most studies in the second category, in most cases, the quality of the analytical studies will not be affected by these factors.

**Results interpretation.** Another important advantage of the analytical approaches is being able to provide closed-form and interpretable solutions for most problems. But, as mentioned before, this is not the case for mathematical approaches and optimization algorithms. It should be noted that one of the important criteria of research is to obtain solutions that are understandable by both academia and managers/owners of the transit systems. Because of the high complexity of the developed mathematical modeling, they must be solved using optimization algorithms and computers, this is a big disadvantage for these methods for some reasons. First of all, most of the proposed algorithms are not able to find optimal solutions for the large case studies. Second, although there are some measures to evaluate the performance of the proposed models and algorithms, still one can call to question the reliability of the results, considering that the ability of the coding also affects significantly the performance of the algorithm. Since transportation systems are large and expensive systems, still there are significant differences between the performance of the network under optimal solutions and near to optimal solutions. On the other hand, while the ability of multi-objective optimization has been mentioned as an advantage for the mathematical modeling and optimization algorithms, ranking Pareto solutions is always a big challenge that leads to another decision and optimization process that adds extra complexity to the decision making process. Also, the obtained results by the mathematical modeling and optimization algorithms in most cases are not interpretable, because the codes and algorithms are acting as a black box and the outcome is just a numerical solution, not an analyzable closed-form formulation. For the same proposed model and algorithm, different computer programming can be written that provides different results for the same problem. Because of this feature, one can call into question the reliability of the obtained results by mathematical modeling.

**Conclusion.** In conclusion, none of the methodologies outperform another one. While analytical approaches mostly provide interpretable and insightful results with less optimization complexity, mathematical programming and optimization algorithms capture more real-world conditions in their studies and provide more flexibility in the modeling process. As mentioned before, using hybrid methodologies, where analytical approaches are integrated with mathematical programming, can provide more benefits in dealing with the PBTNDP&OP.

## 7. Conclusions

PTSs are always assumed to be the most efficient and sustainable transit mode. Especially, when it comes to daily intercity communications, public bus transit services (PBTSs) are assumed to be one of the popular transit services among city planners. As PBTS is one of the oldest public transit modes, and at the same time, one of the cheapest and easiest public transits to add to an urban transportation network, it has been turned to the backbone of the transportation network of many small, medium and large cities all around the world. PBTS, with its different versions (regular buses, BRT, feeder, shuttle, school bus, etc.), is applicable in any city with any transportation network size/structure and any population density. Because of this popularity and applicability, the public bus transit network design problem and operations planning (PBTNDP&OP) have attracted a lot of researchers' and transportation managers' attention. As a result, there is a good number of publications that have focused on different problems related to the PBTNDP&OP, using different methodologies. Analytical approaches, mathematical programming, and optimization algorithms are some of the most popular approaches that researchers have widely applied in this area.

In this paper, a comprehensive literature review has been conducted on the application of analytical approaches in PBTNDP&OPs. The paper includes four major sections: first, a precise statistical analysis has been done on the identified published papers that was the first of its kind. The objectives of this section were providing general information to both experts and junior researchers interested in PBTNDP&OP. Second, the existing literature related to the applications of the analytical approaches in PBTNDP&OPs was reviewed in 13 subcategories. The identified subcategories for this method are: 1) Feeder transit planning, 2) Dispatching policy, 3) Transit network design, 4) Fleet size and capacity of the vehicle, 5) Routing, 6) Rapid transit and shuttle transit system planning, 7) Travel time reliability, 8) Bus priority, 9) Bus holding strategies, 10) Emerging technologies in public bus transit systems, 11) Macroscopic Fundamental Diagram and traffic assignment, 12) Flexible and fixed transit services, 13) Sustainability. Forth, after comprehensively reviewing the selected papers, the identified future directions and possible extensions were reported in the paper. Finally, a critical discussion about the advantages and disadvantages of analytical approaches versus another alternative approach, i.e. mathematical programming, was presented in this paper.

The authors believe that provided taxonomy, review, and research directions in this paper can highlight different extensions and gaps for future studies and will inspire new research on the applications of analytical approaches in the PBTNDP&OPs.